\pdfoutput=1
\documentclass[12pt]{article}

\usepackage[T1]{fontenc}
\usepackage[utf8]{inputenc}
\usepackage[english]{babel}

\usepackage{amsmath}
\usepackage{amssymb}
\usepackage{amsthm}
\usepackage{mathtools}
\mathtoolsset{showonlyrefs,showmanualtags}

\allowdisplaybreaks

\usepackage[natbibapa,nodoi]{apacite}
\usepackage[affil-it]{authblk}
\usepackage{bbm}
\newcommand{\bbone}{\mathbbm{1}}
\usepackage{booktabs}
\usepackage[format=plain,font=small,labelfont=bf,textfont=up]{caption}
\usepackage{enumitem}
\setlist[enumerate]{itemsep=0.1ex}
\setlist[itemize]{itemsep=0.1ex}
\usepackage[letterpaper,margin=1.05in,bottom=1.65in]{geometry}
\usepackage{graphicx}
\usepackage[hidelinks]{hyperref}
\usepackage{setspace}
\usepackage{suppmake}

\usepackage[nohints,tight]{minitoc}

\newcommand\mainref\ref
\newcommand\suppref\ref

\DeclarePairedDelimiter\paren\lparen\rparen
\DeclarePairedDelimiter\bracket\lbrack\rbrack
\DeclarePairedDelimiter\braces\lbrace\rbrace

\DeclarePairedDelimiter\abs\lvert\rvert
\DeclarePairedDelimiter\floor\lfloor\rfloor

\providecommand{\bbone}{\mathbf{1}}
\DeclarePairedDelimiterXPP\indicator[1]{\bbone}{\lbrack}{\rbrack}{}{#1}

\DeclarePairedDelimiterXPP\expf[1]{\exp}{\lparen}{\rparen}{}{#1}
\DeclarePairedDelimiterXPP\logf[1]{\log}{\lparen}{\rparen}{}{#1}
\DeclarePairedDelimiterXPP\maxf[1]{\max}{\lparen}{\rparen}{}{#1}
\DeclarePairedDelimiterXPP\minf[1]{\min}{\lparen}{\rparen}{}{#1}

\DeclarePairedDelimiterXPP\func[2]{#1}{\lparen}{\rparen}{}{#2}

\newcommand{\set}[1]{#1}

\DeclarePairedDelimiter\setb\lbrace\rbrace

\newcommand{\Reals}{\mathbb{R}}

\newcommand{\Naturals}{\mathbb{N}}
\newcommand{\bset}{\setb{0, 1}}

\DeclareMathOperator*{\argmin}{arg\,min}

\DeclareMathOperator*{\esssup}{ess\,sup}

\renewcommand{\vec}[1]{\boldsymbol{#1}}

\makeatletter
\newcommand*{\tran}{{\mathpalette\@tran{}}}
\newcommand*{\@tran}[2]{\raisebox{\depth}{$\m@th#1\intercal$}}
\makeatother

\DeclarePairedDelimiter\norm\lVert\rVert
\DeclarePairedDelimiterXPP\tnorm[1]{}{\lVert}{\rVert_{1}}{}{#1}
\DeclarePairedDelimiterXPP\enorm[1]{}{\lVert}{\rVert_{2}}{}{#1}
\DeclarePairedDelimiterXPP\inorm[1]{}{\lVert}{\rVert_{\infty}}{}{#1}
\DeclarePairedDelimiterXPP\pnorm[2]{}{\lVert}{\rVert_{#1}}{}{#2}

\DeclarePairedDelimiterXPP\detf[1]{\det}{\lparen}{\rparen}{}{#1}

\DeclareMathOperator{\trsym}{tr}
\DeclarePairedDelimiterXPP\tr[1]{\trsym}{\lparen}{\rparen}{}{#1}

\DeclareMathOperator{\diagsym}{diag}
\DeclarePairedDelimiterXPP\diag[1]{\diagsym}{\lparen}{\rparen}{}{#1}

\DeclareMathOperator{\vectorizesym}{vec}
\DeclarePairedDelimiterXPP\vectorize[1]{\vectorizesym}{\lparen}{\rparen}{}{#1}

\providecommand\given{}
\newcommand\givensymbol[1]{\nonscript\:#1\vert\allowbreak\nonscript\:\mathopen{}}

\let\Pr\relax
\DeclareMathOperator{\Prsym}{pr}
\DeclarePairedDelimiterXPP\Pr[1]{\Prsym}{\lparen}{\rparen}{}{%
	\renewcommand\given{\givensymbol{\delimsize}}%
	#1}

\DeclarePairedDelimiterXPP\Prs[1]{\Prsym}{\lbrace}{\rbrace}{}{%
	\renewcommand\given{\givensymbol{\delimsize}}%
	#1}

\DeclarePairedDelimiterXPP\Prt[1]{\Prsym}{\lbrack}{\rbrack}{}{%
	\renewcommand\given{\givensymbol{\delimsize}}%
	#1}

\DeclareMathOperator{\Esym}{E}
\DeclarePairedDelimiterXPP\E[1]{\Esym}{\lparen}{\rparen}{}{%
	\renewcommand\given{\givensymbol{\delimsize}}%
	#1}

\DeclarePairedDelimiterXPP\Es[1]{\Esym}{\lbrace}{\rbrace}{}{%
	\renewcommand\given{\givensymbol{\delimsize}}%
	#1}

\DeclarePairedDelimiterXPP\Et[1]{\Esym}{\lbrack}{\rbrack}{}{%
	\renewcommand\given{\givensymbol{\delimsize}}%
	#1}

\DeclareMathOperator{\Varsym}{var}
\DeclarePairedDelimiterXPP\Var[1]{\Varsym}{\lparen}{\rparen}{}{%
	\renewcommand\given{\givensymbol{\delimsize}}%
	#1}

\DeclarePairedDelimiterXPP\Vars[1]{\Varsym}{\lbrace}{\rbrace}{}{%
	\renewcommand\given{\givensymbol{\delimsize}}%
	#1}

\DeclarePairedDelimiterXPP\EstVar[1]{\widehat{\Varsym}}{\lparen}{\rparen}{}{%
	\renewcommand\given{\givensymbol{\delimsize}}%
	#1}

\DeclareMathOperator{\Covsym}{cov}
\DeclarePairedDelimiterXPP\Cov[1]{\Covsym}{\lparen}{\rparen}{}{%
	\renewcommand\given{\givensymbol{\delimsize}}%
	#1}

\makeatletter
\newcommand{\indep}{\protect\mathpalette{\protect\@indep}{\perp}}
\newcommand*{\@indep}[2]{\mathrel{\rlap{$#1#2$}\mkern3mu{#1#2}}}
\makeatother

\newcommand{\darrow}{\overset{d}{\longrightarrow}}

\newcommand{\bigOsym}{\mathcal{O}}
\DeclarePairedDelimiterXPP\bigO[1]{\bigOsym}{\lparen}{\rparen}{}{#1}

\newcommand{\littleOsym}{o}
\DeclarePairedDelimiterXPP\littleO[1]{\littleOsym}{\lparen}{\rparen}{}{#1}

\newcommand{\bigOpsym}{\bigOsym_p}
\DeclarePairedDelimiterXPP\bigOp[1]{\bigOpsym}{\lparen}{\rparen}{}{#1}

\newcommand{\littleOpsym}{\littleOsym_p}
\DeclarePairedDelimiterXPP\littleOp[1]{\littleOpsym}{\lparen}{\rparen}{}{#1}

\newcommand{\bigOmegasym}{\Omega}
\DeclarePairedDelimiterXPP\bigOmega[1]{\bigOmegasym}{\lparen}{\rparen}{}{#1}

\newcommand{\littleOmegasym}{\omega}
\DeclarePairedDelimiterXPP\littleOmega[1]{\littleOmegasym}{\lparen}{\rparen}{}{#1}

\newcommand{\bigThetasym}{\Theta}
\DeclarePairedDelimiterXPP\bigTheta[1]{\bigThetasym}{\lparen}{\rparen}{}{#1}

\newcommand{\sumin}{\sum_{i=1}^n}
\newcommand{\sumjn}{\sum_{j=1}^n}
\newcommand{\avgin}{\frac{1}{n}\sum_{i=1}^n}

\newcommand{\qqquad}{\qquad\qquad}
\newcommand{\qqqquad}{\qqquad\qquad}
\newcommand{\qqqqquad}{\qqqquad\qquad}

\newcommand{\quadtext}[1]{\quad\text{#1}\quad}
\newcommand{\qquadtext}[1]{\qquad\text{#1}\qquad}

\newcommand{\quadand}{\quadtext{and}}
\newcommand{\qquadand}{\qquadtext{and}}

\theoremstyle{plain}

\newtheorem{corollary}{Corollary}
\newtheorem{lemma}{Lemma}
\newtheorem{proposition}{Proposition}

\newtheorem{innerrefproposition}{Proposition}
\newenvironment{refproposition}[1]
{\renewcommand\theinnerrefproposition{#1}\innerrefproposition}
{\endinnerrefproposition}

\theoremstyle{definition}

\newtheorem{condition}{Condition}
\newtheorem{definition}{Definition}

\theoremstyle{remark}

\renewcommand{\set}[1]{\mathcal{#1}}

\newcommand{\constarb}{k}
\newcommand{\constpo}{k_1}
\newcommand{\constpi}{k_2}
\newcommand{\constpred}{k_3}

\newcommand{\Sample}{\set{U}}

\newcommand{\zs}{z}
\newcommand{\z}{\vec{\zs}}
\newcommand{\ze}[1]{\zs_{#1}}
\newcommand{\zei}{\ze{i}}
\newcommand{\zej}{\ze{j}}

\newcommand{\Zs}{Z}
\newcommand{\Z}{\vec{\Zs}}
\newcommand{\Ze}[1]{\Zs_{#1}}
\newcommand{\Zei}{\Ze{i}}
\newcommand{\Zej}{\Ze{j}}

\newcommand{\zall}{\mathcal{Z}}

\newcommand{\posym}{y}
\newcommand{\pos}[1]{\posym_{#1}}
\newcommand{\posi}{\pos{i}}
\newcommand{\posj}{\pos{j}}

\DeclarePairedDelimiterXPP\po[2]{\pos{#1}}{\lparen}{\rparen}{}{#2}
\DeclarePairedDelimiterXPP\poi[1]{\posi}{\lparen}{\rparen}{}{#1}
\DeclarePairedDelimiterXPP\poj[1]{\posj}{\lparen}{\rparen}{}{#1}

\newcommand{\poiz}{\poi{\z}}

\newcommand{\oosym}{Y}
\newcommand{\oo}[1]{\oosym_{#1}}
\newcommand{\ooi}{\oo{i}}
\newcommand{\ooj}{\oo{j}}

\newcommand{\exmsym}{d}
\newcommand{\exms}[1]{\exmsym_{#1}}
\newcommand{\exmsi}{\exms{i}}
\DeclarePairedDelimiterXPP\exm[2]{\exms{#1}}{\lparen}{\rparen}{}{#2}
\DeclarePairedDelimiterXPP\exmi[1]{\exmsi}{\lparen}{\rparen}{}{#1}
\newcommand{\exmiz}{\exmi{\z}}

\newcommand{\exosym}{D}
\newcommand{\exo}[1]{\exosym_{#1}}
\newcommand{\exoi}{\exo{i}}
\newcommand{\exoj}{\exo{j}}

\newcommand{\exi}[2]{\exo{#1#2}}
\newcommand{\exii}[1]{\exi{i}{#1}}
\newcommand{\exij}[1]{\exi{j}{#1}}
\newcommand{\exiie}{\exii{\exe}}
\newcommand{\exiia}{\exii{\exea}}
\newcommand{\exiib}{\exii{\exeb}}
\newcommand{\exije}{\exij{\exe}}
\newcommand{\exija}{\exij{\exea}}
\newcommand{\exijb}{\exij{\exeb}}

\newcommand{\exe}{\exmsym}
\newcommand{\exeone}{\exe_{1}}
\newcommand{\exetwo}{\exe_{2}}
\newcommand{\exea}{a}
\newcommand{\exeb}{b}

\newcommand{\exall}{\Delta}

\newcommand{\presym}{\pi}
\newcommand{\pres}[1]{\presym_{#1}}

\DeclarePairedDelimiterXPP\pre[2]{\pres{#1}}{\lparen}{\rparen}{}{#2}
\DeclarePairedDelimiterXPP\prei[1]{\pres{i}}{\lparen}{\rparen}{}{#1}
\DeclarePairedDelimiterXPP\prej[1]{\pres{j}}{\lparen}{\rparen}{}{#1}
\newcommand{\preie}{\prei{\exe}}
\newcommand{\preia}{\prei{\exea}}
\newcommand{\preib}{\prei{\exeb}}
\newcommand{\preio}{\prei{\exeone}}

\newcommand{\preje}{\prej{\exe}}

\newcommand{\prejt}{\prej{\exetwo}}

\DeclarePairedDelimiterXPP\preij[1]{\pres{ij}}{\lparen}{\rparen}{}{#1}

\newcommand{\preijot}{\preij{\exeone, \exetwo}}

\newcommand{\pocsexs}[1]{\tilde{\posym}_{#1}}
\DeclarePairedDelimiterXPP\pocsex[2]{\pocsexs{#1}}{\lparen}{\rparen}{}{#2}
\DeclarePairedDelimiterXPP\pocsexi[1]{\pocsexs{i}}{\lparen}{\rparen}{}{#1}

\newcommand{\poexs}[1]{\bar{\posym}_{#1}}
\DeclarePairedDelimiterXPP\poex[2]{\poexs{#1}}{\lparen}{\rparen}{}{#2}
\DeclarePairedDelimiterXPP\poexi[1]{\poexs{i}}{\lparen}{\rparen}{}{#1}
\DeclarePairedDelimiterXPP\poexj[1]{\poexs{j}}{\lparen}{\rparen}{}{#1}

\newcommand{\poexie}{\poexi{\exe}}
\newcommand{\poexia}{\poexi{\exea}}
\newcommand{\poexib}{\poexi{\exeb}}
\newcommand{\poexio}{\poexi{\exeone}}
\newcommand{\poexit}{\poexi{\exetwo}}
\newcommand{\poexje}{\poexj{\exe}}
\newcommand{\poexja}{\poexj{\exea}}
\newcommand{\poexjb}{\poexj{\exeb}}
\newcommand{\poexjo}{\poexj{\exeone}}
\newcommand{\poexjt}{\poexj{\exetwo}}

\newcommand{\porexs}[1]{\bar{\posym}_{#1}}
\DeclarePairedDelimiterXPP\porex[2]{\porexs{#1}}{\lparen}{\rparen}{}{#2}
\DeclarePairedDelimiterXPP\porexij[1]{\porexs{ij}}{\lparen}{\rparen}{}{#1}
\DeclarePairedDelimiterXPP\porexji[1]{\porexs{ji}}{\lparen}{\rparen}{}{#1}

\newcommand{\efs}{\tau}
\DeclarePairedDelimiterXPP\ef[1]{\efs}{\lparen}{\rparen}{}{#1}
\newcommand{\efd}{\ef{\exea, \exeb}}

\newcommand{\efcss}{\tilde{\efs}}
\DeclarePairedDelimiterXPP\efcs[1]{\efcss}{\lparen}{\rparen}{}{#1}

\DeclarePairedDelimiterXPP\est[1]{\hat{\efs}}{\lparen}{\rparen}{}{#1}
\newcommand{\estd}{\est{\exea, \exeb}}

\newcommand{\sersym}{\varepsilon}
\newcommand{\ser}[1]{\sersym_{#1}}
\newcommand{\seri}{\ser{i}}
\newcommand{\serj}{\ser{j}}

\newcommand{\seresym}{e}
\newcommand{\sere}[1]{\seresym_{#1}}
\DeclarePairedDelimiterXPP\seref[2]{\sere{#1}}{\lparen}{\rparen}{}{#2}
\DeclarePairedDelimiterXPP\serefij[1]{\sere{ij}}{\lparen}{\rparen}{}{#1}
\DeclarePairedDelimiterXPP\serefji[1]{\sere{ji}}{\lparen}{\rparen}{}{#1}
\newcommand{\serefije}{\serefij{\exe, \exe}}

\newcommand{\serefijab}{\serefij{\exea, \exeb}}

\newcommand{\serefijot}{\serefij{\exeone, \exetwo}}

\newcommand{\serefjie}{\serefji{\exe, \exe}}

\newcommand{\serefjito}{\serefji{\exetwo, \exeone}}

\newcommand{\sererandsym}{E}
\newcommand{\sererand}[1]{\sererandsym_{#1}}
\newcommand{\sereij}{\sererand{ij}}
\newcommand{\sereji}{\sererand{ji}}

\newcommand{\serusym}{U}
\newcommand{\seru}[1]{\serusym_{#1}}
\newcommand{\seruij}{\seru{ij}}
\newcommand{\seruji}{\seru{ji}}

\newcommand{\ddepsym}{\bar{c}}
\newcommand{\ddep}[1]{\ddepsym_{#1}}
\newcommand{\ddepe}{\ddep{\exe}}
\newcommand{\ddepa}{\ddep{\exea}}
\newcommand{\ddepb}{\ddep{\exeb}}

\newcommand{\odepsym}{\bar{\varepsilon}}
\newcommand{\odep}[1]{\odepsym_{#1}}
\newcommand{\odepe}{\odep{\exe}}
\newcommand{\odepa}{\odep{\exea}}
\newcommand{\odepb}{\odep{\exeb}}

\newcommand{\edepsym}{\bar{e}}
\newcommand{\edep}[1]{\edepsym_{#1}}
\newcommand{\edepe}{\edep{\exe}}
\newcommand{\edepa}{\edep{\exea}}
\newcommand{\edepb}{\edep{\exeb}}

\newcommand{\udepsym}{\bar{u}}
\newcommand{\udep}[1]{\udepsym_{#1}}
\newcommand{\udepe}{\udep{\exe}}
\newcommand{\udepa}{\udep{\exea}}
\newcommand{\udepb}{\udep{\exeb}}

\newcommand{\tdepsym}{\bar{t}}
\newcommand{\tdep}[1]{\tdepsym_{#1}}
\newcommand{\tdepe}{\tdep{\exe}}
\newcommand{\tdepa}{\tdep{\exea}}
\newcommand{\tdepb}{\tdep{\exeb}}

\newcommand{\limdist}{Q}

\newcommand{\convseq}{R_n}

\newcommand{\zpresym}{s}
\newcommand{\zpres}[1]{\zpresym_{#1}}
\DeclarePairedDelimiterXPP\zpre[2]{\zpres{#1}}{\lparen}{\rparen}{}{#2}
\DeclarePairedDelimiterXPP\zprei[1]{\zpres{i}}{\lparen}{\rparen}{}{#1}
\DeclarePairedDelimiterXPP\zprej[1]{\zpres{j}}{\lparen}{\rparen}{}{#1}
\newcommand{\zpreie}{\zprei{\exe}}
\newcommand{\zpreia}{\zprei{\exea}}
\newcommand{\zpreib}{\zprei{\exeb}}
\newcommand{\zpreje}{\zprej{\exe}}

\DeclarePairedDelimiterXPP\zpreij[1]{\zpres{ij}}{\lparen}{\rparen}{}{#1}
\DeclarePairedDelimiterXPP\zpreji[1]{\zpres{ji}}{\lparen}{\rparen}{}{#1}

\newcommand{\zpreijaa}{\zpreij{\exea, \exea}}
\newcommand{\zpreijbb}{\zpreij{\exeb, \exeb}}
\newcommand{\zpreijab}{\zpreij{\exea, \exeb}}
\newcommand{\zpreijba}{\zpreij{\exeb, \exea}}

\newcommand{\zpreijot}{\zpreij{\exeone, \exetwo}}

\newcommand{\zprejito}{\zpreji{\exetwo, \exeone}}

\newcommand{\exppimom}{p}

\newcommand{\zpravg}[1]{\bar{\zpresym}_{#1}}
\newcommand{\zpravge}{\zpravg{\exe}}
\newcommand{\zpravga}{\zpravg{\exea}}
\newcommand{\zpravgb}{\zpravg{\exeb}}

\DeclarePairedDelimiterXPP\zprmom[1]{\Pi}{\lparen}{\rparen}{}{#1}
\newcommand{\zprmome}{\zprmom{\exe, \exppimom}}
\newcommand{\zprmoma}{\zprmom{\exea, \exppimom}}
\newcommand{\zprmomb}{\zprmom{\exeb, \exppimom}}

\DeclarePairedDelimiterXPP\ddepext[2]{\ddepsym_{#1}}{\lparen}{\rparen}{}{#2}
\newcommand{\ddepexte}{\ddepext{\exe}{q}}

\newcommand{\ddepextstde}{\ddepsym_{\exe\exppimom}}
\newcommand{\ddepextstda}{\ddepsym_{\exea\exppimom}}
\newcommand{\ddepextstdb}{\ddepsym_{\exeb\exppimom}}

\DeclarePairedDelimiterXPP\EstVarAS[1]{\widehat{\Varsym}_{\normalfont\textsc{as}}}{\lparen}{\rparen}{}{%
	\renewcommand\given{\givensymbol{\delimsize}}%
	#1}

\newcommand{\EstVarHTAS}{\EstVarAS[\big]{\estd}}

\DeclarePairedDelimiterXPP\varwei[1]{P_{ij}}{\lparen}{\rparen}{}{#1}
\newcommand{\varbiassym}{B}
\newcommand{\varbias}[1]{\varbiassym_{#1}}
\DeclarePairedDelimiterXPP\varbiasf[2]{\varbias{#1}}{\lparen}{\rparen}{}{#2}

\newcommand{\covsym}{\vec{x}}
\newcommand{\cov}[1]{\covsym_{#1}}
\newcommand{\covi}{\cov{i}}

\newcommand{\coefsym}{\boldsymbol{\beta}}
\DeclarePairedDelimiterXPP\coef[1]{\coefsym}{\lparen}{\rparen}{}{#1}
\newcommand{\coefe}{\coef{\exe}}

\DeclarePairedDelimiterXPP\coefest[1]{\hat{\coefsym}}{\lparen}{\rparen}{}{#1}
\newcommand{\coefeste}{\coefest{\exe}}
\newcommand{\coefesta}{\coefest{\exea}}
\newcommand{\coefestb}{\coefest{\exeb}}

\newcommand{\poexests}[1]{\hat{\posym}_{#1}}
\DeclarePairedDelimiterXPP\poexest[2]{\poexests{#1}}{\lparen}{\rparen}{}{#2}
\DeclarePairedDelimiterXPP\poexesti[1]{\poexests{i}}{\lparen}{\rparen}{}{#1}
\DeclarePairedDelimiterXPP\poexestj[1]{\poexests{j}}{\lparen}{\rparen}{}{#1}

\newcommand{\poexestie}{\poexesti{\exe}}
\newcommand{\poexestia}{\poexesti{\exea}}
\newcommand{\poexestib}{\poexesti{\exeb}}
\newcommand{\poexestje}{\poexestj{\exe}}

\newcommand{\pdepsym}{p}
\newcommand{\pdep}[1]{\pdepsym_{#1}}
\newcommand{\pdepe}{\pdep{\exe}}
\newcommand{\pdepa}{\pdep{\exea}}
\newcommand{\pdepb}{\pdep{\exeb}}

\DeclarePairedDelimiterXPP\estha[1]{\hat{\efs}_{\normalfont\textsc{h{\'a}}}}{\lparen}{\rparen}{}{#1}
\newcommand{\esthad}{\estha{\exea, \exeb}}

\DeclarePairedDelimiterXPP\estde[1]{\hat{\efs}_{\normalfont\textsc{de}}}{\lparen}{\rparen}{}{#1}
\newcommand{\estded}{\estde{\exea, \exeb}}

\DeclarePairedDelimiterXPP\estgr[1]{\hat{\efs}_{\normalfont\textsc{gr}}}{\lparen}{\rparen}{}{#1}
\newcommand{\estgrd}{\estgr{\exea, \exeb}}

\newcommand{\musym}{\mu}
\newcommand{\mutr}[1]{\musym_{#1}}
\newcommand{\mutre}{\mutr{\exe}}
\newcommand{\mutra}{\mutr{\exea}}
\newcommand{\mutrb}{\mutr{\exeb}}

\newcommand{\muest}[1]{\hat{\musym}_{#1}}
\newcommand{\mueste}{\muest{\exe}}
\newcommand{\muesta}{\muest{\exea}}
\newcommand{\muestb}{\muest{\exeb}}

\newcommand{\nest}[1]{\hat{n}_{#1}}
\newcommand{\neste}{\nest{\exe}}
\newcommand{\nesta}{\nest{\exea}}
\newcommand{\nestb}{\nest{\exeb}}

\title{\Large\textbf{Causal inference with misspecified exposure mappings: separating definitions and assumptions}}

\author{Fredrik Sävje}
\affil{\small Department of Political Science and Department of Statistics \& Data Science\\Yale University}

\date{\today}

\begin{document}

\makeatletter%
\begin{NoHyper}\gdef\@thefnmark{}\@footnotetext{\hspace{-1em}I thank P Aronow, Avi Feller, Chad Hazlett, Gareth Nellis, Cyrus Samii, Drew Stommes and Panos Toulis for helpful comments and discussions.}\end{NoHyper}%
\makeatother%

\maketitle

\bigskip
\begin{abstract}
\begin{singlespace}
\noindent
Exposure mappings facilitate investigations of complex causal effects when units interact in experiments.
Current methods require experimenters to use the same exposure mappings both to define the effect of interest and to impose assumptions on the interference structure.
However, the two roles rarely coincide in practice, and experimenters are forced to make the often questionable assumption that their exposures are correctly specified.
This paper argues that the two roles exposure mappings currently serve can, and typically should, be separated, so that exposures are used to define effects without necessarily assuming that they are capturing the complete causal structure in the experiment.
The paper shows that this approach is practically viable by providing conditions under which exposure effects can be precisely estimated when the exposures are misspecified.
Some important questions remain open.

\vspace{0.2in}

\end{singlespace}
\end{abstract}

\doparttoc
\faketableofcontents

\clearpage

\section{Introduction}

Experimenters use exposure mappings to investigate complex causal effects involving interference between units.
An exposure mapping is a terse representation of the nominal treatments assigned to the units under study.
The representation facilitates both definition and estimation of causally relevant exposure effects, provided that the exposures are correctly specified.
The exposures are correctly specified when they capture all causal information pertaining to the nominal treatments.
Recognizing that it is difficult to construct correctly specified exposures, this paper considers estimation of exposure effects when the exposures are misspecified.

The central insight is that exposure mappings conventionally fill two roles.
The first role is to capture aspects of the nominal treatments deemed relevant or interesting for the question at hand.
To serve this purpose, exposure mappings do not need to be correctly specified; they can successfully capture relevant aspects of the causality operating in a certain context, including aspects of the interference, without necessarily capture all causal information.
The second role is to encode assumptions about the causal structure in the experiment.
It is convenient when an exposure mapping fills both roles simultaneously, as this allows experimenters to use standard causal inference techniques also in the presence of interference.
However, this can typically only be achieved by making untenable assumptions.
The current paper suggests that experimenters should focus primarily on the first role, making the exposures as relevant and interpretable as possible.
However, doing so would mean that standard techniques and results no longer apply because the exposures would become misspecified.

The paper demonstrates that separating the two roles nevertheless is practically viable by showing that conventional estimators of exposure effects are consistent for a generalization of the exposure effects under misspecification given relatively mild conditions on the specification errors.
Like the assumption of correctly specified exposure mappings, these conditions are generally not testable.
Their advantage is instead that they are considerably more tenable than the prevailing assumptions.
Assuming that the exposures are correctly specified is equivalent to assuming that the specification errors are zero uniformly.
To achieve consistency, it is sufficient to assume that the errors are not strongly dependent.
Weak dependence allows for potentially grave misspecification as long as the units' exposures are not misspecified in the same way.
This condition does not require that the interference takes any particular structure, making the results widely applicable.
This comes at the cost of being somewhat opaque, but the condition is sufficiently approachable to allow experimenters to reason about it in practical situations, as illustrated by several examples throughout the paper.

\section{Illustration}\label{sec:illustration}

\newcommand{\Bogota}{Bogot\'a}

\citet*{Blattman2021Place} conduct an experiment in the city of \Bogota\ to investigate whether intensive policing reduces crime.
Among 1,919 streets identified to be crime hot spots, the authors randomly selected 756 streets to be patrolled by police twice as much as the other hot spots.
An important concern was displacement of crime.
Intensive policing in some streets could cause criminals to move to other streets that are less patrolled, in which case crime is reduced in the targeted streets only by increasing crime elsewhere.
The authors addressed this concern by estimating spillover effects.
Among all non-hot spot streets that were within 250 meters of one of the 1,919 hot spots, the authors compared streets for which at least one of the neighboring hot spots was assigned intensive policing against streets for which none of the neighboring hot spots were assigned such policing.
If there was a displacement effect, streets neighboring treated hot spots would be expected to experience more crime than other streets.

The authors effectively defined two exposures for the non-hot spot streets.
The first exposure is to have at least one neighboring hot spot that is treated, and the other exposure is to have only untreated neighboring hot spots.
Following the current literature, we would need to assume that these two exposures are correctly specified to estimate their effect.
This means that as long as we hold the exposure of a non-hot spot street fixed, it cannot at all matter for its crime level how we otherwise assign policing resources in the city.
For example, it cannot matter which hot spots within the 250 meter radius are treated.
A treated hot spot 10 meters away must have the same effect as a hot spot 250 meters away.
A treated hot spot where a local gang is known to gather must have the same effect as a treated hot spot that is a busy commercial street with many pickpockets.
Furthermore, it cannot matter whether streets outside the 250 meter radius are treated.
A treated hot spot 260 meters away must be the same as one 5,000 meters away, and the same as when no distant hot spots are treated.
The number of treated hot spot streets within the 250 meter radius must also not matter.
One treated neighboring hot spot must have the same effect as ten.
The list goes on.

The assumption that the binary 250 meter radius exposure is correctly specified is not reasonable.
A superficial solution is refine the exposures.
For example, one could consider a set of concentric circles with radii $1, 2, 3, \dotsc$ meters centered at each non-hot spot street with a corresponding set of exposures capturing whether a treated hot spot falls in between two of the circles.
In this case, a non-hot spot street would be assigned exposure 239 if there is a treated hot spot street 239 meters away.
This will typically be an incomplete solution, for several reasons.
First, many of the problems listed above still remain: the refined exposures still ignore potential interaction effects; they still ignore the identity of the streets; and it is doubtful whether bee-line distance is the relevant metric of closeness between streets.
To make the exposures truly correctly specified, we would need thousands, if not millions, of them.
Second, it is difficult to interpret and present the results from a study with more than a handful of exposures.
Experimenters can hardly include tables with thousands or millions of estimated exposure effects.
Third, an experiment with many exposures will often suffer from positivity violations.
That is, most units will have no probability of being assigned most of the exposures.
For example, a non-hot spot street that does not have a neighboring hot spot in the interval 100 to 200 meters cannot be assigned exposures 100 to 200.
This means that the outcome for each exposure can be estimated only for a small subset of the experimental units.
The subsets will differ between exposures, making it difficult, or impossible, to estimate exposure effects because comparisons between exposures will be confounded by differences in the composition of the subsets.
While it sometimes is reasonable to make some refinements to the exposures, such refinements will rarely ensure that the exposures are correctly specified.

The suggestion of the current paper is to shift the focus from creating exposures that are correctly specified to creating exposures that are relevant for the question at hand.
In the experiment in \Bogota, we should ask what exposures would be informative to the policymaker, rather than asking what exposures are correctly specified.
If we suspect that potential displacement effects are strongest locally, the binary 250 meter radius exposure is a reasonable way to capture and measure those effects.
In doing so, we do not necessarily mean to say that the displacement effects are all the same within that 250 meter radius, nor that there is no displacement effect beyond 250 meters.
We simply want to measure the average difference in crime level between when at least one neighboring hot spot is assigned to be policed intensively and when all neighboring hot spots are assigned ordinary policing.
While these exposures do not capture all of the interference dynamics related to policing and crime in \Bogota, they are relatively easy to interpret, also for laypeople, and they are informative for policy.
The current literature would have us believe that it is not possible to precisely estimate this type of effect because the exposures are undoubtedly misspecified.
The purpose of this paper is to show that it is possible.

\section{Related work}\label{sec:related-work}

The idea of exposure mappings has its origin in \citet{Halloran1995Causal}, who consider causal inference under interference and provide foundational definitions.
This initial work was extended by \citet{Sobel2006What} and \citet{Hudgens2008}, who describe effect estimators for exposures based on proportions of treated units in disjoint groups of units.
\citet{Manski2013Identification} and \citet{Aronow2017Estimating} recognized that the key methodological tool in the prior literature was a terse description of the full treatment vector.
They used this insight to generalize the approach beyond proportions of treated units in disjoint groups to arbitrary summaries of the nominal treatments, as formalized by exposure mappings.
The terminology of exposures and exposure mappings used in this paper comes from \citet{Aronow2017Estimating}, but the underlying idea is essentially identical to the concept of effective treatments in \citet{Manski2013Identification}.
This literature assumes that the exposures are correctly specified.
The assumption transforms the interference problem into a causal inference problem without interference at the level of the exposures, which can be solved using standard techniques.

The assumption of correctly specified exposures has a direct parallel to the conventional no-interference assumption.
The necessity of the no-interference assumption when estimating average treatment effects has recently been investigated by \citet*{Saevje2021Average}.
These authors show that a generalization of the average treatment effect can be precisely estimated even in presence of unmodeled interference.
The current paper connects these ideas with the ideas in \citet{Manski2013Identification} and \citet{Aronow2017Estimating}, extending the results to exposure mappings of arbitrary complexity.
Additionally, the current paper derives its results using quantitative interference measures, whereas \citet{Saevje2021Average} used a qualitative measure.
That is, the interference measure used here takes into account not only whether interference exists between units, but also the strength of that interference.
The differences between these two assumptions are discussed in detail in Section~\suppref{sec:conditions-sah} in the online supplement.

To the best of my knowledge, the general case of misspecified exposures has not previously been investigated, but restricted forms of misspecification have.
Partly building on insights of the current paper, \citet{Leung2022Causal} develops methods to estimate exposure effects under misspecification when the interference structure is approximately known.
Leung assumes that the strength of the interference between units decays in the geodesic distance between vertices representing the units in a known graph.
This assumption is considerably weaker than the conventional assumption that the exposures are correctly specified, but it is more restrictive than the setting considered in this paper.
However, \citet{Leung2022Causal} can leverage the additional structure imposed by the interference graph to answer to some of the questions that remain open in the setting considered here.
Section~\suppref{sec:leung-connection} in the online supplement describes the connection to \citet{Leung2022Causal} in more detail.

\citet{Egami2021Spillover} studies estimation of spillover effects in partially unobserved interference networks, which is a way to formalize certain forms of misspecification.
\citet{Li2022Random} consider estimation of direct and spillover effects when the interference graph is generated by a graphon.
They do not require knowledge of the graph or the graphon for some of their results, but such knowledge is required for the estimation of spillover effects.
\citet{Wager2021Experimenting} and \citet{Munro2021Treatment} consider experimentation when units interfere through an equilibrium mechanism, such as a market price.
Section~\suppref{sec:example-market} in the online supplement explains how this setting can be understood in the framework of the current paper.
\citet*{Wang2022Design} describe a new type of estimand that provides an alternative way to describe spillover effects, and they show how this effect can be precisely estimated under certain types of misspecification.

\section{Misspecified exposures}

\subsection{Preliminaries}

Consider a sample of $n$ units indexed by $\setb{1, \dotsc, n}$.
Each unit is assigned one of two possible treatments $\zei \in \bset$.
The assignments of all units are collected in $\z = \paren{\ze{1}, \dotsc, \ze{n}}$, and the set of all possible assignments is $\zall = \bset^n$.
A function $\posi \colon \zall \to \Reals$ gives the realized outcome for unit $i$ under a specific, potentially counterfactual, assignment.
That is, $\poiz$ is the response of unit $i$ when the treatments are assigned as $\z \in \zall$.
The elements of the image of the function are called \emph{potential outcomes}.
The potential outcomes are assumed to be well-defined throughout the paper, which requires that no hidden versions exist of the treatments in $\zall$.
The function $\posi$ may depend on the full vector $\z$, so no restrictions are made at this point on the interference between units.

The treatments are assigned at random.
Let $\Z$ be a random vector denoting the randomly selected treatments.
The distribution of $\Z$ is called the \emph{assignment mechanism} or \emph{experimental design}.
The support of $\Z$ may be smaller than $\zall$, so that some assignment vectors are never realized by the design.
The design is the sole source of randomness under consideration in this paper, and the sample of units and their potential outcomes are considered non-random and fixed.
The observed outcome $\ooi$ for unit $i$ is defined as the potential outcome corresponding to the randomly selected treatment vector: $\ooi = \poi{\Z}$.

The asymptotic regime used in the large sample investigation considers a sequence of fixed samples implicitly indexed by $n$.
All quantities pertaining to the sample will thus have their own sequences.
This type of regime has been used extensively in the literature on design-based sampling.
It has more recently seen use in the design-based causal inference literature \citep[see, e.g.,][]{Lin2013Agnostic}.

\subsection{Exposures}

The potential outcomes contain all causal information pertaining to the treatments, and any causal quantity can be expressed solely using them.
However, definitions of such causal quantities may be complex, and it is often difficult to formulate relevant and interesting effects when one is working with $2^n$ distinct potential outcomes for each unit.
Exposures and exposure mappings are used to make the definitions more interpretable and intuitive.
The idea is that different assignment vectors often share a similar interpretation.
For example, in the illustration in Section~\ref{sec:illustration}, having a treated hot-spot street at a distance of 50 meters might have the same causal interpretation as a treated street at 75 meters.
The purpose of an exposure mapping is to encode this type of similarity.
In particular, two assignment vectors are mapped to the same exposure if they are deemed similar with respect to the application at hand.
The exposures can therefore be seen as nothing more than labels on subsets of $\zall$ that have the same or similar interpretation.

To state this formally, consider a set of exposure labels indexed by $\exall \subset \Naturals$.
An exposure mapping is then a function $\exmsi \colon \zall \to \exall$ for each unit that maps from all possible assignments to the exposures.
The exposure of unit $i$ is $\exmiz$ when the treatment assignments are $\z$.
If $\exmiz = \exmi{\z'}$, then $\z$ has a similar interpretation as $\z'$ with respect to unit $i$.
The realized exposure is a random variable because the treatments are randomly assigned.
Let $\exoi = \exmi{\Z}$ denote the realized exposure for unit $i$, and let $\preie = \Pr{\exoi = \exe}$ describe its marginal distribution.

\subsection{Exposure effects}\label{sec:exposure-effects}

The current convention is to assume that exposure mappings are correctly specified.
The assumption states that $\poiz = \poi{\z'}$ whenever $\exmi{\z} = \exmi{\z'}$ for all units $i \in \setb{1, \dotsc, n}$ and assignments $\z, \z' \in \zall$.
This implies that a function $\pocsexs{i} \colon \exall \to \Reals$ exists for each unit such that $\pocsexi{\exmi{\z}} = \poiz$ for all $\z \in \zall$, meaning that the exposures are assumed to capture the complete causal structure in the experiment.
As noted in the introduction, this forces the exposures to serve two roles: to encode which assignment vectors have similar interpretation and to encode the (presumed) causal structure in the experiment.

If the exposures are indeed correctly specified, the full treatment vector provides no causal information in addition to what a unit's exposure already provides.
This implies that we can use $\pocsexi{\exe}$ defined on $\exall$ rather than the more cumbersome potential outcomes $\poiz$ defined on the full $\zall$ without loss of information.
The reduction in complexity can be considerable, because $\abs{\zall}$ grows exponentially in $n$ while $\abs{\exall}$ typically is fixed.
We can then define causal effects in the usual way, as contrasts between potential outcomes produced by the exposures.
For example, the average causal effect of exposure $\exea \in \exall$ relative to $\exeb \in \exall$ would be $n^{-1} \sumin \braces[\big]{\pocsexi{\exea} - \pocsexi{\exeb}}$.
This is the definition in \citet{Aronow2017Estimating}.
The interpretation of these effects tends to be straightforward, because the exposures are typically chosen to be easy to interpret.

It is not possible to construct a function $\pocsexs{i} : \exall \to \Reals$ such that $\pocsexi{\exmi{\z}} = \poiz$ when the exposures are misspecified.
One alternative is to extend $\exall$ with more labels until the exposure mappings become correctly specified, which might require that $\abs{\exall} \approx \abs{\zall} = 2^n$.
But, as noted in Section~\ref{sec:illustration}, this will generally not be a feasible or desirable solution.

A more productive alternative is to make the definition of the exposure effects robust to misspecification.
We can achieve this by creating analogues of the exposure-based potential outcomes that remain unambiguous even when the exposures are misspecified.
Let $\poexs{i} : \exall \to \Reals$ be a function such that $\poexie = \Es{\poi{\Z} \given \exoi = \exe}$, where the expectation is taken over the design.
The interpretation of $\poexs{i}$ is essentially the same as for $\pocsexs{i}$.
The function captures the expected potential outcome under each exposure for each unit, so $\poexi{\exe}$ is the potential outcome we expect to be realized when unit $i$ is assigned to exposure $\exe$ under the current design.
A definition of an exposure effect under misspecification is now immediate.

\begin{definition}\label{def:exposure-effect}
The expected exposure effect for exposures $\exea$ and $\exeb$ is
\begin{equation}
	\efd
	= \frac{1}{n} \sumin \braces[\big]{\poexia - \poexib}.
\end{equation}
\end{definition}

Effects building on this idea have previously been discussed in the literature.
The earliest examples were introduced by \citet{Hudgens2008}.
The authors derive their main results assuming that the implicit exposures are correctly specified, but they define their effects allowing for some misspecification.
They achieve this by marginalizing over all assignments that map to the same exposure, exactly as in Definition~\ref{def:exposure-effect}.
\citet{Aronow2017Estimating} derive the expectation of their estimator after relaxing their assumption that the exposures are correctly specified.
They show the expectation is a particular weighted average of the potential outcomes defined on the full treatment vector, and this average can be shown to coincide with Definition~\ref{def:exposure-effect}.
When the exposures are correctly specified, the expected exposure effect defined here coincides with the conventional exposure effect defined by \citet{Aronow2017Estimating}.

Because the expected potential outcomes $\poexie$ use the distribution of $\Z$ conditional on $\exoi = \exe$ in the marginalization, the expected exposure effect $\efd$ captures two aspects.
The obvious captured aspect is the difference in the outcomes produced by changing the assigned exposure from $\exea$ to $\exeb$.
The second aspect is more subtle.
When the exposures are misspecified, a unit's outcome could differ when other units' exposures or treatments are changed even holding its own exposure fixed.
Because of the conditioning on $\exoi = \exe$, the distribution of other units' treatments might be different in the contrast, and the expected exposure effect could capture this difference.
\citet{Saevje2021Average} discuss an alternative estimand that only captures the first aspect when the exposure mapping is $\exmiz = \zei$.
However, as discussed in Section~\suppref{sec:other-estimands} in the online supplement, the nature of many exposure mappings prevents similar effects to be defined more generally, because the definition would involve potential outcomes that are inherently unrealizable.
Similarly, when $\exmiz = \zei$, one can eliminate the second aspect by assigning treatments independently between units, but that approach is not typically viable here because most exposure mappings introduce strong dependencies between exposures even if the nominal treatments are independent.

Consider Definition~\ref{def:exposure-effect} in the context of the experiment in \Bogota\ discussed in Section~\ref{sec:illustration}.
If there are spillover effects of crime between neighboring streets, we would expect the outcomes under the first exposure to be different than under the second exposure.
However, because we allow for misspecification, there will be a distribution of outcomes for each of the exposures.
The expected exposure effect compares the centers of these two outcome distributions.
This comparison provides insights to the policymaker about whether displacement of crime is an important concern in \Bogota.
If the effect is found to be small, the policymaker might have reason to target crime fighting measures exclusively on hot spot streets.
A more targeted approach would likely be more effective in preventing crime in the hot spots themselves, and the small expected exposure effect indicates that we do not need to be concerned about crime displacement, at least not locally.
But if the effect is found to be positive and large, the policymaker might have reason to consider a less targeted approach that would be better at preventing some of the displacement.

The \Bogota\ experiment also provides a context in which we can understand the second, more subtle aspect captured by the expected exposure effect, as discussed above.
Consider a clustered design where \Bogota\ is divided into a grid of squares with  one kilometer sides, and treatment is assigned to either all or none of the hot spot streets within each square.
With this design, a street with no treated hot spots within a 250 meter radius tends to have fewer treated hot spots also within the annulus (i.e., ring) with inner and outer radii of 250 and 500 meters.
Therefore, the exposures here act as proxies for treatment assignments beyond 250 meters.
If we find a large effect using this design, we cannot conclude that it necessarily is the 250 meter radius that is causally relevant; the effect could have been zero using the same exposure mapping if we had used a different design, because the exposures may then not act as proxies in the same way.
This illustrates that experimenters must carefully consider how their exposures interact with their design when they interpret exposure effects.

\subsection{Specification errors}\label{sec:spec-errors}

Misspecification introduces specification errors.
The errors can be formalized as differences between the potential outcomes based on the full treatment vector, which is known to be correctly specified, and the potential outcomes based on the exposures, which may be misspecified.

\begin{definition}[Specification error]\label{def:specification-error}
$\seri = \poi{\Z} - \poexi{\exoi}$.
\end{definition}

The assumption that the exposures are correctly specified is the same as assuming that the specification errors are zero with probability one for all units.
This insight suggests a way to weaken the assumption.
Rather than assuming that the specification errors are zero, it may be sufficient to ensure that they are small, or perhaps only that they are sufficiently controlled in some other way.

Small specification errors are indeed sufficient to precisely estimate exposure effects, but such a condition is unnecessarily strong.
Instead, the critical aspect is the dependence of errors between units.
There are several ways to formalize this dependence.
The route explored here is to measure the dependence by the expectation of the product of two units' specification errors conditional on the event that the units are assigned the same exposure: $\E{\seri\serj \given \exoi = \exoj = \exe}$.
The expectation is defined to be zero if the event $\exoi = \exoj = \exe$ has probability zero.
The overall dependence is captured by the following quantity.

\begin{definition}\label{def:error-dependence}
The overall error dependence for exposure $\exe \in \exall$ is
\begin{equation}
	\odepe = \frac{1}{n^2} \sumin \sum_{j \neq i} \E{\seri\serj \given \exoi = \exoj = \exe}.
\end{equation}
\end{definition}

We can understand this quantity to capture two sources of dependence.
The first is the conditioning event itself, capturing the fact that knowledge about $j$'s exposure can provide information about $i$'s outcome in excess of the information provided by $i$'s exposure.
For illustration, consider again the experiment in \Bogota.
It is not reasonable to assume the binary exposures used here are correctly specified, so $\poi{\Z}$ will vary even when $\exoi$ is fixed, and $\seri$ will not be zero.
This means that other units' exposures could provide information about a unit's specification error.
If $i$ and $j$ are two non-hot spot streets with partially overlapping 250 meter radii, then knowing that $j$ has one or more neighboring treated hot spots gives us reason to suspect that $\seri > 0$ even if we already knew that $i$ had one or more neighboring treated hot spots.
This is because when both $i$ and $j$ have neighboring hot spots that are treated, we have reason to believe that more than one of unit $i$'s neighboring hot spots are treated, indicating a greater local police presence that potentially could displace crime.

The second source is dependence in excess of what can be explained by the conditioning event.
This captures the fact that two units’ errors can be dependent if misspecified in the same way even if the exposures themselves provide no information about the outcomes.
This might be information about the assignment vector $\Z$ that is altogether lost by the exposure mapping, or information that can only be captured by intricate combinations of many exposures.
An example of this is general equilibrium effects.
Returning once more to the experiment in Section~\ref{sec:illustration}, we can imagine two or more equilibria for crime in \Bogota.
The first could be an overall low-crime equilibrium achieved when sufficiently many hot spots are policed intensively.
This could perhaps be because crime becomes too costly, so criminals find other ways to make a living or move to other cities.
The second could be a high-crime equilibrium achieved when few hot spots are policed.
The exposures might be correctly specified within each equilibrium, so that learning other units' exposures provide no more information about a unit's potential outcome as long as we remain in the same equilibrium.
However, if the experimental design induces variation in which equilibrium is realized, the specification errors will typically be highly dependent between units.
The equilibrium will act as a coordinating force for the specification errors, making $\odepe$ large.
This can be remedied either by ensuring that there are no such global equilibrium effects, or by picking an experimental design that does not induce such variation, so only one equilibrium is realized with high probability.
The numerical example in Section~\ref{sec:num-example} considers a setting with this type of equilibrium phenomenon.

An extended discussion about the specification errors is provided in Section~\suppref{sec:specification-errors} in the online supplement.
This includes a decomposition of the specification error that formalizes the two sources of dependence discussed above, an investigation of how the specification errors relate to the interference conditions used by \citet{Saevje2021Average} and \citet{Leung2022Causal}, and simulation studies that elucidate what the specification errors are and how they relate to estimation of exposure effects in concrete settings.

The generality of the definitions in this section can sometimes make them difficult to reason about.
One of the goals of the paper is to highlight the wide applicability of the idea that one can separate the two roles that exposure mappings traditionally have served, and the definitions were made to be general to emphasize this.
If one were to find the definitions too opaque to be helpful in practice, one could see them as templates for constructing more concrete and situational definitions in particular contexts.
Section~\suppref{sec:error-examples} in the online supplement contains several illustrations of such concrete constructions.

\section{The precision of conventional estimators under misspecification}\label{sec:precision-estimators}

\subsection{Estimator}

Commonly used estimators for exposure effects build on ideas originally introduced in the survey sampling literature.
The focus here is the Horvitz--Thompson estimator \citep{Horvitz1952}:
\begin{equation}
	\estd = \avgin \frac{\exiia \ooi}{\preia} - \avgin \frac{\exiib \ooi}{\preib},
\end{equation}
where $\exiie = \indicator{\exoi = \exe}$ is an indicator denoting whether unit $i$ is assigned exposure $\exe \in \exall$.

Proposition~8.1 in \citet{Aronow2017Estimating} implies that the Horvitz--Thompson estimator is unbiased for the expected exposure effect $\efd$.
This is true no matter how severe the misspecification might be.
See also Lemma~\suppref{lem:unbiasedness} in the online supplement.

The results in this paper can be extended to several other common estimators, including the difference-in-means and H\'ajek estimators, and estimators making covariate adjustments.
The investigation of these other estimators is relegated to Section~\suppref{sec:other-estimators} in the online supplement in the interest of space.

\subsection{Regularity conditions}

Three regularity conditions, which are not directly related to misspecification, facilitate the following investigation.
The first of these conditions is an assumption that the potential outcomes are bounded.
This assumption can be weakened without materially changing the results, but it eases the exposition.

\begin{condition}\label{cond:bounded-pos}
For all $i \in \setb{1, \dotsc, n}$ and $\z \in \zall$, $\abs{\poi{\z}} \leq \constpo < \infty$.
\end{condition}

The second and third conditions concern the design and the distribution it induces on the exposures.
The second condition is a positivity assumption.
This assumption can be weakened as well, as explored in Section~\suppref{sec:lack-positivity} in the online supplement.

\begin{condition}\label{cond:positivity}
An exposure $\exe \in \exall$ satisfies \emph{positivity} if $1/\preie \leq \constpi < \infty$ for all $i \in \setb{1, \dotsc, n}$.
\end{condition}

The third condition concerns the dependence between exposures.
The exposure mappings could potentially introduce strong dependencies between the realized exposures even if the design is well-behaved for the nominal treatments.
For example, if the exposures were defined as the proportion of treated units in the whole sample, it would follow that $\exo{1} = \exo{2} = \dotsb = \exo{n}$, which would make precise estimation impossible.
Note that the concern here is not dependencies in the outcomes, which would relate to the actual interference structure in the experiment, but dependencies between the exposure labels selected by the experimenter.
The following condition limits the dependence between exposures introduced by the design.

\begin{condition}\label{cond:limited-design-dependence}
An exposure $\exe \in \exall$ satisfies \emph{limited design dependence} if $\ddepe = \littleO{1}$, where
\begin{equation}
	\ddepe = \frac{1}{n^2} \sumin \sum_{i \neq j} \abs[\big]{\Cov[\big]{\exiie, \exije}}.
\end{equation}
\end{condition}

Condition~\ref{cond:limited-design-dependence} might appear less familiar than Conditions~\ref{cond:bounded-pos}~and~\ref{cond:positivity}, but it is a somewhat common assumption in the recent literature.
Indeed, it corresponds to Condition~4 in \citet{Aronow2017Estimating}.
The current condition is slightly weaker, as it considers quantitative dependence between exposures rather than the binary dependence concept used by \citet{Aronow2017Estimating}, but the interpretation carries over unchanged.

Note that the experimental design and the exposure mappings are known.
It is therefore possible to verify the positivity condition and to measure the amount of design dependence in a particular sample.
Indeed, design conditions like these are often seen as innocuous, because experimenters control their designs and can ensure that they hold.
This may not be the case when estimating exposure effects.
Exposure mappings tend to be complex, and it is often difficult to construct a design that would induce a desired distribution over the exposures.
In particular, positivity as stipulated by Condition~\ref{cond:positivity} will sometimes be difficult to achieve.
This is the motivation for the relaxation of the positivity assumption explored in Section~\suppref{sec:lack-positivity} in the online supplement.

\subsection{Variance bound and consistency}\label{sec:var-bound-consistency}

We now have the components needed to characterize the precision of the estimator.
The proof of the following proposition is presented in the online supplement.

\newcommand{\propvariancebound}{%
Provided that Conditions~\mainref{cond:bounded-pos}~and~\mainref{cond:positivity} hold for exposures $\exea$ and $\exeb$, the variance of the Horvitz--Thompson estimator is upper bounded by
\begin{equation}
	\Var[\big]{\estd} \leq 8 \constpo^2 \constpi / n + 20 \constpo^2 \constpi^2 \paren{\ddepa + \ddepb} + 4 \paren{\odepa + \odepb}.
\end{equation}
}

\begin{proposition}\label{prop:variance-bound}
\propvariancebound
\end{proposition}

The proposition demonstrates that the variance of the estimator is governed by three aspects, as captured by the three terms in the bound.
The first term captures variability induced by the fact that the exposures are randomly assigned.
That is, even when the exposures are correctly specified and independent, the estimator would still vary because different potential outcomes would be observed for different assignments.
The second term captures variability induced by dependence between exposures.
That is, even when the exposures are correctly specified, the estimator tends to be less precise when exposures are highly dependent.

The final term of the bound captures variance stemming from misspecification.
Recall that Definition~\ref{def:error-dependence} captures the dependence between the specification errors.
If the specification errors are strongly positively correlated, the estimator will exhibit more variability.
The bound makes clear that the magnitude of the specification errors in itself is less of a concern.
Large specification errors will affect the precision, but their effect is absorbed by the first term.
The intuition for this is that large but uncorrelated specification errors tend to cancel when averaged.

In order for the estimator to become more precise as the sample increases in size, the overall dependence between the specification errors must be limited.
This is captured by the following condition and corollary.

\begin{condition}\label{cond:limited-error-dependence}
An exposure $\exe \in \exall$ satisfies \emph{limited specification error dependence} if there exists some positive sequence $r_n = \littleO{1}$ such that $\odepe \leq r_n$.
\end{condition}

\newcommand{\coroconsistency}{%
Provided that Conditions~\mainref{cond:bounded-pos},~\mainref{cond:positivity},~\mainref{cond:limited-design-dependence}~and~\mainref{cond:limited-error-dependence} hold for exposures $\exea$ and $\exeb$, the Horvitz--Thompson estimator $\estd$ is consistent for the expected exposure effect $\efd$.%
}

\begin{corollary}\label{coro:consistency}
\coroconsistency
\end{corollary}

Because the estimator is unbiased, the variance bound directly describes the estimator's asymptotic behavior in a mean square sense.
Control over the terms in the bound thus provides consistency through Chebyshev's inequality, resulting in Corollary~\ref{coro:consistency}.

Condition~\ref{cond:limited-error-dependence} states that the overall dependence between specification errors becomes smaller as the sample grows.
One situation in which this is satisfied is when the interference is local, in the sense that $\E{\seri\serj \given \exoi = \exoj = \exe}$ is zero for most pairs of units.
This is essentially the condition used by \citet{Saevje2021Average}.
However, Condition~\ref{cond:limited-error-dependence} allows for global interference as long as it is limited.
That is, $\E{\seri\serj \given \exoi = \exoj = \exe}$ can be non-zero for all pairs of units as long as it is not large and positive on average.
Section~\ref{sec:num-example} below and Section~\suppref{sec:error-examples} in the online supplement provide concrete examples of such limited global interference.

Consider Condition~\ref{cond:limited-error-dependence} in the context of the experiment in \Bogota\ described in Section~\ref{sec:illustration}.
As discussed in Section~\ref{sec:spec-errors}, one source of the error dependence is the conditioning event $\exoi = \exoj = \exe$.
This is unlikely to be an issue in this setting, because only a small fraction of the streets will have overlapping 250 meter radii.
The fraction of overlapping streets will also diminish as the sample grows, presuming that the 250 meter radius is held fixed and the streets do not become more dense.
While exposures of more distant streets could occasionally be informative, that will not be the case for most pairs of streets.

The second source of error dependence could potentially be more problematic in \Bogota.
If equilibrium effects of the type discussed in Section~\ref{sec:spec-errors} exist, and design induces variation in which equilibrium is realized, then Condition~\ref{cond:limited-error-dependence} is unlikely to hold.
For example, we could imagine that intensive policing in a certain combination of hot spots could lead to a prominent figure in one of the city's criminal gangs being apprehended, which in turn could cause an overall low-crime environment that otherwise would not have happened.
However, global interference does not need to be a concern, and we can allow for equilibrium effects as long as the same equilibrium is realized with high probability.
In fact, we can allow for unstable equilibria if the units' outcomes are not affected in the same way, meaning that we avoid strong positive dependencies.
The concern with unstable equilibria is that the outcomes for most units will typically be affected in the same way when the equilibrium changes.

Note that positively correlated errors are the concern here.
If the errors are negatively correlated, they act to stabilize the estimator, making it more precise than under correctly specified exposures.
While this is theoretically possible, overall negatively correlated errors do not appear to be practically relevant.


\section{Numerical example}\label{sec:num-example}

As a concrete illustration of the ideas in this paper, consider an example study with global interference, meaning that all units can potentially interfere with all other units.
Additional numerical examples and details about this example appear in Section~\suppref{sec:error-examples} in the online supplement.

Suppose the experimenter is interested in evaluating the effect of a vaccine against a communicable disease, meaning that interference is expected to be present.
The experimenter has access to an observed network through which units are hypothesized to interact, and this network is used to define the exposures.
Denote an edge from unit $i$ to unit $j$ in the network as $g_{ij} = g_{ji} = 1$, and otherwise $g_{ij} = 0$.
For simplicity, the graph is a cycle graph, meaning that $g_{ij} = 1$ if and only if $\abs{i - j} \in \braces{1, n - 1}$.
A conventional Bernoulli design is used, meaning that treatments are assigned independently with equal probability.

There are two exposures of interest.
The first exposure, denoted $\exoi = 1$, is when the unit itself is untreated and both neighbors in the cycle graph are treated.
The second exposure, denoted $\exoi = 0$, is when the unit itself as well as both neighbors are untreated.
Because of the Bernoulli design and cycle graph structure, the exposure probabilities are $12.5\%$ for both exposures and for all units.
The treatment effect of interest is $\ef{1,0}$, which is the expected indirect effect of having two treated neighbors when being untreated compared to having no treated neighbors.

The potential outcomes are such that units interfere locally through the observed network, but they also interfere globally through behavior akin to herd immunity.
The fact that herd immunity could occur is not known or hypothesized by the experimenter.
This means that the exposures are misspecified.
If the number of vaccinated units passes some threshold, there is no transmission of the disease in the community under study, and the outcome (e.g., viral load) is zero for all units.
In more detail, the potential outcome functions are
\begin{equation}
	\poiz = \begin{cases}
		0 & \text{if } \sum_{j=1}^n \zej \geq \phi n,
		\\
		\alpha_i - \beta_i \zei - \gamma_i \sum_{j=1}^n g_{ij} \zej & \text{else},
	\end{cases}
\end{equation}
where $\alpha_i$ is drawn uniformly at random from $[15, 25]$, $\beta_i$ is drawn from $[1, 10]$ and $\gamma_i$ is drawn from $[1, 2]$.
The coefficients are drawn once and held fixed over the simulation rounds, mirroring the fact that all randomness under consideration comes from treatment assignment.
The parameter $\phi \in (0, 1)$ sets the herd immunity threshold; if a share of $\phi$ units are treated, then herd immunity occurs, and the outcome is zero for all units.

Three versions of this data generating process will be considered, corresponding to three values of the herd immunity threshold $\phi$.
The first version, labelled ``Rarely,'' uses $\phi = 0.52$, meaning that herd immunity occurs when at least $52\%$ of the units are treated.
The second version, labelled ``Infrequently,'' uses $\phi = 0.51$, and the third version, labelled ``Frequently,'' uses $\phi = 0.5$.
Note that $\sum_{j=1}^n \Zej$ follows a binomial distribution of $n$ trials with $0.5$ success probability.
Herd immunity will therefore occur somewhat rarely under the first version of the data generating process ($\phi = 0.52$), and it grows increasingly rare as the sample size grows.
However, when $\phi = 0.5$, herd immunity will be frequent also in large samples, because the distribution of $\sum_{j=1}^n \Zej$ is centered around $n / 2$.
Indeed, herd immunity occurs with roughly $50\%$ probability under all three versions of the process when $n = 100$.
But when $n = 10,000$, herd immunity occurs with a probability of only $0.003\%$ under the first version ($\phi = 0.52$), a probability of $2.33\%$ under the second version ($\phi = 0.51$), and a probability of $50.4\%$ under the third version ($\phi = 0.50$).

Interference is global under all three versions, because there are always situations in which changing a single unit's treatment assignment changes the outcome of essentially all other units.
That is, the event $\sum_{i=1}^n \Zei < \phi n \leq 1 + \sum_{i=1}^n \Zei$ has a non-zero probability of occurring under all three versions, and changing any untreated unit to be treated in this case induces herd immunity.
However, such global interference takes place only very rarely under the first version when the sample is large, because changing a single unit's treatment assignment under most assignments will not have global effects.
In this sense, global interference does exist under all three versions, but it is not practically relevant for the first two versions of the data generating process, $\phi \in \braces{0.51, 0.52}$, provided that the sample is sufficiently large.
This is not the case when $\phi = 0.5$, because there will be non-negligible variation in whether herd immunity occurs no matter how large the sample gets.

\begin{figure}
	\centering
	\includegraphics[width=\textwidth]{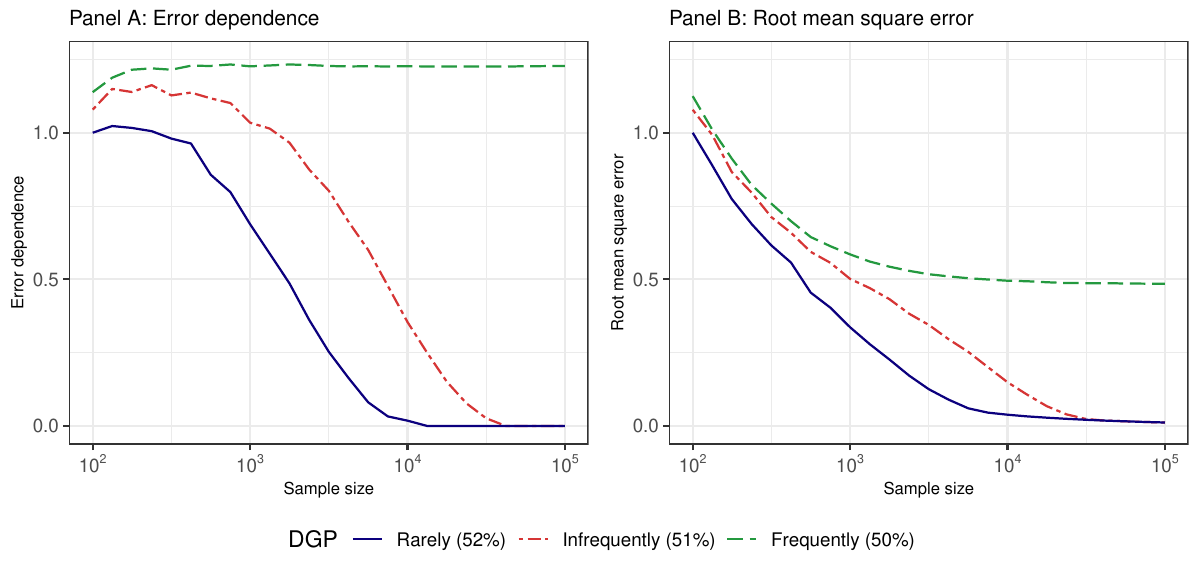}
	\caption{Overall error dependence (A) and root mean square error (B) under interference through herd immunity.}\label{fig:herd-interference}
\end{figure}

The results from the numerical example are presented in Figure~\ref{fig:herd-interference}.
Panel A presents the root of the sum of the magnitude of the overall error dependence for the two exposures, $\sqrt{\abs{\odep{1}} + \abs{\odep{0}}}$, as described in Definition~\ref{def:error-dependence}.
As the sample grows in size, the error dependence decreases toward zero for both ``Rarely'' ($\phi = 0.52$) and ``Infrequently'' ($\phi = 0.51$).
This mirrors the fact that variation in whether herd immunity occurs decreases as the sample size grows under both of these versions.
However, the error dependence does not decrease for the version ``Frequently,'' instead flattening out slightly below the 1.25 mark.
Therefore, Condition~\mainref{cond:limited-error-dependence} does not hold under the third version, and the consistency result of the current paper does not apply.

Panel B presents the root mean square errors of the H{\'a}jek estimator, which in this case coincides with the difference-in-means estimator.
The H{\'a}jek estimator is presented here because it is typically more precise than the Horvitz--Thompson estimator.
Root mean square errors of the Horvitz--Thompson estimator are reported in the online supplement.
The results are qualitatively the same, but differences between the data generating processes are less pronounced for the Horvitz--Thompson estimator because of the overall lower level of precision.

The figure shows that the precision of the H{\'a}jek estimator initially improves for all three versions of the data generating process.
The precision continues to improve for the first two versions (``Rarely'' and ``Infrequently''), and the mean square error is essentially zero when $n=100,000$.
However, under the third version (``Frequently''), the mean square error flattens out at around the $0.5$ mark, indicating that the estimator is not consistent.
This mirrors what we learned about the behavior of the error dependence measure in Panel A.
The reason the precision of the estimator initially improves also under the third version, despite no decrease in the error dependence, is that the first two terms of the variance bound in Proposition~\mainref{prop:variance-bound} approach zero as the sample size grows even if the error dependence is large.


\section{Discussion}

The motivating idea of this paper is that exposure mappings are primarily used to collect and describe sets of assignment vectors that share a similar interpretation with respect to the application at hand.
It is rare that exposures that serve this purpose are also correctly specified in the sense that they provide a complete description of the causal structure.
The results herein show that conventional point estimators perform well even if the exposures are misspecified, provided that the specification errors are only weakly dependent.
This gives reassurance to experimenters studying complex causal effects under interference, including spillover effects, that their analyses remain informative even in the event their exposures omit some aspects of the causal structure.
An expected exposure effect depends on the implemented design, which can make its interpretation challenging.
However, an expected exposure effect will typically be more interpretable than the alternative: an exposure effect based on extremely refined exposures.
These insights should prompt experimenters to focus on defining exposure mappings that are relevant for the question at hand.
If they instead followed the conventional recommendation and focused on making the exposures correctly specified, the exposures would often be too granular to be relevant and useful.

The investigation highlights several open questions.
The first question is whether estimators can be constructed to fully separate the two roles traditionally served by exposure mappings.
That is, estimators that can incorporate knowledge about the interference structure in the estimation of an exposure effect without necessarily changing the effect that is being estimated.
The results in this paper show that estimation is possible without such knowledge, but precision could perhaps be improved if we can take advantage of all available information about the interference, even if that information is imperfect or incomplete.
A step in this direction is the effect definition based on linear functionals described by \citet{Harshaw2022Design}.

The second open question concerns the limiting distribution of the point estimator under misspecification.
Section~\suppref{sec:limiting-distribution} in the online supplement describes two approaches to characterize the limiting distribution in this setting, but both approaches require considerably stronger assumptions than the limited error dependence condition used for the convergence results in Section~\ref{sec:precision-estimators}.
It remains an open question whether the limiting distribution can be characterized under conditions resembling limited error dependence.
A relevant result here is the central limit theorem by \citet{Kojevnikov2021Limit} for network dependent observations under the assumption that the dependence decays in geodesic distance, which was used by \citet{Leung2022Causal} to investigate the limiting distribution of estimators of expected exposure effects under approximate neighborhood interference.
While this approach requires stronger assumptions than limited error dependence and requires that the interference structure is known approximately, it is undoubtedly an important step in the right direction.

A final set of open questions concerns variance estimation.
Section~\suppref{sec:variance-estimation} in the online supplement notes that variance estimation under both correctly specified and misspecified exposures is difficult, because complex exposure mappings can make many joint exposure probabilities small or zero.
Experimenters can address this problem by making the estimator conservative, but the estimator will often be excessively conservative.
Recent work by \citet{Harshaw2021Optimized} describes methods that can mitigate the conservativeness when the exposures are correctly specified.
However, variance estimation appears to be particularly challenging under misspecification, because the units' specification errors can be dependent in such a way to make variance estimators anti-conservative.
Also in this case do \citet{Kojevnikov2021Limit} and \citet{Leung2022Causal} provide an encouraging result by describing a variance estimator to be used under misspecification when the interference structure is known approximately.
Relatedly, the online supplement describes ways to leverage unstructured partial knowledge of the interference in an experiment to partially mitigate the problem.
It nevertheless remains an open question whether practically useful variance estimators exist when little is known about the interference structure.

\begin{singlespace}
\bibliography{references}
\end{singlespace}

\clearpage
\startsupp{S}{Additional results and proofs}{supp:proofs}

\section{Understanding the specification errors}\label{sec:specification-errors}

This section provides additional results, explanations and illustrations to aid the understanding of the specification errors and the conditions imposed on them.
Section~\ref{sec:decomposition} provides a decomposition of the specification error into two parts to elucidate when error dependence arises and when it is controlled.
Section~\ref{sec:conditions-sah} relates the condition on the error dependence in the current paper to the interference conditions in \citet{Saevje2021Average}.
Section~\ref{sec:leung-connection} does the same for the interference conditions in \citet{Leung2022Causal}.
Section~\ref{sec:error-examples} provides several concrete numerical examples that illustrate how the specification errors behave in specific settings.

\subsection{Decomposition}\label{sec:decomposition}

The main paper defines the specification errors as $\seri = \poi{\Z} - \poexi{\exoi}$.
It is possible to decompose this quantity into two parts, which may provide more intuition for the error.
This section explores this decomposition.

Dependence between errors can be separated into two components.
The first is the conditioning event itself, capturing the fact that knowledge about $j$'s exposure can provide information about $i$'s outcome in excess of the information provided by $i$'s exposure.
An example is when unit $j$ interferes with unit $i$ in a way that is not captured in $i$'s exposure.
The second source is dependence in excess of what can be explained by the conditioning event.
This captures the fact two units' errors can be dependent if misspecified in the same way even if the exposures themselves provide no information about the outcomes.

We may gain a better understanding about the two components after realizing that the specification errors to some degree are in our control, because we decide how to define the exposures.
Consider when $j$'s exposure provide information about $i$'s outcome in excess of the information provided by its own exposure, meaning that the conditioning event is informative.
A simple way to eliminate this misspecification is to redefine $i$'s exposure to include also the exposure of $j$.
If $i$'s redefined exposure is $\paren{\exoi, \exoj}$, no part of $i$'s specification error can be explained by $j$'s exposure because $i$'s exposure already contains this information.

Conceptually, it is straightforward to gradually remove misspecification by redefining the exposures, but such an approach will often prove impractical.
The exposures would in that case increasingly depart from their primary purpose of producing an effect that is interpretable and relevant for theory or policy.
In particular, if applied to all units in the sample, the redefined exposures would be the intersection of all units' nominal exposures, and much of the reduction in complexity the original exposures provided is lost.
However, the idea of redefined exposures suggests a way to formalize the decomposition of the specification error that will prove useful.

Let $\porexs{ij} : \exall \times \exall \to \Reals$ be a function such that $\porexij{\exeone, \exetwo} = \E{\poi{\Z} \given \exoi = \exeone, \exoj = \exetwo}$.
That is, $\porexij{\exeone, \exetwo}$ is the potential outcome of unit $i$ when defined over the exposures of both $i$ and $j$.
Such a function may not be unambiguously defined if the event $\exoi = \exeone$ and $\exoj = \exetwo$ has measure zero; that is, when $\preijot = \Pr{\exoi = \exeone, \exoj = \exetwo} = 0$.
To accommodate such cases, let the full definition be
\begin{equation}
	\porexij{\exeone, \exetwo} =
		\begin{cases}
		\E{\poi{\Z} \given \exoi = \exeone, \exoj = \exetwo} & \text{if } \preij{\exeone, \exetwo} > 0,
		\\
		\poexi{\exeone} & \text{if } \preij{\exeone, \exetwo} = 0,
		\end{cases}
\end{equation}
This captures the intuition that learning $\exoj = \exetwo$ provides no information about unit $i$'s outcome under $\exoi = \exeone$ if $\exoi = \exeone$ is not simultaneously realizable with $\exoj = \exetwo$.

Because the combination of $\exoi$ and $\exoj$ provides more information about the treatment assignments than $\exoi$ alone, the potential outcome $\porexij{\exeone, \exetwo}$ based on the refined exposures is a more precise representation of unit $i$'s outcome than the potential outcome $\poexi{\exeone}$ based on the original exposures.
We may therefore interpret the difference between $\porexij{\exeone, \exetwo}$ and $\poexi{\exeone}$ as the part of the specification error for unit $i$ explainable by $j$'s exposure.

\begin{definition}[Explainable specification error]\label{def:expl-specification-error}
$\serefij{\exeone, \exetwo} = \porexij{\exeone, \exetwo} - \poexi{\exeone}$.
\end{definition}

While $\porexij{\exeone, \exetwo}$ provides more information than $\poexi{\exeone}$, it will generally not be correctly specified.
That is, $\porexij{\exeone, \exetwo}$ will not provide complete causal information, in the sense that it does not provide the same information as $\poi{\z}$.
The remaining error is that which cannot be explained by $j$'s exposure.
This part is strictly speaking not unexplainable, because the full treatment vector will always perfectly explain the potential outcomes, but it is unexplainable with respect to pairwise refinements of the exposures.
Similar to Definition~\mainref{def:specification-error}, we may define the error not explainable by $j$'s exposure as the difference between the actual potential outcome and the outcome predicted by the redefined exposures.

\begin{definition}[Unexplainable specification error]\label{def:unexpl-specification-error}
$\seruij = \poi{\Z} - \porexij{\exoi, \exoj}$.
\end{definition}

The overall specification error can now be decomposed using the explainable and unexplainable specification errors.
In particular, we have $\seri = \serefij{\exoi, \exoj} + \seruij$ with probability one for any pair of units $i$ and $j$.

Definitions~\ref{def:expl-specification-error}~and~\ref{def:unexpl-specification-error} capture the specification errors pertaining to any particular unit.
The following definition aggregates these errors to an overall description of the specification errors in the experiment as a whole.

\begin{definition}\label{def:decomp-error-dependence}
The \emph{average explainable error dependence} for exposure $\exe \in \exall$ is
\begin{equation}
	\edepe = \frac{1}{n^2} \sumin \sum_{j \neq i} \serefije\serefjie,
\end{equation}
and the \emph{average unexplainable error dependence} for the same exposure is
\begin{equation}
	\udepe = \frac{1}{n^2} \sumin \sum_{j \neq i} \Cov{\seruij, \seruji \given \exoi = \exoj = \exe},
\end{equation}
where $\Cov{\seruij, \seruji \given \exoi = \exoj = \exeone}$ is defined to be zero when $\preij{\exeone, \exetwo} = 0$.
\end{definition}

Definition~\ref{def:decomp-error-dependence} captures pair-wise dependencies between errors of units.
To understand the definitions, first consider the average explainable error dependence.
If $\serefijab = 0$, then knowing that $\exoj = \exeb$ provides no insights about $\ooi$ in excess of knowing that $\exoi = \exea$.
Thus, $\serefije\serefjie$ is non-zero only when the exposures of $i$ and $j$ both provide information about the other unit's outcome.
This means that the magnitude of the explainable errors $\serefije$ matters only insofar that the dependence make them large simultaneously.
If the explainable errors are perfectly symmetric, so that $\serefije = \serefjie$ for all pairs of units, then $\edepe$ collapses to a measure of magnitude.
However, without perfect symmetry, $\edepe$ is a measure of both magnitude and between-unit coordination in the explainable errors.
Indeed, $\edepe$ will be small, or even negative, if the pair-wise explainable errors tend to have opposite signs.
These insights are perhaps made clear by the bound
\begin{equation}
	\edepe \leq \frac{1}{n^2} \sumin \sumjn \bracket{\serefije}^2,
\end{equation}
which shows that the average explainable error dependence is upper bounded by the average magnitude of the explainable errors.

Consider a vaccination trial as an example.
Unit $j$ in this trial is an asymptomatic potential carrier, meaning that $j$ would not get sick if infected but could potentially spread the pathogen to other units.
Unit $i$ on the other hand will show symptoms if infected.
Here, the exposure assigned to $j$ provides information about $i$'s outcome in excess of knowing $i$'s exposure, because $j$'s exposure provides information about whether unit $i$ is infected, and thus shows symptoms.
Part of $i$'s error is thus explainable by $j$'s exposure, and $\serefije$ is non-zero.
However, $i$'s exposure contains no information about $j$'s outcome, because $j$ never shows symptoms, so $\serefjie = 0$.
The lack of symmetry means that there is no dependence between the explainable errors of units $i$ and $j$ according to Definition~\ref{def:decomp-error-dependence}.

Next, consider the average unexplainable error dependence.
The fact that $\udepe$ captures dependence is immediate by the use of a covariance in its definition.
To build intuition, consider the vaccine trial again.
Consider when the exposures capture whether units close to the unit in question are vaccinated (e.g., in their household, or in a neighborhood in a social network).
For illustration, assume that the experiment is so large that the vaccinations in the experiment have the potential to induce herd immunity.
The exposures of any pair of units will in this case provide little information about whether herd immunity is achieved, because whether or not a particular household is vaccinated matters little in that context.
However, if the design of the experiment induces variation in whether head immunity is achieved, then the units' errors will exhibit great dependence even in cases where the explainable errors are small or zero, because pair-wise exposure cannot capture the global behavior.
This example is explored numerically in Section~\ref{sec:example-herd}.

We can now extend the convergence result in the main paper to use the decomposition of the errors.

\newcommand{\propconsistencydecomp}{%
Provided that Conditions~\mainref{cond:bounded-pos},~\mainref{cond:positivity},~\mainref{cond:limited-design-dependence}~and~\mainref{cond:limited-error-dependence} hold for exposures $\exea$ and $\exeb$, the Horvitz--Thompson estimator is consistent for the expected exposure effect and converges at the rate
\begin{equation}
	\estd - \efd
	= \bigOp[\big]{n^{-0.5} + \ddep{\exea}^{0.5} + \ddep{\exeb}^{0.5} + \edep{\exea}^{0.5} + \edep{\exeb}^{0.5} + \udep{\exea}^{0.5} + \udep{\exeb}^{0.5}}.
\end{equation}
}

\begin{proposition}\label{prop:consistency-decompose}
\propconsistencydecomp
\end{proposition}

\subsection{Connection to interference condition in Sävje et al.\ (2021)}\label{sec:conditions-sah}

As noted in Section~\mainref{sec:related-work} in the main paper, the key difference between the interference conditions in this paper and the corresponding conditions in \citet{Saevje2021Average} is that the current paper uses a quantitative concept of interference, while \citet{Saevje2021Average} used a qualitative (or counting) concept.
It is possible to relate the two concepts.
As the counting measure in \citet{Saevje2021Average} is often easier to understand than a quantitative measure, the exercise provides understanding of the measure in the current paper.
Furthermore, the exercise shows that the current measure is a relaxation of the type of measure used in \citet{Saevje2021Average}.
To abstract away from complications that do not provide useful insights for the current purpose, I will consider experimental designs that assign the treatments $\zei$ independently (i.e., Bernoulli designs).

\citet{Saevje2021Average} defined an indicator for each pair of units, $i$ and $j$, that captured what the authors referred to as ``interference dependence.''
The authors used $d_{ij} \in \bset$ to indicate such dependence, but that notation collides with the exposure mapping notation used in this paper, so I will use $\delta_{ij} \in \bset$ to denote the interference dependence indicator in this discussion.
Interference dependence exists between units $i$ and $j$, denoted $\delta_{ij} = 1$, in two situations.
The first situation is when $i$ and $j$ are interfering directly with each other; that is, when $i$'s treatment affects $j$'s outcome, or vice versa.
The second situation is when there exists a third unit $k$ such that $k$'s treatment assignment affects the outcome of both $i$ and $j$.
Both these situations will induce dependence between the outcomes of $i$ and $j$ that is due to the interference.
This could prevent the estimator from concentrating, and \citet{Saevje2021Average} highlight that restricting this type of the interference ensures consistency.
They impose the condition that $\sumin \sumjn \delta_{ij}$ is dominated by $n^2$, meaning that a diminishing fraction of the $n^2$ possible pairs of units are interference dependent.

\citet{Saevje2021Average} considered estimation of an ordinary average treatment effect, corresponding the exposure mapping $\exmiz = \zei$.
We must extend the idea of interference dependence to general exposure mappings to employ a counting concept of interference in the setting considered in this paper.
Also in this setting, there are two situations to consider, corresponding to the two situations above.
The first situation corresponds to the case of direct interference between $i$ and $j$.
In particular, we say that $i$ and $j$ are interference dependent if $j$'s exposure affects $i$'s outcome while holding $i$'s exposure fixed.
Using the notation from Section~\ref{sec:decomposition}, we have $\delta_{ij} = 1$ when $\porexij{\exeone, \exetwo}$ depends on $\exetwo$ while we hold $\exeone$ fixed, or vice versa.
What this tells us is that $j$'s exposure contains information about the treatment of at least one unit that is interfering with $i$, and in this sense, $j$'s exposure interferes directly with $i$'s outcome.
The second situation corresponds to the case when a third unit is interfering with both $i$ and $j$.
That is, if changing some unit $k$'s treatment affects both the outcomes of $i$ and $j$, while holding their exposures fixed, then we set $\delta_{ij} = 1$.
(Often, but not always, this will be equivalent to if changing some unit $k$'s \emph{exposure} affects both the outcomes of $i$ and $j$.)

To connect the concept of interference dependence to specification errors, I will show that $\delta_{ij} = 0$ implies that the expected product of the units' specification errors conditional on their exposures is zero: $\E{\seri\serj \given \exoi = \exoj = \exe} = 0$.
Recall that the expected product of the error is the basis of the interference measure used in the current paper.

First, recall from Section~\ref{sec:decomposition} that we can decompose the error as $\seri = \serefij{\exoi, \exoj} + \seruij$, where $\serefij{\exoi, \exoj}$ is the explainable error and $\seruij$ is the unexplainable error.
If $\delta_{ij} = 0$, then $i$'s exposure is not affecting $j$'s outcome, and vice versa.
This means that $\serefij{\exeone, \exetwo} = \porexij{\exeone, \exetwo} - \poexi{\exeone} = 0$ for all exposures $\exeone$ and $\exetwo$.
Hence, when $\delta_{ij} = 0$, we can write $\seri = \seruij$.
As shown by Lemma~\ref{lem:bound-error-cov-decomp} in Section~\ref{sec:consistency-decompose-proof}, we have $\E{\seruij \given \exoi = \exoj = \exe} = 0$ by construction, so
\begin{equation}
	\E{\seri\serj \given \exoi = \exoj = \exe} = \E{\seruij\seruji \given \exoi = \exoj = \exe} = \Cov{\seruij, \seruji \given \exoi = \exoj = \exe}.
\end{equation}

Recall that the unexplainable error is defined as $\seruij = \poi{\Z} - \porexij{\exoi, \exoj}$.
Hence, under a Bernoulli design, if there is not a third unit interfering with both $i$ and $j$, then $\seruij$ and $\seruji$ are uncorrelated conditional on the exposures of $i$ and $j$:
\begin{equation}
	\Cov{\seruij, \seruji \given \exoi = \exoj = \exe} = 0,
\end{equation}
meaning that $\E{\seri\serj \given \exoi = \exoj = \exe}$ is zero in that case.

Next, consider when $\delta_{ij} = 1$.
The expected product of the units' errors, $\E{\seri\serj \given \exoi = \exoj = \exe}$, will generally not be zero in this setting.
The reason why the quantitative concept of interference used in the current paper is useful is that the expected product can be small (perhaps very small) even if it isn't zero, and that is sufficient for the estimator to be precise in large samples.
What \citet{Saevje2021Average} effectively are doing is to impose a worse-case bound on expected product whenever it is non-zero.
While the approach is somewhat more sophisticated in \citet{Saevje2021Average}, a simple such worse-case bound is to use Condition~\mainref{cond:bounded-pos} in the current paper (i.e., bounded potential outcomes).
We then have $\abs{\E{\seri\serj \given \exoi = \exoj = \exe}} \leq 4 \constpo^2$, as shown in Lemma~\ref{lem:outcome-bounds}.

This worst-case bound allows us to upper bound the interference measure in the current paper with the interference measure in \citet{Saevje2021Average}:
\begin{equation}
\odepe = \frac{1}{n^2} \sumin \sum_{j \neq i} \E{\seri\serj \given \exoi = \exoj = \exe} \leq \frac{4\constpo^2}{n^2} \sumin \sumjn \delta_{ij}.
\end{equation}
The upper bound on the right-hand side is proportional to the quantity $d_{\textsc{avg}} / n$ in the notation of \citet{Saevje2021Average}.
The key condition in \citet{Saevje2021Average} is that $d_{\textsc{avg}} = \littleO{n}$.
Hence, if this condition holds, then $\odepe = \littleO{1}$, which implies that Condition~\mainref{cond:limited-error-dependence} (limited specification error dependence) in the current paper holds.
In other words, the condition in \citet{Saevje2021Average} implies the condition in the current paper, and in that sense, the current condition is weaker.
The intuition for this is exactly that the condition in \citet{Saevje2021Average} only considers whether there is interference between units, but not how strong that interference is.
That is, \citet{Saevje2021Average} require the interference to be strictly local, in the sense that sufficiently few pairs of units can interfere with each other.
The current paper shows that the estimators under study concentrate also when the interference is global, as long as it is either sufficiently rare or sufficiently weak.
This is illustrated in the examples in Section~\ref{sec:error-examples}.

\subsection{Connection to interference condition in Leung (2022)}\label{sec:leung-connection}

\citet{Leung2022Causal} uses a restricted interference assumption based on a graph, referred to as ``approximate neighborhood interference.''
The idea is that units in the experiment are represented as vertices in a graph, and as the shortest path between two vertices grows longer, the interference between the corresponding two units becomes weaker.
This idea is formalized in Assumption~4 in \citet{Leung2022Causal}.

Some additional notation is required to formally state this assumption.
To the greatest degree possible, I use the same notation as in \citet{Leung2022Causal}, with the exception that I suppress the indexing on the sample size ($n$) and the graphs ($A$).
Let $\mathcal{N}(i, s)$ be the subset of unit indices for which the shortest path between unit $i$ and units $j \in \mathcal{N}(i, s)$ is $s$ or less.
In the vernacular of the interference literature, $\mathcal{N}(i, s)$ is the $s$-hop neighborhood around $i$ in the graph.
Let $\Z'$ be a vector of assignments drawn from the experimental design, independent of the actual assignments $\Z$.
Define $\Z^{(i,s)}$ to be equal to $\Z$ on all coordinates $j \in \mathcal{N}(i, s)$ and equal to $\Z'$ on all coordinates $j \not\in \mathcal{N}(i, s)$.
That is, the treatment assignments in $\Z^{(i,s)}$ are the same as the actual treatment assignment for all units in the $s$-hop neighborhood around $i$, and new assignments have been drawn for all units outside of the neighborhood.
Critical to the investigation in \citet{Leung2022Causal} is that the design is Bernoulli, so that all coordinates in $\Z$ are independent.
This ensures that we can construct $\Z^{(i,s)}$ without being concerned about dependencies between assignments inside and outside the $s$-hop neighborhood around $i$.

Now, consider the expected difference in outcomes for unit $i$ under the actual assignment and the artificial assignment just constructed:
\begin{equation}
	\theta_{i,s} = \Es[\Big]{\abs[\big]{\poi[\big]{\Z} - \poi[\big]{\Z^{(i,s)}}}}.
\end{equation}
The quantity $\theta_{i,s}$ is a measure how much units outside of $i$'s $s$-hop neighborhood matters for its outcome.
If the interference is strictly local, in the sense that only units in some $k$-hop neighborhood around $i$ affect $i$'s outcome, then $\theta_{i,s} = 0$ for all $s > k$.
The approximate neighborhood interference assumption stipulates that
\begin{equation}
	\lim_{s\to\infty}\sup_n \max_i \theta_{i,s} = 0,
\end{equation}
where $\sup_n$ denotes the supremum over all samples in the asymptotic sequence.
The assumption ensures that the interference is local in the graph as long as the neighborhood is made sufficiently large.
By itself, the approximate neighborhood interference assumption is somewhat vacuous, because we have $\theta_{i,s} = 0$ for all $i$, $n$ and $s \geq 2$, if we use the complete graph as the interference graph.
Most of the action in \citet{Leung2022Causal} is instead in Assumption~5, which restricts both the topology of the graph and the amount of interference within $s$-hop neighborhoods.
We will return to this assumption below.

To connect the interference conditions in the current paper with those in \citet{Leung2022Causal}, I will derive an upper bound for the overall error dependence (Definition~\mainref{def:error-dependence} in the current paper) in terms of quantities used by \citet{Leung2022Causal}.
I will then show that Assumption~5 in \citet{Leung2022Causal} implies that Condition~\mainref{cond:limited-error-dependence} in the current paper holds; that is, that the overall error dependence diminishes.
This connects the two papers, and specifically shows that the conditions in the current paper are weaker and more general.
This is not meant as a critique of the results in \citet{Leung2022Causal}, as his results are sharper than those in the current paper, and his assumptions might be easier to interpret in some contexts.

Recall the definition of the overall error dependence in the current paper (Definition~\mainref{def:error-dependence}):
\begin{equation}
	\odepe = \frac{1}{n^2} \sumin \sum_{j \neq i} \E{\seri\serj \given \exoi = \exoj = \exe}.
\end{equation}
Let $\mathcal{N}^\partial(i, s)$ be the subset of unit indices for which the shortest path between unit $i$ and units $j \in \mathcal{N}^\partial(i, s)$ is \emph{exactly} $s$.
\citet{Leung2022Causal} calls this the $s$-neighborhood boundary of $i$.
Note that $\mathcal{N}^\partial(i, s)$ for $s \in \braces{1, 2, \dotsc, n}$ partitions the unit indices with the exception of $i$.
(The index $i$ is excluded because $\mathcal{N}^\partial(i, 0) = \braces{i}$.)
Using this partition, we can write the overall error dependence as
\begin{equation}
	\odepe = \frac{1}{n^2} \sum_{s = 1}^n \sumin \sum_{j \in \mathcal{N}^\partial(i, s)} \E{\seri\serj \given \exoi = \exoj = \exe}.
\end{equation}

Similar to the current paper, \citet{Leung2022Causal} considers generic exposure mappings.
However, unlike the current paper, it is critical in \citet{Leung2022Causal} that the exposure mappings are local to the observed graph.
This is formalized in Assumption~1 in \citet{Leung2022Causal}, which states that unit $i$'s exposure can only depend on the treatment assignments in a $K$-hop neighborhood around $i$ in the graph, where $K$ is a constant in the asymptotic sequence.
We will break up the summation in the overall error dependence into two parts, the first over $1, \dotsc, 2K$ and the second over $2K + 1, \dotsc, n$:
\begin{multline}
	\odepe = \frac{1}{n^2} \sum_{s = 1}^{2K} \sumin \sum_{j \in \mathcal{N}^\partial(i, s)} \E{\seri\serj \given \exoi = \exoj = \exe}
	\\
	+ \frac{1}{n^2} \sum_{s = 2K + 1}^n \sumin \sum_{j \in \mathcal{N}^\partial(i, s)} \E{\seri\serj \given \exoi = \exoj = \exe}.
	\label{eq:leung-split}
\end{multline}

Consider the first term in expression~\eqref{eq:leung-split}.
Because the potential outcomes are bounded, which is an assumption made both in the current paper and in \citet{Leung2022Causal}, we can apply Lemma~\ref{lem:outcome-bounds} of the current paper, which states that $\E{\seri\serj \given \exoi = \exoj = \exe} \leq 4\constpo^2$.
This allows us to write
\begin{equation}
	\frac{1}{n^2} \sum_{s = 1}^{2K} \sumin \sum_{j \in \mathcal{N}^\partial(i, s)} \E{\seri\serj \given \exoi = \exoj = \exe} \leq
	\frac{1}{n^2} \sum_{s = 1}^{2K} \sumin \sum_{j \in \mathcal{N}^\partial(i, s)} 4\constpo^2 =
	\frac{4\constpo^2}{n} \sum_{s = 1}^{2K} M^\partial(s),
\end{equation}
where $M^\partial(s) = n^{-1} \sumin \abs{\mathcal{N}^\partial(i, s)}$ is the average size of the $s$-neighborhood boundary.
The quantity $M^\partial(s)$ is used by \citet{Leung2022Causal}, in Assumption~5 and elsewhere.

Next, consider the second term in expression~\eqref{eq:leung-split}.
Recall the definition of the specification errors (Definition~\mainref{def:specification-error}): $\seri = \poi{\Z} - \poexi{\exoi}$.
Let
\begin{equation}
	A_i(s) = \poi[\big]{\Z} - \poi[\big]{\Z^{(i,s)}}
	\qquadand
	B_i(s) = \poi[\big]{\Z^{(i,s)}} - \poexi[\big]{\exoi},
\end{equation}
so that the specification error can be written $\seri = A_i(s) + B_i(s)$ for any $s$.
For any $i$ and any $j \in \mathcal{N}^\partial(i, s)$ such that $s \geq 2K + 1$, we can write
\begin{equation}
	\E{\seri\serj \given \exoi = \exoj = \exe}
	= \Es[\Big]{B_i\big(\floor{s/2}\big) B_j\big(\floor{s/2}\big) + C_{ij}(s) \given \exoi = \exoj = \exe},
\end{equation}
where
\begin{equation}
	C_{ij}(s) = \serj A_i\big(\floor{s/2}\big) + A_j\big(\floor{s/2}\big) B_i\big(\floor{s/2}\big).
\end{equation}

Note that $B_i(\floor{s/2})$ depends on $\Z^{(i,\floor{s/2})}$ and $\exoi$.
When $s \geq 2K + 1$, both of these random variables only depend on assignments in $\Z$ within a $\floor{s/2}$-hop neighborhood around $i$.
Furthermore, if $j \in \mathcal{N}^\partial(i, s)$, then the $\floor{s/2}$-hop neighborhood around $i$ is disjoint from the $\floor{s/2}$-hop neighborhood around $j$.
The implication is that $B_i(\floor{s/2})$ and $B_j(\floor{s/2})$ are independent, because the design is Bernoulli.
This means that we can write
\begin{equation}
	\Es[\Big]{B_i\big(\floor{s/2}\big) B_j\big(\floor{s/2}\big) \given \exoi = \exoj = \exe}
	=
	\Es[\Big]{B_i\big(\floor{s/2}\big) \given \exoi = \exe} \Es[\Big]{B_j\big(\floor{s/2}\big) \given \exoj = \exe}.
\end{equation}
These two factors are both zero, because
\begin{equation}
\Es[\Big]{B_i\big(\floor{s/2}\big) \given \exoi = \exe}
= \Es[\Big]{\poi[\big]{\Z^{(i,s)}} \given \exoi = \exe} - \poexi{\exe}
= \poexi{\exe} - \poexi{\exe} = 0.
\end{equation}
We therefore have
\begin{equation}
	\E{\seri\serj \given \exoi = \exoj = \exe}
	= \Es[\big]{C_{ij}(s) \given \exoi = \exoj = \exe}.
\end{equation}

Because \citet{Leung2022Causal} assumes that the exposures are defined within a $K$-hop neighborhood and because he assumes that the design is Bernoulli, two units that are at a geodesic distance of at least $2K + 1$ will have independent exposures.
This tells us that, for any $i$ and any $j \in \mathcal{N}^\partial(i, s)$ such that $s \geq 2K + 1$,
\begin{equation}
	\Pr{\exoi = \exoj = \exe} = \Pr{\exoi = \exe} \Pr{\exoj = \exe},
\end{equation}
where the probabilities on the right-hand side are bounded away from zero by the positivity assumption made in both the current paper and in \citet{Leung2022Causal}.
We can therefore write
\begin{multline}
	\E{\seri\serj \given \exoi = \exoj = \exe}
	\leq \Es[\big]{\abs{C_{ij}(s)} \given \exoi = \exoj = \exe}
	\\
	\leq \constpi^2 \Pr{\exoi = \exoj = \exe} \Es[\big]{\abs{C_{ij}(s)} \given \exoi = \exoj = \exe}
	\leq \constpi^2 \Es[\big]{\abs{C_{ij}(s)}},
\end{multline}
where the bound $1 / \Pr{\exoi = \exe} \leq \constpi$ was used.

Using Hölder's inequality, we can write
\begin{equation}
	\Es[\big]{\abs{C_{ij}(s)}} \leq \Es[\big]{\abs{A_i(\floor{s/2})}} \esssup \abs{\serj} + \Es[\big]{\abs{A_j(\floor{s/2})}} \esssup \abs{B_i(\floor{s/2})},
\end{equation}
where $\esssup X$ denotes the essential supremum of the random variable $X$.
Note that $\Es[\big]{\abs{A_i(\floor{s/2})}}$ is equal to $\theta_{i,s}$, as defined in the beginning of this section.
Also note that bounded potential outcomes implies that both $\esssup \abs{\serj}$ and $\esssup \abs{B_i(\floor{s/2})}$ are bounded by $2\constpo$.
Let $\theta_{s} = \max_{i} \theta_{i,s}$, which corresponds to the definition in the beginning of Section~3 in \citet{Leung2022Causal}.
This allows us to write
\begin{equation}
	\Es[\big]{\abs{C_{ij}(s)}} \leq 4 \constpo \theta_{\floor{s/2}}.
\end{equation}

Taken together, we have
\begin{equation}
	\E{\seri\serj \given \exoi = \exoj = \exe}
	\leq 4 \constpo \constpi^2 \theta_{\floor{s/2}}
\end{equation}
for any $i$ and any $j \in \mathcal{N}^\partial(i, s)$ as long as $s \geq 2K + 1$.
In turn, this gives us that the second term in expression~\eqref{eq:leung-split} is bounded as
\begin{multline}
	\frac{1}{n^2} \sum_{s = 2K + 1}^n \sumin \sum_{j \in \mathcal{N}^\partial(i, s)} \E{\seri\serj \given \exoi = \exoj = \exe}
	\\
	\leq
	\frac{4 \constpo \constpi^2}{n^2} \sum_{s = 2K + 1}^n \sumin \sum_{j \in \mathcal{N}^\partial(i, s)} \theta_{\floor{s/2}}
	=
	\frac{4 \constpo \constpi^2}{n} \sum_{s = 2K + 1}^n M^\partial(s) \theta_{\floor{s/2}},
\end{multline}
where, as above, $M^\partial(s) = n^{-1} \sumin \abs{\mathcal{N}^\partial(i, s)}$.

Combining the results for the two terms in expression~\eqref{eq:leung-split}, we have
\begin{equation}
	\odepe \leq \frac{4\constpo^2}{n} \sum_{s = 1}^{2K} M^\partial(s)
	+ \frac{4 \constpo \constpi^2}{n} \sum_{s = 2K + 1}^n M^\partial(s) \theta_{\floor{s/2}}.
\end{equation}
To connect the bound to the assumptions in \citet{Leung2022Causal}, it will prove convenient to make it less tight by harmonizing the constants for the two terms and to extend the summation to include the case $s = 0$:
\begin{equation}
	\odepe \leq \frac{4\constpo^2 \constpi^2}{n} \sum_{s = 0}^{2K} M^\partial(s)
	+ \frac{4 \constpo^2 \constpi^2}{n} \sum_{s = 2K + 1}^n M^\partial(s) \theta_{\floor{s/2}}.
\end{equation}
Next, define $\tilde{\theta}_s$ as
\begin{equation}
	\tilde{\theta}_s = \begin{cases}
		1 & \text{if } s \leq 2K, \\
		\theta_{\floor{s/2}} & \text{if } s > 2K.
	\end{cases}
\end{equation}
Note that this corresponds exactly to the definition of the same symbol in Theorem~1 in \citet{Leung2022Causal}.
This allows us to write
\begin{equation}
	\odepe \leq \frac{4\constpo^2 \constpi^2}{n} \sum_{s = 0}^{2K} M^\partial(s)
	+ \frac{4 \constpo^2 \constpi^2}{n} \sum_{s = 2K + 1}^n M^\partial(s) \theta_{\floor{s/2}}
	=
	\frac{4\constpo^2 \constpi^2}{n} \sum_{s = 0}^n M^\partial(s) \tilde{\theta}_s.
\end{equation}

Now, Assumption~5 in \citet{Leung2022Causal} is that
\begin{equation}
	\sum_{s = 0}^n M^\partial(s) \tilde{\theta}_s = \littleO{n}.
\end{equation}
Given the bound we have just derived, this implies that $\odepe = \littleO{1}$, which in turn implies that Condition~\mainref{cond:limited-error-dependence} in the current paper holds.
In other words, in the context studied by \citet{Leung2022Causal}, the assumptions he makes imply the interference condition used in the current paper.
In this sense, the conditions in the current paper are weaker and more general than those in \citet{Leung2022Causal}.
As noted above, this is not meant as a critique; instead, it highlights that the scopes of the results are different.
Note that \citet{Leung2022Causal} uses Assumption~6, which is a stronger version of Assumption~5, to prove asymptotic normality.
Hence, he uses even stronger assumptions than the conditions in the current paper to derive this result, similar to how stronger assumptions are used in Section~\ref{sec:limiting-distribution} of the current paper to show asymptotic normality.

\subsection{Numerical examples}\label{sec:error-examples}

\subsubsection{Setup}

This section provides concrete, numerical examples to illustrate the specification errors and their connection to the precision of the estimators.

The basic setup is the same in all examples.
There is an observed network through which units are hypothesized to interact.
This network is used to define the exposures.
Denote an edge from unit $i$ to unit $j$ in the network as $g_{ij} = g_{ji} = 1$, and the absence of an edge is denoted as $g_{ij} = 0$.
For simplicity, the graph is a cycle graph, meaning that $g_{ij} = 1$ if and only if $\abs{i - j} \in \braces{1, n - 1}$.
That is, unit $i$ is connected to units $i - 1$ and $i + 1$, with the exception of units $1$ and $n$, who are connected to each other.
The qualitative results of this simulation do not depend on the use of a cycle graph; it is used purely for convenience.

Each unit is assigned a binary treatment independently at random with equal probability: $\Pr{\Zei = 1} = 1/2$ for all $i$.
That is, the experimental design is Bernoulli.
There are two exposures of interest.
The first exposure, denoted $1$, is when the ego is untreated and both neighbors in the cycle graph are treated.
The second exposure, denoted $0$, is when the ego as well as both neighbors are untreated.
All other treatment vectors are mapped to a residual exposure, denoted $2$, which will not be used in the analysis.
That is, the exposure mapping for unit $i$ is
\begin{equation}
	\exmiz = \begin{cases}
		0 & \text{if } \zei = 0 \text{ and } \sum_{j=1}^n g_{ij} \zej = 0, \\
		1 & \text{if } \zei = 0 \text{ and } \sum_{j=1}^n g_{ij} \zej = 2, \\
		2 & \text{else.}
	\end{cases}
\end{equation}
Recall that $\exoi = \exmi{\Z}$ denotes the realized exposure for unit $i$.
Because of the Bernoulli design and cycle graph structure, the exposure probabilities are $12.5\%$ for both exposures and for all units:
\begin{equation}
	\prei{1} = \Pr{\exoi = 1} = 1/8
	\qquadand
	\prei{0} = \Pr{\exoi = 0} = 1/8.
\end{equation}

The treatment effect of interest is
\begin{equation}
	\ef{1,0}
	= \frac{1}{n} \sumin \braces[\big]{\poexi{1} - \poexi{0}},
\end{equation}
which is the expected indirect effect of having two treated neighbors when being untreated compared to having no treated neighbors.
Because of its improved small sample performance, the main focus of the simulation study will be the H{\'a}jek estimator, which in this case coincide with the difference-in-means estimator.
The results for the Horvitz--Thompson estimator are reported in Table~\ref{tab:simulation-results} at the end of the section on page~\pageref{tab:simulation-results}.
The simulation consists of 50,000 rounds (for each setting and sample size), ensuring that the Monte Carlo errors are negligible, and no uncertain measures for Monte Carlo imprecision are reported.

\subsubsection{Local interference}\label{sec:example-local}

The interference is local in the first example.
This largely mirrors the setting studied by \citet{Saevje2021Average}, and the discussion in Section~\ref{sec:conditions-sah} tells us that the estimators should be consistent here.
Two version of the data generating process will be considered.

In the first version, labelled ``Correct,'' the exposures are correctly specified with respect to the two exposures of interest.
In particular, the potential outcome function is
\begin{equation}
	\poiz = \alpha_i - \beta_i \zei - \gamma_i \sum_{j=1}^n g_{ij} \zej,
\end{equation}
where $\alpha_i$ is drawn uniformly at random from $[15, 25]$, $\beta_i$ is drawn from $[1, 10]$ and $\gamma_i$ is drawn from $[1, 2]$.
The coefficients are drawn once and held fixed over the simulation rounds, mirroring the fact that all randomness under consideration comes from treatment assignment.
Here, we have $\poiz = \poi{\z'}$ whenever $\exmi{\z} = \exmi{\z'} \in \braces{0, 1}$, meaning that the exposures of interest, $d \in \braces{0, 1}$, are correctly specified.
That is, the overall error dependence in Definition~\mainref{def:error-dependence} is zero.

We can interpret this data generating process to describe a vaccination study.
The treatment is whether a unit has received the vaccine, and the outcome is some measure of viral load.
The coefficient $\alpha_i$ captures $i$'s baseline exposure and susceptibility to the virus, $\beta_i$ captures the direct effect of taking the vaccine oneself, and $\gamma_i$ captures the indirect effect of having vaccinated neighbors.
Of course, this is a highly stylized data generating process that is unlikely to capture the dynamics of virus transmission in reality, but the process fills its role as an illustration.

The second version, labelled ``Local,'' has unmodeled interference, meaning that the exposures are not correctly specified.
An additional network is added here, which is directed and unobserved by the investigator.
The unobserved network is generated by randomly and independently drawing an arc from $i$ to $j$ with probability $n^{-4/5}$.
Let $u_{ij} \in \bset$ denote whether there is an arc from $i$ to $j$ in this network.
If an unobserved arc overlaps with an edge in the observed network, the unobserved arc is deleted, ensuring that $g_{ij}u_{ij}=0$.
The unobserved network is generated once and held fixed over the simulation rounds.
Because the observed network is completely unrelated to the unobserved network, approximate neighborhood interference, as defined in \citet{Leung2022Causal} and discussed in Section~\ref{sec:leung-connection} above, does not hold here.

The potential outcome functions under the second version is
\begin{equation}
	\poiz = \alpha_i - \beta_i \zei - \gamma_i \sum_{j=1}^n \paren{g_{ij} + u_{ij}} \zej.
\end{equation}
That is, vaccinated neighbors have the same effect in the observed and unobserved networks.
The coefficients used in the first version of the data generating process are used here as well.
We can interpret the second data generating process as also describing a vaccination study, but with unobserved transmission paths.

Note that the exposures are not correctly specified under the second version of the process.
We have that the expected potential outcomes for unit $i$ under the two exposures of interest, $d \in \bset$, are
\begin{equation}
	\poexi{d} = \alpha_i - 2 \gamma_i d - \frac{\gamma_i}{2} \sum_{j=1}^n u_{ij},
\end{equation}
meaning that the specification error for unit $i$ conditional on $\exoi \in \bset$ is
\begin{equation}
	\seri = \frac{\gamma_i}{2} \sum_{j=1}^n u_{ij} \paren{1 - 2 \zej}.
\end{equation}
However, because the unobserved network is sparse, the overall error dependence in Definition~\mainref{def:error-dependence} will diminish as the sample size grows.
Indeed, the overall error dependence will be of order $n^{-16/25}$ in this setting, ensuring that Condition~\mainref{cond:limited-error-dependence} is satisfied.
The rate of the overall error dependence can be derived either by direct calculation or by using the approach described in Section~\ref{sec:conditions-sah} above.

\begin{figure}
	\centering
	\includegraphics[width=\textwidth]{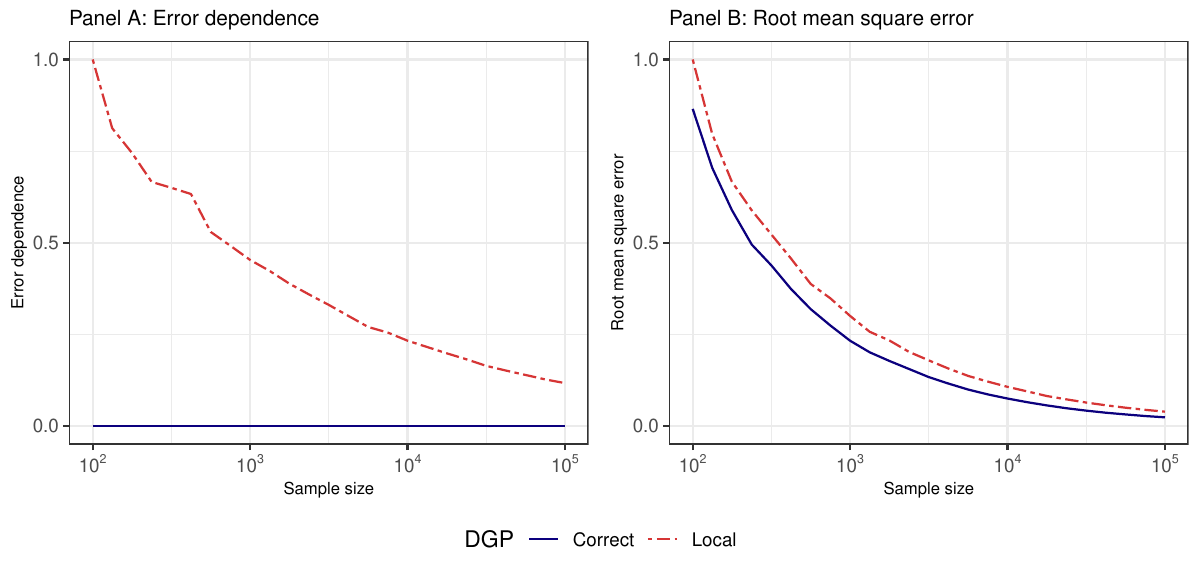}
	\caption{Root of the sum of overall error dependence (A) and root mean square error of the H{\'a}jek estimator (B) under local interference.  The lines correspond to different versions of the data generating process, as described in the text. The measures are normalized so that they start at $1$ for one of the versions. The simulation is based on 50,000 draws from the experimental design, and simulation error is negligible.}\label{fig:local-interference}
\end{figure}

The results are presented in Figure~\ref{fig:local-interference}.
Panel A presents the root of the sum of the magnitude of the overall error dependence measures defined in the paper for the two exposures, $\sqrt{\abs{\odep{1}} + \abs{\odep{0}}}$.
For the first data generating process (``Correct''), the error dependence is constant at zero, indicating that the exposures indeed are correctly specified under this process.
The measure is not zero under the second data generating process (``Local'').
However, the error dependence decreases at a reasonably fast rate, indicating that Condition~\mainref{cond:limited-error-dependence} indeed is satisfied, and that we can expect the estimators to be precise in large samples.

Panel B presents the root mean square errors of under the two versions of the data generating process.
The estimator is somewhat less precise under the second version for all sample sizes, but the mean square error approaches zero at a reasonable rate under both versions.
This is what we would expect given Proposition~\mainref{prop:variance-bound} in the main paper and the diminishing error dependence shown in Panel A.

\subsubsection{Global interference: Herd immunity}\label{sec:example-herd}

The second example introduces global interference.
That is, all units can potentially interfere with all other units.
The example builds on the first data generating process in the previous section by adding behavior akin to herd immunity.
If the number of vaccinated units passes some threshold, there is no transmission of the virus in the community under study, and the outcome (e.g., viral load) is zero for all units.
This is again a highly stylized example, but it suffices to illustrate the key insights.

The potential outcome function is here
\begin{equation}
	\poiz = \begin{cases}
		0 & \text{if } \sum_{j=1}^n \zej \geq \phi n,
		\\
		\alpha_i - \beta_i \zei - \gamma_i \sum_{j=1}^n g_{ij} \zej & \text{else},
	\end{cases}
\end{equation}
where the same coefficients as in the previous data generating processes are used, and $\phi \in (0, 1)$.
This process captures behavior reminiscent of herd immunity because if a share of $\phi$ units are treated, then the outcome is zero for all units.

Three versions of this data generating process will be considered, corresponding to three values of $\phi$.
The first version, labelled ``Rarely,'' uses $\phi = 0.52$, meaning that herd immunity occurs when at least $52\%$ of the units are treated.
The second version, labelled ``Infrequently,'' uses $\phi = 0.51$, and the third version, labelled ``Frequently,'' uses $\phi = 0.5$.

Note that $\sum_{j=1}^n \zej$ follows a binomial distribution of $n$ trials with $0.5$ success probability.
Hence, herd immunity will occur somewhat rarely under the first version of the data generating process ($\phi = 0.52$), and it grows increasingly rare as the sample size grows.
Herd immunity is more common, but still fairly infrequent, when $\phi = 0.51$.
However, when $\phi = 0.5$, herd immunity will be frequent also in large samples, because the distribution of $\sum_{j=1}^n \zej$ is centered around $n / 2$.
Indeed, when $n = 100$, herd immunity occurs with roughly the same probability under all three versions of the process.
But when $n = 10,000$, herd immunity occurs with a probability of only $0.003\%$ under the first version ($\phi = 0.52$), a probability of $2.33\%$ under the second version ($\phi = 0.51$), and a probability of $50.4\%$ under the third version ($\phi = 0.50$).

Interference is global under all three versions, because there are situations in which changing a single unit's treatment assignment changes the outcome of all other units under all three versions.
That is, the event $\sum_{i=1}^n \Zei < \phi n \leq 1 + \sum_{i=1}^n \Zei$ has a non-zero probability of occurring under all three versions, and changing any untreated unit to be treated in this case induces herd immunity, changing the outcomes of all other units.
However, such global interference takes place only very rarely under the first version when the sample is large, because changing a single unit's treatment assignment under most assignments will not have global effects.
In this sense, global interference does exist under all three versions, but it is not practically relevant for the first two versions of the data generating process, $\phi \in \braces{0.51, 0.52}$, if the sample is sufficiently large, because it will happen with a probability so low that it can be ignored.
This is not the case when $\phi = 0.5$, because there will be non-negligible variation in whether herd immunity occurs no matter how large the sample gets.

To see this formally, let $H = \indicator[\big]{\sum_{i=1}^n \Zei \geq \phi n}$ be an indicator whether herd immunity occurs, and let $\pi_h = \E{H \given \Ze{1} = \Ze{2} = \Ze{3} = 0}$ be the probability of that event conditional on that three units are untreated.
The expected potential outcomes for unit $i$ under exposure $d = 0$ is
\begin{equation}
	\poexi{0} = \paren{1 - \pi_h} \alpha_i,
\end{equation}
and the specification error for unit $i$ conditional on $\exoi = 0$ is
\begin{equation}
	\seri = \paren{\pi_h - H} \alpha_i.
\end{equation}
This implies that the overall error dependence will be of the same order as $\Var{H}$.
That is, if there is variability in whether herd immunity occurs also in large samples, then the overall error dependence will not approach zero.
The variance will be larger when $\phi = 0.51$ compared to $\phi = 0.52$, but we have $\Var{H} \to 0$ in both cases.
However, when $\phi = 0.5$, $\Var{H}$ is asymptotically bounded away from zero, meaning that the overall error dependence does not approach zero in this setting.
Hence, given Proposition~\mainref{prop:variance-bound}, we expect convergence of the estimator when $\phi = 0.52$ and $\phi = 0.51$, but possibly not when $\phi = 0.5$.

\begin{figure}
	\centering
	\includegraphics[width=\textwidth]{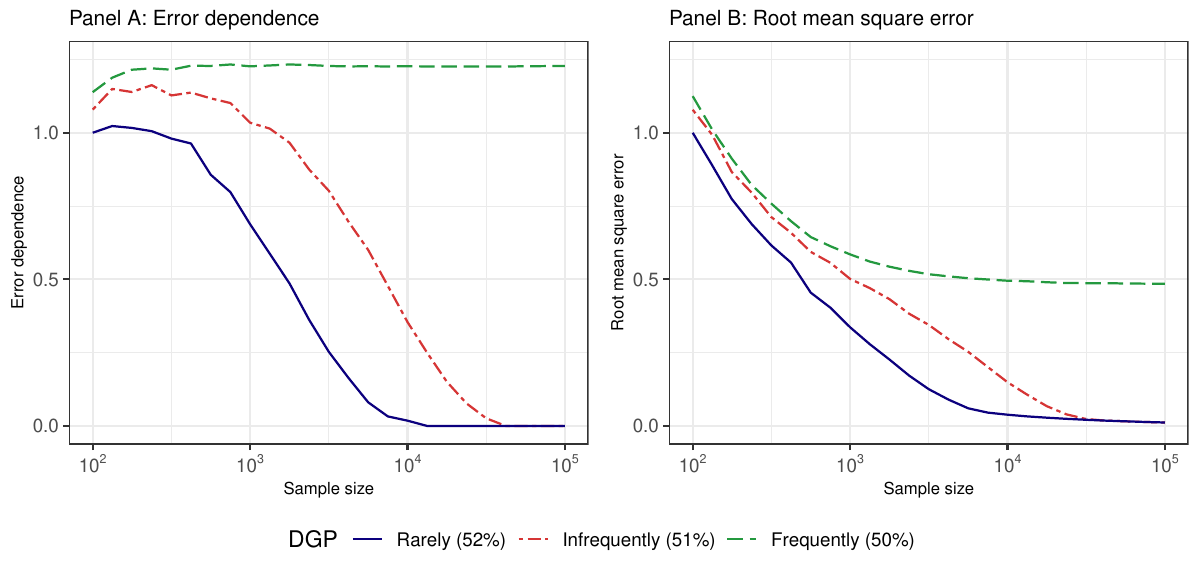}
	\caption{Root of the sum of overall error dependence (A) and root mean square error of the H{\'a}jek estimator (B) under interference through herd immunity. See the note of Figure~\ref{fig:local-interference} for additional details.}\label{fig:herd-interference-supp}
\end{figure}

The results are presented in Figure~\ref{fig:herd-interference-supp}.
As in the previous figure, Panel A presents the sum of the magnitude of the overall error dependence measures for the two exposures.
The error dependence increases initially, but as the sample grows in size, the error dependence decreases for both ``Rarely'' ($\phi = 0.52$) and ``Infrequently'' ($\phi = 0.51$).
This mirrors the fact that variation in whether herd immunity occurs decreases as the sample size grows under both of these versions.
However, the error dependence does not decrease for the version ``Frequently,'' instead flattening out slightly below the 1.25 mark.
Therefore, Condition~\mainref{cond:limited-error-dependence} does not hold, and the consistency result of the current paper does not apply.

Panel B presents the root mean square errors of under the three version of the data generating process.
The precision of the estimator initially improves for all three versions.
The precision continues to improve for the first two versions (``Rarely'' and ``Infrequently''), and the mean square error is very close to zero when $n=100,000$.
However, under the third version (``Frequently''), the mean square error flattens out at around the $0.5$ mark, indicating that the estimator is not consistent.
This mirrors what we learned about the behavior of the error dependence measure in Panel A.

The reason the precision of the estimator initially improves also under the third version, despite not seeing a corresponding decrease in the error dependence, is that two first terms of the variance bound in Proposition~\mainref{prop:variance-bound} approach zero as the sample size grows even if the error dependence is large.

\subsubsection{Global interference: General equilibrium}\label{sec:example-market}

The third example considers interference transmitted through a market price.
There are $n$ companies producing some good.
The treatment under study is an improved production technology.
There are spillovers in the production technology, in the sense that a company's production is improved if other neighboring companies have access to the improved technology.
This could, for example, arise because technology-specific inputs (e.g., skilled labor) become more accessible in the local market if many companies use the same technology, because of technology transmission, or because of some other cluster effect in production.

Given market price $p$ and treatments $\z$, the production function for company $i$ is
\begin{equation}
	q_i(p, \z) = \gamma_i p + \gamma_i p^2 \paren[\bigg]{\zei + \frac{1}{2}\sum_{j=1}^n g_{ij} \zej}.
\end{equation}
That is, in absence of treatment, $\z = \boldsymbol{0}$, the production function is linear.
If the company itself or any of its neighbors have access to the improved technology, the production function adds a quadratic component, making the production more efficient.
The parameter $\gamma_i$ captures the size of the company, scaling the production function so that larger companies produce more goods, all else equal.
The size of the supply side of the market is normalized to $100$, in the sense that $\sumin \gamma_i = 100$.

The supply curve for the market is given by
\begin{equation}
	Q_s(p, \z) = \sumin q_i(p, \z) = 100 p + a(\z) p^2,
\end{equation}
where
\begin{equation}
	a(\z) = \sumin \gamma_i \zei + \frac{1}{2}\sumin \sum_{j=1}^n g_{ij} \gamma_i \zej.
\end{equation}
The demand curve is
\begin{equation}
	Q_d(p) = 100 / p.
\end{equation}
Given treatments $\z$, the equilibrium price $p^*(\z)$ is the one that clears the market:
\begin{equation}
	p^*(\z) = \argmin_{p \in \Reals^+} \abs{Q_s(p, \z) - Q_d(p)}.
\end{equation}
Note that when no company is treated, $\z = \boldsymbol{0}$, the equilibrium price is one, $p^*(\z) = 1$, because the market then clears when $100 p = 100 / p$.
This was the price at baseline, before the experiment was run.
If one or more companies are treated, the price will be lower than one.
The exact price depends on both the number of treated companies and which companies are treated.
As more large companies are treated (or large companies have treated neighbors), the lower the equilibrium price will be.
Changing any company's treatment will always change the equilibrium price at least slightly.

Given treatments $\z$, the production in company $i$ is
\begin{equation}
q_i\braces{p^*(\z), \z} = \gamma_i p^*(\z) + \gamma_i \braces{p^*(\z)}^2 \paren[\bigg]{\zei + \frac{1}{2}\sum_{j=1}^n g_{ij} \zej}.
\end{equation}
Note that treatment has two effects.
Holding the market price fixed, companies with access to the improved technology will have an increased production.
But as more companies have access to the technology, the market price will decrease, leading to decreased production.
The overall production will increase as more companies have access to the technology, but some of the companies might still decrease their production.

Company $i$'s baseline production (when $\z = \boldsymbol{0}$) was
\begin{equation}
q_i\braces{p^*(\boldsymbol{0}), \boldsymbol{0}} = \gamma_i.
\end{equation}
The outcome of interest is the relative increase in production compared to baseline:
\begin{equation}
	\poiz = \frac{q_i\braces{p^*(\z), \z}}{q_i\braces{p^*(\boldsymbol{0}), \boldsymbol{0}}} = p^*(\z) + \braces{p^*(\z)}^2 \paren[\bigg]{\zei + \frac{1}{2}\sum_{j=1}^n g_{ij} \zej}.
\end{equation}

Unlike the example in Section~\ref{sec:example-local}, the interference is global here.
That is, changing a single unit's treatment could change the outcome of all other units.
Unlike the example with herd immunity in Section~\ref{sec:example-herd}, which also considered global interference, it is common that interference occurs globally.
That is, changing a single unit's treatment will \emph{always} change the outcome of all other units.
However, it is still possible to achieve consistency here, because the interference could be limited.
In this case, all global interference is mediated through the market price.
If the market price stabilizes in large samples, in the sense that it is close to some value with high probability, then the interference between most units will be negligible.
This is because the production functions are smooth in the market price, so minor price disturbances will be inconsequential for production.

Three versions of this data generating process will be considered, differing in the distribution of the companies' sizes.
In the first version, labelled ``All small,'' all companies are reasonably small: $\gamma_i$ for unit $i$ is proportional to the $i / (n + 1)$-th percentile of the standard log-normal distribution.
When $n=1000$, the largest company has $\gamma_i = 1.35$, corresponding to a $1.35\%$ share of the total market at baseline.
The five largest companies have a $5.07\%$ share of the market at baseline when $n=1000$.
The market share of the top five companies will approach zero as the sample size grows.

In the second version, labelled ``Some outliers,'' most companies are small but there are some large outliers.
The coefficients $\gamma_i$ are here proportional to the $i / (n + 1)$-th percentile of the log-normal distribution with mean parameter zero and standard deviation parameter three.
When $n=1000$, the largest company has $\gamma_i = 19.6$, corresponding to a $19.6\%$ share of the total market at baseline.
The five largest companies have a $46.3\%$ share of the market at baseline.
The market share of the top five companies will approach zero as the sample size grows also in this case, but at a slower rate than in the first version.

In the third version, labelled ``One large,'' most companies are small but there is one company that has half of the market: $\gamma_i$ is proportional to the percentiles of the standard log-normal distribution, as in the first version, except for the first company, which instead has $\gamma_1 = 50$.
For all sample sizes, the largest company has $\gamma_i = 50$, meaning a $50\%$ share of the total market at baseline.
When $n=1000$, the five largest companies have a $52.1\%$ baseline share of the market.
The market share of the top five companies will not approach zero as the sample size grows in this case.

The market price will stabilize in large markets with many small companies.
Treatment assignment will introduce variability for individual companies, but the aggregated supply curve will be similar (but never identical) under most assignments.
However, when there is one or a few very large companies in the market, the treatment assigned to those companies, or their neighbors, will have a large effect of the aggregated supply curve, and thus the market price.
This is true even if the market as a whole is large.
Hence, we expect the market price to stabilize under the first two versions of the data generating process (``All small'' and ``Some outliers''), but the market price might stabilize slowly when there are outliers as in the second version.
The price will not stabilize under the third version (``One large'').
The treatments assigned to the one large company and its neighbors will have a large effect on the market price, even if the market is very large.

To see this formally, consider the expected potential outcomes for unit $i$ under exposure $d = 0$:
\begin{equation}
	\poexi{0} = \Es[\big]{\poi{\Z} \given \exoi = 0} = \Es[\big]{p^*(\Z) \given \exoi = 0},
\end{equation}
which is the expected market price conditional on that unit $i$ has exposure $d = 0$.
For ease of exposition, assume that $i$ is a small company with neighbors that also are small.
This means that the conditioning event $\exoi = 0$ is inconsequential for the distribution of the market price in large samples, and we can approximate the conditional expectation with the unconditional expected market price $\Et[\big]{p^*(\Z)}$.
That is, we have $\poexi{0} \approx \Et[\big]{p^*(\Z)}$ with a negligible approximation error in large samples.
Note that most companies will be small with small neighbors, even if there is one or a few large companies, so this case covers the almost all units in all three versions of this data generating process.

The specification error for unit $i$ conditional on $\exoi = 0$ then becomes
\begin{equation}
	\seri = p^*(\Z) - \Es[\big]{p^*(\Z)}.
\end{equation}
This means that the error dependence between two small companies $i$ and $j$ approximately is
\begin{equation}
	\E{\seri\serj \given \exoi = \exoj = 0}
	\approx \E[\Big]{\bracket[\big]{ p^*(\Z) - \Es[\big]{p^*(\Z)}}^2 } = \Vars[\big]{ p^*(\Z) }.
\end{equation}
Like above, the approximation here uses the fact that the conditioning event $\exoi = \exoj = 0$ is inconsequential when companies $i$ and $j$, and their neighbors, are small.
The number of companies that are large or have large neighbors will be a diminishing share of all companies under all three versions of the data generating process.
The outcome is bounded for all companies, including large ones, meaning that $\E{\seri\serj \given \exoi = \exoj = 0}$ is bounded for all units (see Lemma~\ref{lem:outcome-bounds}).
This means the error dependence involving large companies can be ignored in large samples for all three versions, and we have
\begin{equation}
	\odep{0} = \frac{1}{n^2} \sumin \sum_{j \neq i} \E{\seri\serj \given \exoi = \exoj = 0} = \Vars[\big]{ p^*(\Z) } + \littleO{1}.
\end{equation}

\begin{figure}
	\centering
	\includegraphics[width=\textwidth]{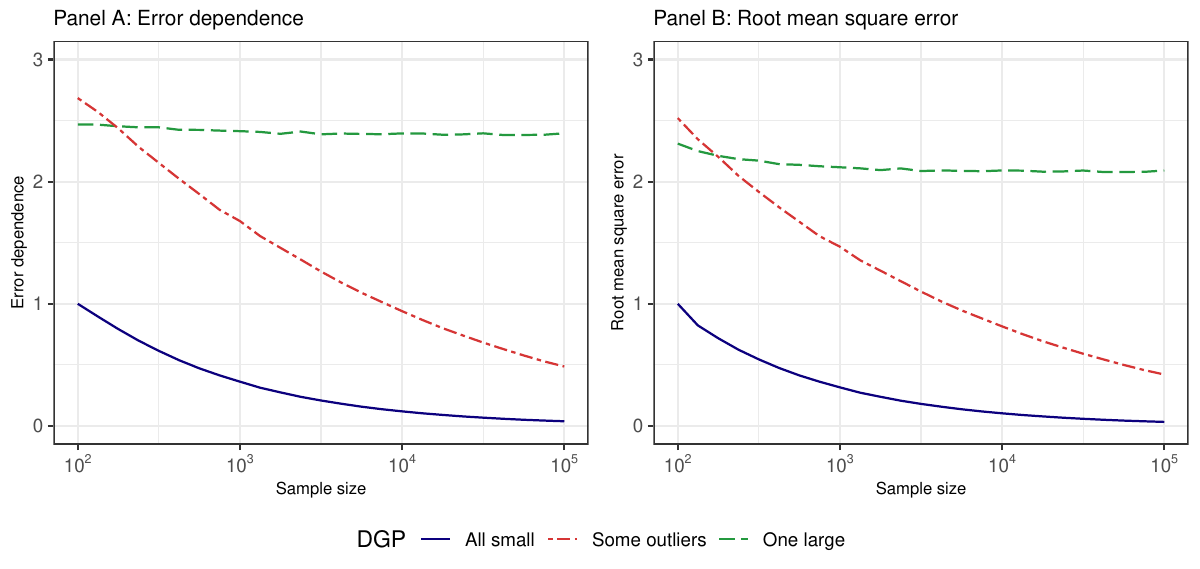}
	\caption{Root of the sum of overall error dependence (A) and root mean square error of the H{\'a}jek estimator (B) under general equilibrium effects. See the note of Figure~\ref{fig:local-interference} for additional details.}\label{fig:market-interference}
\end{figure}

We need that $\Vars[\big]{ p^*(\Z) } \to 0$ for limited specification error dependence (Condition~\mainref{cond:limited-error-dependence}) to hold.
The variance of the market price will approach zero under the first two versions of the data generating process (``All small'' and ``Some outliers''), by the law of large numbers.
However, the experimental design will induce variability in the price under the third version (``One large''), no matter the sample size.
That is, $\Vars[\big]{ p^*(\Z) }$ is asymptotically bounded away from zero.

The results are presented in Figure~\ref{fig:market-interference}.
Panel A shows that the error dependence diminishes quickly under the first version (``All small'').
It also diminishes under the second version (``Some outliers''), but at a slower rate.
The error dependence appears to be unrelated to the sample size under the third version (``One large''), indicating approximately the same level of market price variability for all sample sizes in this case.
Condition~\mainref{cond:limited-error-dependence} does not hold under the third version, and the consistency result in the current paper do not apply.

Panel B largely mirrors the results in the first panel.
The mean square error diminishes under the first two versions of the data generating process (``All small'' and ``Some outliers''), and the rate is faster under the first version.
While the precision initially improves slightly under the third version, it quickly flattens out, and there appears to be no improvements after about $n = 1,000$.
The estimator does not appear to be consistent under the third version.

\begin{table}
\caption{Additional simulation results}\label{tab:simulation-results}
\centering
\renewcommand{\arraystretch}{1.2}
\vspace{0.08in}
\resizebox{\textwidth}{!}{ %
\setlength{\tabcolsep}{12pt}
\hspace{-0.5in}
\begin{tabular}{l p{0.05in} rrrr p{0.1in} rrrr}
\multicolumn{11}{l}{Panel A: Normalized overall error dependence} \\ \cmidrule[\heavyrulewidth]{1-6}
 & & \multicolumn{1}{l}{$n=10^2$} & \multicolumn{1}{l}{$n=10^3$} & \multicolumn{1}{l}{$n=10^4$} & \multicolumn{1}{l}{$n=10^5$} \\ \cmidrule[\lightrulewidth]{1-6}
Local: Correct        &   & 0.000 & 0.000 & 0.000 & 0.000 \\ 
Local: Misspecified   &   & 1.000 & 0.453 & 0.233 & 0.117 \\ 
Herd: Rarely          &   & 1.000 & 0.688 & 0.018 & 0.000 \\ 
Herd: Infrequently    &   & 1.079 & 1.034 & 0.354 & 0.000 \\ 
Herd: Frequently      &   & 1.138 & 1.226 & 1.227 & 1.227 \\ 
Market: All small     &   & 1.000 & 0.363 & 0.120 & 0.039 \\ 
Market: Some outliers &   & 2.684 & 1.677 & 0.940 & 0.487 \\ 
Market: One large     &   & 2.468 & 2.413 & 2.395 & 2.396\\
\cmidrule[\heavyrulewidth]{1-6} \\[0.6em]
\multicolumn{11}{l}{Panel B: Normalized root mean square error} \\ \toprule
& & \multicolumn{4}{c}{H\'ajek estimator} & & \multicolumn{4}{c}{Horvitz--Thompson estimator} \\ \cmidrule{3-6} \cmidrule{8-11}
 & & \multicolumn{1}{l}{$n=10^2$} & \multicolumn{1}{l}{$n=10^3$} & \multicolumn{1}{l}{$n=10^4$} & \multicolumn{1}{l}{$n=10^5$} & & \multicolumn{1}{l}{$n=10^2$} & \multicolumn{1}{l}{$n=10^3$} & \multicolumn{1}{l}{$n=10^4$} & \multicolumn{1}{l}{$n=10^5$} \\ \midrule
Local: Correct        &   & 0.865 & 0.233 & 0.075 & 0.024 &   &  6.884 &  2.184 &  0.691 &  0.219 \\ 
Local: Misspecified   &   & 1.000 & 0.300 & 0.107 & 0.039 &   &  6.228 &  1.861 &  0.527 &  0.137 \\ 
Herd: Rarely          &   & 1.000 & 0.337 & 0.038 & 0.012 &   &  2.800 &  1.030 &  0.348 &  0.110 \\ 
Herd: Infrequently    &   & 1.079 & 0.501 & 0.149 & 0.012 &   &  2.709 &  1.013 &  0.366 &  0.110 \\ 
Herd: Frequently      &   & 1.125 & 0.585 & 0.495 & 0.485 &   &  2.633 &  1.030 &  0.619 &  0.520 \\ 
Market: All small     &   & 1.000 & 0.316 & 0.104 & 0.033 &   & 17.113 &  5.397 &  1.717 &  0.545 \\ 
Market: Some outliers &   & 2.520 & 1.467 & 0.815 & 0.421 &   & 17.404 &  5.605 &  1.896 &  0.687 \\ 
Market: One large     &   & 2.311 & 2.118 & 2.091 & 2.091 &   & 17.355 &  5.824 &  2.701 &  2.160\\
\bottomrule \\[0.6em]
\multicolumn{11}{l}{Panel C: Squared bias as share of mean square error (\%)} \\ \toprule
& & \multicolumn{4}{c}{H\'ajek estimator} & & \multicolumn{4}{c}{Horvitz--Thompson estimator} \\ \cmidrule{3-6} \cmidrule{8-11}
 & & \multicolumn{1}{l}{$n=10^2$} & \multicolumn{1}{l}{$n=10^3$} & \multicolumn{1}{l}{$n=10^4$} & \multicolumn{1}{l}{$n=10^5$} & & \multicolumn{1}{l}{$n=10^2$} & \multicolumn{1}{l}{$n=10^3$} & \multicolumn{1}{l}{$n=10^4$} & \multicolumn{1}{l}{$n=10^5$} \\ \midrule
Local: Correct        &   &  0.0 &  0.0 &  0.0 &  0.0 &   & 0.0 & 0.0 & 0.0 & 0.0 \\ 
Local: Misspecified   &   &  0.3 &  0.1 &  0.0 &  0.0 &   & 0.0 & 0.0 & 0.0 & 0.0 \\ 
Herd: Rarely          &   & 69.0 & 15.1 &  0.0 &  0.0 &   & 0.0 & 0.0 & 0.0 & 0.0 \\ 
Herd: Infrequently    &   & 72.6 & 24.4 &  0.7 &  0.0 &   & 0.0 & 0.0 & 0.0 & 0.0 \\ 
Herd: Frequently      &   & 74.9 & 28.7 &  4.1 &  0.4 &   & 0.0 & 0.0 & 0.0 & 0.0 \\ 
Market: All small     &   &  5.5 &  0.5 &  0.1 &  0.0 &   & 0.0 & 0.0 & 0.0 & 0.0 \\ 
Market: Some outliers &   &  1.2 &  0.0 &  0.0 &  0.0 &   & 0.0 & 0.0 & 0.0 & 0.0 \\ 
Market: One large     &   &  1.4 &  0.0 &  0.0 &  0.0 &   & 0.0 & 0.0 & 0.0 & 0.0\\
\bottomrule \multicolumn{11}{l}{\footnotesize Note: Each cell presents the result from 50,000 draws from the experimental design.} \\
\end{tabular}
}
\end{table}

\clearpage

\section{Estimand without conditioning}\label{sec:other-estimands}

\citet{Saevje2021Average} investigate estimation of average treatment effects under interference.
They distinguish between two types of effects.
The effect of primary interest in \citet{Saevje2021Average} is the \emph{expected average treatment effect}, or EATE, defined as
\begin{equation}
	\frac{1}{n} \sumin \bracket[\Big]{\Es[\big]{\poi{1; \Z_{-i}}} - \Es[\big]{\poi{0; \Z_{-i}}}},
\end{equation}
where $\poi{z; \z_{-i}}$ denotes the potential outcome for unit $i$ when $i$'s treatment assignment is $z$ and the assignment of all other units is $\z_{-i}$.
This estimand captures the effect of changing a unit's own treatment while holding the treatments of all other units fixed, averaged over all units and marginalized over the experimental design.
In this sense, it captures the ``direct'' treatment effect.

\citet{Saevje2021Average} compare the EATE estimand with an estimand described by \citet{Hudgens2008}.
This alternative estimand is often referred to as the direct treatment effect, but \citet{Saevje2021Average} refer to it as the \emph{average distributional shift effect} (ADSE) to differentiate it with EATE.
The ADSE estimand is defined as
\begin{equation}
	\frac{1}{n} \sumin \bracket[\Big]{\Es[\big]{\poi{1; \Z_{-i}} \given \Zei = 1} - \Es[\big]{\poi{0; \Z_{-i}} \given \Zei = 0}}.
\end{equation}

The difference is that EATE marginalizes $\poi{z; \Z_{-i}}$ over the unconditional distribution of $\Z_{-i}$, while ADSE marginalizes $\poi{z; \Z_{-i}}$ over the conditional distribution of $\Z_{-i}$ given $\Zei = z$.
Because the two conditional distributions, corresponding to $\Zei = 1$ and $\Zei = 0$, potentially are different, the ADSE estimand can capture both the direct of effect of $\Zei$ on the outcome of unit $i$ and effects due to the distributional shift of the marginalization.

For example, consider a sample with two units, for which $\poi{z_i; z_j} = z_j$ for both units.
That is, a unit's outcome does not depend on its own treatment but it does depend on the other unit's treatment.
The unit-level direct treatment effects are all zero here:
\begin{equation}
	\po{1}{1; z_2} - \po{1}{0; z_2} = z_2 - z_2 = 0
	\qquad\qquad
	\po{2}{1; z_1} - \po{2}{0; z_1} = z_1 - z_1 = 0.
\end{equation}
This means that the EATE estimand is zero, no matter the experimental design.
However, the ADSE will not always be zero, because it could capture a distributional shift.
For example, if the experimental design is such that $\Ze{1} = 1 - \Ze{2}$, so the two units have opposite treatments, the ADSE estimands is
\begin{equation}
	\frac{1}{n} \sumin \bracket[\Big]{\Es[\big]{\poi{1; 0} \given \Zei = 1} - \Es[\big]{\poi{0; 1} \given \Zei = 0}}
	=
	\frac{1}{n} \sumin \bracket[\big]{0 - 1}
	=
	-1.
\end{equation}

The expected exposure effect, as defined in the current paper, is
\begin{equation}
	\frac{1}{n} \sumin \bracket[\Big]{\Es[\big]{\poi{\Z} \given \exoi = \exea} - \Es[\big]{\poi{\Z} \given \exoi = \exeb}}.
\end{equation}
As this estimand is marginalizing over conditional distributions, it is closer to ADSE than to EATE.
Indeed, with the exposure mapping $\exmiz = \zei$, expected exposure effect is the same as ADSE.

A relevant question in this setting is whether it is possible to define an estimand corresponding to EATE for exposure effects.
The answer is that, in most cases, it will not be possible.
The issue with an EATE-type exposure effect is that units' exposures typically cannot be independently manipulated.
This is an inherent limitation due to the definition of the exposures, rather than a consequence of the experimental design.

To see this, note that in \citet{Saevje2021Average}, the EATE estimand asks ``What is the (expected, average) effect of changing a unit's treatment, holding all other units' treatments fixed?''
This question makes sense because it is generally possible to change a unit's treatment without changing other units' treatment, even if that configuration of treatments is not in the support of the experimental design.
An EATE-type exposure effect would ask ``What is the (expected, average) effect of changing a unit's exposure, holding all other units' exposures (or treatments) fixed?''
The issue is that we can essentially never change a unit's exposure while holding all other units' exposures/treatments fixed.

To illustrate this, consider a network setting with two exposures of interest: when no neighbors are treated ($D_i = 0$) and when at least one neighbor is treated ($D_i = 1$).
The effect of interest is the contrast of the potential outcomes produced by these two exposures.
Now, consider a group of three units ($A$, $B$, $C$), where there are two edges, $(A, B)$ and $(B, C)$, and there are no other units.
That is, the network is $A \leftrightarrow B \leftrightarrow C$.
To define an EATE-type exposure effect in this setting, we would need to compare the situation when unit $A$ is assigned exposure $D_A = 0$ to the situation when it is assigned $D_A = 1$, holding the exposures for $B$ and $C$ fixed.
When $D_A = 0$, we know that unit $B$ is assigned to the control treatment (because $A$ does not have any treated neighbors), which means that unit $C$ also must be assigned to $D_C = 0$ (because it doesn't have any treated neighbors either).
We would now want to consider a setting where we change unit $A$'s exposure while holding the exposures of all other units the same.
But this is not possible.
The only way to change unit $A$'s exposure to $D_A = 1$ is to change unit $B$'s treatment assignment from control to active treatment, ensuring that unit $A$ has a treated neighbor.
As an inescapable consequence, unit $C$ will then also have a treated neighbor, so it will also switch to exposure $D_C = 1$.
In this setting, it is inherently impossible to manipulate the units exposures independently.

The problem in this example is that we cannot separate the effect of having a neighbor that is treated from the effect of having a neighbor that in turn has a neighbor that is treated.
This is similar to the problem that arises for the standard ADSE estimand, but the problem is more fundamental here.
In the setting considered by \citet{Saevje2021Average}, where the implicit exposure mapping is $\exmiz = \zei$, we can almost always independently manipulate the treatments, at least in principle.
That is, we can imagine changing the treatment assigned to one unit without changing any other units' treatment (even if that type of change is not realized by the experimental design).
The problem is therefore entirely driven by dependence in treatment assignments introduced by the experimental design.
Under a Bernoulli design, the ADSE and EATE estimands coincide when the exposure mapping is $\exmiz = \zei$.
The problem with exposure effects more generally is due to dependencies introduced by the nature of the exposures; the exposures are by construction not independently manipulable.
In the example above, there is no way to assign the treatments (i.e., there exists no $\z$) so that unit $A$ has a treated neighbor while unit $C$ does not; units $A$ and $C$ have partially linked exposures by construction.

It is possible to define EATE-type exposure effects for specific types of exposures.
For example, an EATE-type estimand can be defined if the exposures at least in principle are independently manipulable.
One such case is when $\exmiz = \zei$.
More generally, if no two exposures depend on the same treatment assignment $\zei$, they are independently manipulable.
An example when this is the case is when $\exmiz = \ze{\rho(i)}$ for some permutation $\rho(i)$ of the unit indices.
Another example is if the sample is partitioned into clusters and one ``focal'' unit is selected within each cluster so that its exposure only depends on the treatments within its own cluster.
If the analysis is restricted to only the focal units, their exposures will be independently manipulable.
Another possibility is to extend treatment to include interventions to the interference structure itself, which could be used to isolate units from spillover effects.
For example, the question above could be reformulated to ask what the effect is of having at least one treated neighbor versus being completely isolated (i.e., having no neighbors at all).

Slightly more generally, an EATE-type effect can often be defined if the exposure mapping is such that the exposure is given by the treatment assignments of a subset of the units.
Specifically, let $\mathbf{m}_i \in \bset^n$ be a binary vector of length $n$ acting as a masking vector, meaning that for some $\z \in \bset^n$, the $j$th coordinate of the Hadamard product $\z \odot \mathbf{m}_i$ is equal to $\zej$ if and only if $m_i = 1$.
If the exposure mappings are such that $\exmi{\z} = \exea$ if and only if $\z \odot \mathbf{m}_i = \mathbf{r}_{i,\exea}$ and $\exmi{\z} = \exeb$ if and only if $\z \odot \mathbf{m}_i = \mathbf{r}_{i,\exeb}$, where $\mathbf{m}_i$ is common for both exposures, then an EATE-type exposure effect for exposures $\exea$ and $\exeb$ can be defined.
This is done by fixing the coordinates of $\Z$ corresponding to coordinates with value one in $\mathbf{m}_i$ so that $\Z \odot \mathbf{m}_i$ is equal to either $\mathbf{r}_{i,\exea}$ or $\mathbf{r}_{i,\exeb}$.
The coordinates of $\Z$ corresponding to coordinates with value zero in $\mathbf{m}_i$ are marginalized over using their marginal distributions given the experimental design, without regard for other coordinates of $\Z$.
This estimand would avoid the distributional shift artifact that affects the expected exposure effect.
It is possible to slightly extend this effect to setting where $\exmi{\z} = \exea$ if and only if $\z \odot \mathbf{m}_i \in \mathcal{R}_{i,\exea}$ for some set $\mathcal{R}_{i,\exea} \subset \bset^n$.
This would be done by repeating the exercise above for each element in $\mathcal{R}_{i,\exea}$, and then marginalize of these elements using an approriate distribution.
Precise estimation of these EATE-type exposure effects would require additional conditions in line with those in \citet{Saevje2021Average}; it is beyond the scope of this paper to describe and investigation those conditions.

\section{Precision of other estimators}\label{sec:other-estimators}

\subsection{The H{\'a}jek estimator}

The Horvitz--Thompson estimator is rarely a good choice in practice because of its instability in small samples.
The analysis of the Horvitz--Thompson estimator serves as a foundation on which we can build understanding about the behavior of other estimators.
This section investigates the behavior of common refinements of the Horvitz--Thompson estimator, which experimenters often prefer over the original estimator.

The first refinement accounts for the realized number of units assigned to the exposures of interest.
The H{\'a}jek estimator \citep{Hajek1971} does this by dividing each term in the estimator with the sum of the reciprocals of the assignment probabilities for the units assigned to the exposure, rather than dividing by $n$.
The change can absorb some of the variability in the estimator introduced by randomness in the number of units assigned to each exposure.
The ratio structure introduces bias, but the bias is generally small enough to still grant improvements in mean square error sense.
The denominator can generally be shown to be well-behaved, so the estimator's limited behavior can be linked to the Horvitz--Thompson estimator through linearization.

\begin{definition}
The \emph{H{\'a}jek estimator} for expected exposure effect $\efd$ is:
\begin{equation}
	\esthad = \paren[\Bigg]{\sumin \frac{\exiia \ooi}{\preia} \bigg/ \sumin \frac{\exiia}{\preia}} - \paren[\Bigg]{\sumin \frac{\exiib \ooi}{\preib} \bigg/ \sumin \frac{\exiib}{\preib}}.
\end{equation}
\end{definition}

\newcommand{\prophaconsistency}{%
Provided that Conditions~\mainref{cond:bounded-pos},~\mainref{cond:positivity},~\mainref{cond:limited-design-dependence}~and~\mainref{cond:limited-error-dependence} hold for exposures $\exea$ and $\exeb$, the H{\'a}jek estimator is consistent for the expected exposure effect and converges at the rate
\begin{equation}
	\esthad - \efd
	= \bigOp[\big]{n^{-0.5} + \ddep{\exea}^{0.5} + \ddep{\exeb}^{0.5} + \edep{\exea}^{0.5} + \edep{\exeb}^{0.5} + \udep{\exea}^{0.5} + \udep{\exeb}^{0.5}}.
\end{equation}
}

\begin{proposition}\label{prop:ha-consistency}
\prophaconsistency
\end{proposition}

Experimenters often use estimators that implicitly adjust for the assignment probabilities.
One such example is the difference-in-means estimator.
This estimator can be shown to coincide with the H{\'a}jek estimator whenever the assignment probabilities are the same for all units: $\preie = \preje$ for all $i,j\in\Sample$.
Proposition~\ref{prop:ha-consistency} thus implies that the difference-in-means estimator can be used in similar situations also under misspecification.
Experimenters should, however, not blindly use the difference-in-means estimator for exposure effects because the exposure mappings may not induce equal assignment probabilities on the exposures even if they are equal for the nominal treatments.
Another estimator coinciding with the H{\'a}jek estimator is the ordinary least squares (OLS) estimator.
The unweighted version requires equal assignment probabilities just like the difference-in-means estimator, but a weighted OLS estimator is equivalent to the H{\'a}jek estimator also with unequal assignment probabilities.

\subsection{The difference estimator}

A disadvantage of both the Horvitz--Thompson and H{\'a}jek estimators is their inability to take advantage of auxiliary information.
A modification of the Horvitz--Thompson estimator allows us to incorporate such information.
The idea is that information beside the observed potential outcomes themselves might allow us to predict the potential outcomes we do not observe.
If this prediction is sufficiently good, the predicted outcomes can be used to offset chance imbalances introduced by the randomization.
\citet{Saerndal1992Model} call it the difference estimator in a sampling setting, and the name will be used here as well.

\begin{definition}
The \emph{difference estimator} for the expected exposure effect $\efd$ is
\begin{equation}
	\estded = \avgin \bracket[\big]{\poexestia - \poexestib} + \avgin \frac{\paren{\exiia - \exiib} \bracket[\big]{\ooi - \poexesti{\exoi}}}{\prei{\exoi}},
\end{equation}
where $\poexestie$ is a prediction of unit $i$'s potential outcome when assigned to $\exe \in \exall$.
\end{definition}

The definition of the estimator reveals the idea that motivate its use.
The first term is simply the average difference in predicted potential outcomes.
If the predictions are of high quality, this term will be an accurate estimator of the exposure effect.
The issue is that the predictions may have systematic errors.
The second term is included to ensure unbiasedness.
If the predictions are of low quality, this term will compensate for the errors in the first term, and it ensures that the estimator performs well in expectation.
The estimator bears a resemblance in this regard to the class of doubly robust estimators used in observational studies when the assignment mechanism is unknown \citep[see, e.g.,][]{Robins2001Comment}.

The properties of the difference estimator depend on the way the predictions are constructed.
In particular, the estimator can be shown to retain the advantageous properties of the Horvitz--Thompson estimator if the predictions are external to the study.
External here means that they do not depend on the treatment assignment.
As the only randomness under consideration stems from the assignment mechanism, independence between $\poexestie$ and $\Z$ implies that the predictions are non-random.
The probability space can extended to accommodate random predictions if one wants to account for the consequences of external variability.
Such variability could affect the rate of convergence if predictions are sufficiently dependent between units, but it is otherwise inconsequential to the results.

\begin{definition}
A prediction $\poexestie$ is \emph{external} if it is independent of the treatment assignment vector: $\poexestie \indep \Z$.
\end{definition}

\begin{definition}\label{def:average-prediction-dependence}
The \emph{average prediction dependence} for exposure $\exe \in \exall$ is
\begin{equation}
	\pdepe = \frac{1}{n^2} \sumin \sum_{j \neq i} \abs[\big]{\Cov[\big]{\poexestie, \poexestje}}.
\end{equation}
\end{definition}

An alternative to focusing on the dependence between predictions is to consider their convergence.
In particular, the average prediction dependence can be bounded by
\begin{equation}
	\pdepe = \frac{1}{n^2} \sumin \sum_{j \neq i} \abs[\big]{\Cov[\big]{\poexestie, \poexestje}} \leq \paren[\bigg]{\frac{1}{n}\sumin\sqrt{\Var{\poexestie}}}^2,
\end{equation}
which shows that Condition~\ref{cond:limited-prediction-dependence} is satisfied if the predictions on average converges in mean square.

\begin{condition}[Prediction moments]\label{cond:prediction-moments}
$\E[\big]{\abs{\poexestie}^2} \leq \constpred < \infty$ for all $i \in \Sample$ and $\exe \in \exall$.
\end{condition}

\begin{condition}\label{cond:limited-prediction-dependence}
An exposure $\exe \in \exall$ satisfies \emph{limited prediction dependence} if $\pdepe = \littleO{1}$.
\end{condition}

\newcommand{\propdeunbiasedness}{%
Provided that Condition~\mainref{cond:positivity} holds and that the predictions are external, the difference estimator is unbiased for the expected exposure effect: $\E{\estded} = \efd$.
}

\begin{proposition}\label{prop:de-unbiasedness}
\propdeunbiasedness
\end{proposition}

\newcommand{\propdeconsistency}{%
Provided that Conditions~\mainref{cond:bounded-pos},~\mainref{cond:positivity},~\mainref{cond:limited-design-dependence},~\mainref{cond:limited-error-dependence},~\ref{cond:prediction-moments}~and~\ref{cond:limited-prediction-dependence} hold and that the predictions are external, the difference estimator is consistent for the expected exposure effect and converges at the rate
\begin{equation}
	\estded - \efd
	= \bigOp[\big]{n^{-0.5} + \ddep{\exea}^{0.5} + \ddep{\exeb}^{0.5} + \edep{\exea}^{0.5} + \edep{\exeb}^{0.5} + \udep{\exea}^{0.5} + \udep{\exeb}^{0.5} + \pdepa^{0.5} + \pdepb^{0.5}}.
\end{equation}
}

\begin{proposition}\label{prop:de-consistency}
\propdeconsistency
\end{proposition}

The difference estimator seemingly provides advantages at no cost.
Good predictions of the potential outcomes confer improvements in finite samples, but the estimator has the same behavior as the Horvitz--Thompson estimator in large samples.
The no-cost advantages are superficial.
The mean square error may increase when the predictions are poor, so investigators should use the difference estimator only when the predictions are expected to be of reasonably high quality.

However, the quality of the predictions is less of a concern than their construction.
Covariate information can be used to make the predictions, but the assigned exposures and the observed outcomes can generally not be used because it would induce dependence between the predictions and $\Z$.
More precisely, if $\covi$ denotes a vector of covariates describing characteristics of unit $i$, we can form the predictions as $\poexestie = f\paren{\exe, \covi}$ for some function $f$.
The function $f$ can, however, not be constructed using $\paren{\oo{1}, \oo{2}, \dotsc, \oo{n}}$ or $\paren{\exo{1}, \exo{2}, \dotsc, \exo{n}}$.
This illustrates that the construction of $f$ truly needs to be external to treatment assignment when used for the predictions in the difference estimator.
This severely limits its applicability.
Split-sample or leave-one-out approaches \citep[see, e.g.,][]{Williams1961Generating} that often are used to solve the issue cannot be used here because the misspecification may induce dependence between subsamples that otherwise appear isolated.

\subsection{The generalized regression estimator}

An estimator facilitating dependence between the predictions of the potential outcomes and the treatment assignments is inspired by the generalized regression estimator commonly used in the sampling literature.
The estimator has received recent attention in the causal inference literature as well \citep[see, e.g.,][]{Lin2013Agnostic,Middleton2018Unified}.

The estimator uses a linear working model for the relationship between the potential outcomes and the covariates.
The working model is used to construct the predictions.
Generally, $\poexestie = \covi^\tran \coefe$ for some vector of coefficients $\coefe$ indexed by $\exe \in \exall$, so different coefficients are used for different exposures.
No assumptions are made about the validity of the model, but the quality of the predictions are related to how well the model can approximate the potential outcomes.
It remains to pick the coefficients $\coefe$.
The generalized regression estimator allows for dependence between the coefficients and the treatment assignments, so the coefficients can be estimated in the sample.
For example, we may pick them as the minimizing solution to $\sumin \exiie\bracket{\ooi - \covi^\tran \coefe}^2$ as is often done in applications.
But other choices exist, and the estimator is largely agnostic about how the coefficients were constructed.

\begin{definition}
The \emph{generalized regression estimator} for expected exposure effect is
\begin{equation}
	\estgrd = \avgin \covi^\tran \bracket[\big]{\coefesta - \coefestb} + \avgin \frac{\paren{\exiia - \exiib} \bracket[\big]{\ooi - \covi^\tran \coefest{\exoi}}}{\prei{\exoi}},
\end{equation}
where $\coefesta$ and $\coefestb$ are two random vectors of the same dimensions as $\covi$.
\end{definition}

The conventional approach to investigate the properties of the generalized regression estimator is to assume that the vector of coefficients constructed in the sample convergences to some fixed vector asymptotically.
This ensures that the dependence between units' predictions is small in large samples, which provides consistency.
The assumption can be weaken to only require that the magnitude of the vector of coefficients is asymptotically bounded, thereby bypassing the need of assuming a well-defined limit.

\begin{condition}[Bounded covariates]\label{cond:bounded-covs}
$\covi \in \mathcal{X}$ for some bounded $\mathcal{X} \subset \Reals^p$.
\end{condition}

\begin{condition}[Asymptotically bounded regression coefficients]\label{cond:bounded-coefs}
$\E[\big]{\norm{\coefeste}} = \bigO{1}$.
\end{condition}

\newcommand{\propgrconsistency}{%
Provided that Conditions~\mainref{cond:bounded-pos},~\mainref{cond:positivity},~\mainref{cond:limited-design-dependence},~\mainref{cond:limited-error-dependence},~\ref{cond:bounded-covs}~and~\ref{cond:bounded-coefs} hold, the generalized regression estimator is consistent for the expected exposure effect and converges at the rate
\begin{equation}
	\estgrd - \efd
	= \bigOp[\big]{n^{-0.5} + \ddep{\exea}^{0.5} + \ddep{\exeb}^{0.5} + \edep{\exea}^{0.5} + \edep{\exeb}^{0.5} + \udep{\exea}^{0.5} + \udep{\exeb}^{0.5}}.
\end{equation}
}

\begin{proposition}\label{prop:gr-consistency}
\propgrconsistency
\end{proposition}

\section{Lack of positivity}\label{sec:lack-positivity}

Positivity conditions are often seen as innocuous in experiments because the experimenter controls the design and can ensure their validity.
But this is rarely the case when estimating exposure effects.
Exposure mappings tend to be complex, and it may not be feasible to construct a design that would induce a desired distribution over the exposures.
Experimenters will instead settle for heuristic choices for the design at the treatment level, and this could induce violations of Condition~\mainref{cond:positivity}.

The positivity condition can fail in two ways.
The first is when it is fundamentally impossible for a unit to be assigned a particular exposure.
For example, in the experiment in Bogot\'a described in Section~\mainref{sec:illustration} in the main paper, a non-hot spot street without any neighboring hot spot streets cannot be assigned to an exposure requiring that at least one neighboring hot spot street receives intensive policing.
This may be formalized by saying that there is an exposure $\exe \in \exall$ for which no $\z \in \zall$ exists with $\exmiz = \exe$.

The consequences of such a failure are more than just statistical.
If it is nonsensical to talk about some collection of units being assigned to a certain exposure, it is nonsensical to consider exposure effects that include those units in its average.
Unless the experimenter is comfortable stipulating a metaphysical model allowing extrapolation to unrealizable potential outcomes, the only solution is to exclude such units from the average.
The result may be that the number of units included in the analysis is fewer than the length of $\z$, but this is not an issue other than for efficiency.
In the following discussion, it will be assumed that such exclusions have been made if necessary.
That is, if the aim is to estimate the effect of exposures $\exea$ and $\exeb$, then $\setb{\exea, \exeb} \subseteq \setb{\exmiz : \z \in \zall}$ for all units $i \in \Sample$.

The second way the positivity condition can fail is through the design; assignments $\z \in \zall$ exist so that $\exmiz = \exe$, but the design is such that $\preie = 0$.
Statistical issues are the only sequelae in this case, which all have cures.
Two situations must be considered.
The first is when the assignment probability for some exposure is exactly zero, $\preie = 0$.
The second is when the probability approaches zero asymptotically.
Both are problematic, but they have different solutions.

Superficially, the first situation appears most acute.
There are two issues to consider.
The first is that the definition of $\poexie$ conditions on the measure-zero event $\exoi = \exe$, rendering the definition ambiguous.
To address this, extend the definition as follows:
\begin{equation}
	\poexie =
		\begin{cases}
		\E{\poi{\Z} \given \exoi = \exe} & \text{if } \preie > 0,
		\\
		\paren[\Big]{ \sum_{\z \in \zall} \indicator[\big]{\exmiz = \exe} }^{-1} \paren[\Big]{ \sum_{\z \in \zall} \indicator[\big]{\exmiz = \exe} \poiz }  & \text{if } \preie = 0.
		\end{cases}
\end{equation}
That is, if $\preie = 0$, then $\poexie$ is the arithmetic mean of all potential outcomes for which the corresponding $\z$ maps to $\exmiz = \exe$.
The second concern is that the Horvitz--Thompson estimator now involves division by zero, rendering it ill-defined.
However, for each term of the estimator with a zero denominator, the numerator is also zero with probability one.
Hence, a straightforward solution is to define $0/0$ as zero, and that is the solution that will be used here.
But to ensure that the estimator behaves well, the proportion of units with zero assignment probability must be small.
To capture this, let $\zpreie = \indicator{\preie = 0}$ denote whether unit $i$ has a zero assignment probability, and let
\begin{equation}
	\zpravge = \frac{1}{n} \sumin \zpreie,
\end{equation}
be the proportion of such units in the sample.

As for assignment probabilities that approaches zero, we must consider the rate at which they do so.
The following norm-like quantity captures the average rate of convergence towards zero:
\begin{equation}
	\zprmome = \bracket[\bigg]{\frac{1}{n} \sumin \frac{1 - \zpreie}{\bracket{\preie}^{\exppimom} + \zpreie}}^{1/\exppimom}.
\end{equation}

The quantities $\zpravge$ and $\zprmome$ allow us to weaken the positivity assumption in a controllable way.
In particular, Condition~\mainref{cond:positivity} is the same as $\zpravge = 0$ and $\lim_{\exppimom \to \infty} \zprmome \leq \constpi < \infty$.
The following proposition shows that neither part is necessary for consistency.
But the weakening comes at the cost of potentially slower convergence rates.
This is captured by a strengthening of the definition of the design dependence, namely
\begin{equation}
	\ddepexte = \bracket[\Bigg]{\frac{1}{n^2}\sumin \sum_{j \neq i} \abs[\big]{\Cov[\big]{\exiie, \exije}}^{q}}^{1/q}.
\end{equation}
This extended definition collapses to the definition in Condition~\mainref{cond:limited-design-dependence} when $q = 1$, but generally $\ddepe = \littleO[\big]{\ddepexte}$ when $q > 1$.
Thus, $\ddepexte = \littleO{1}$ when $q > 1$ is stronger than the original condition.
Extending the results in \citet{Delevoye2020Consistency} to a setting with interference, this additional machinery admits a proof of consistency without positivity.

\newcommand{\propzhtconsistency}{%
Suppose Conditions~\mainref{cond:bounded-pos}~and~\mainref{cond:limited-error-dependence} hold, and that
\begin{equation}
\zprmome \leq \constpi < \infty,
\qquad
\zpravge  = \littleO{1}
\qquadand
\ddepext[\big]{\exe}{\exppimom / (\exppimom - 2)} = \littleO{1},
\end{equation}
for $\exe \in \setb{\exea, \exeb}$ and some $\exppimom > 2$.
The Horvitz--Thompson estimator is then consistent for the expected exposure effect and converges at the rate
\begin{equation}
	\estd - \efd = \bigOp[\big]{n^{-0.5} + \zpravga + \zpravgb + \ddepextstda^{0.5} + \ddepextstdb^{0.5} + \edep{\exea}^{0.5} + \edep{\exeb}^{0.5} + \udep{\exea}^{0.5} + \udep{\exeb}^{0.5}},
\end{equation}
where $\ddepextstde$ is short-hand for $\ddepext[\big]{\exe}{\exppimom / (\exppimom - 2)}$.%
}

\begin{proposition}\label{prop:zht-consistency}
\propzhtconsistency
\end{proposition}

The proposition states that we can achieve consistency under misspecification even if positivity does not hold as long as the dependence between exposures, as captured by $\ddepexte$, is sufficiently weak.
While experimenters still should try to ensure that their designs and exposure mappings satisfy positivity, Proposition~\ref{prop:zht-consistency} provides some reassurance that the results are not automatically invalidated in the case they are not perfectly successful and small violations to positivity occur.

\section{Limiting distribution}\label{sec:limiting-distribution}

The limiting distribution of the estimator is less tractable than its limit.
Two approaches are explored in this section.
The first is to tie the limiting behavior of the estimator under misspecification to its behavior when the exposures are correctly specified.
Experimenters might find this result useful because familiar results for correctly specified exposures are directly extended to a setting with misspecification.
The second approach is a direct proof of asymptotic normality using Stein's method.
Both approaches require considerably stronger assumptions than those needed for consistency.

\subsection{Connection to correctly specified exposures}\label{sec:limiting-distribution-v1}

A situation where progress can be made is when the specification errors are very small relative to the sample size.
The following aggregated measure of the misspecification is a strengthening of the error dependence measures used for the consistency results.

\begin{definition}\label{def:total-error-dependence}
The \emph{average total error dependence} is
\begin{equation}
	\tdepe = \frac{1}{n^2} \sumin \sumjn \braces[\big]{\E{\seri \serj \given \exoi = \exoj = \exe}}^+,
\end{equation}
where $\braces{x}^+ = \maxf{0, x}$ denotes the positive part of $x$.
\end{definition}

The definition is a strengthening of Definition~\mainref{def:error-dependence} in two ways.
First, the average total error dependence considers only the positive parts of the specification errors, while the dependence measures in Definition~\mainref{def:error-dependence} includes the negative terms as well.
It is possible that the dependence between some units' errors is negative, and this will have a compensatory effect in Definition~\mainref{def:error-dependence}, making the average smaller.
Definition~\ref{def:total-error-dependence} ignores any such compensatory effects.
Second, unlike the previous dependence measures, the average total error dependence includes the diagonal elements of the double sum, $\E{\seri^2 \given \exoi = \exe}$, which generally will be larger than the off-diagonal elements.

To connect the behavior of the estimator under misspecification to its behavior under correctly specified exposures, we will impose a general condition on the experimental design and exposure mappings to ensure that the estimator is well-behaved when the exposures are correctly specified.

\begin{condition}[Design regularity]\label{cond:design-regularity}
The experimental design and exposure mappings are such that for any collection of bounded functions $f_i : \exall \to \Reals$ for $i \in \setb{1, \dotsc, n}$, the following random variable converges in distribution to some random variable $\limdist$:
\begin{equation}
\convseq \sumin \bracket[\bigg]{ \frac{\braces{\exiia - \preia} f_i(\exea) }{n \preia} - \frac{\braces{\exiib - \preib} f_i(\exeb)}{n \preib}} \darrow \limdist,
\end{equation}
where $\convseq$ is some sequence ensures that the overall sequence is on the order of a constant.
\end{condition}

Under Condition~\mainref{cond:bounded-pos} (bounded potential outcomes), the Horvitz--Thompson estimator of the exposure effect between exposures $\exea$ and $\exeb$ has the limiting distribution $\limdist$ when the exposures are correctly specified if and only if Condition~\ref{cond:design-regularity} holds.
Indeed, Condition~\ref{cond:design-regularity} is simply that statement written in a somewhat more general form.
That is, if an experimenter is in a setting where they believe the estimator has some limiting distribution $\limdist$ if the exposures are correctly specified, then Condition~\ref{cond:design-regularity} holds.

\newcommand{\proplimitdist}{%
Suppose Conditions~\mainref{cond:bounded-pos}~and~\mainref{cond:positivity} hold, and that Condition~\ref{cond:design-regularity} holds for some sequence $\convseq$ and some random variable $\limdist$.
If the exposures are misspecified, but the misspecification is sufficiently weak so that $\tdepa + \tdepb = \littleO{\convseq^{-2}}$, then the limiting distribution of the Horvitz--Thompson estimator is $\limdist$:
\begin{equation}
	\convseq \bracket[\big]{\estd - \efd} \darrow \limdist.
\end{equation}
}

\begin{proposition}\label{prop:limit-dist}
\proplimitdist
\end{proposition}

The proposition states that the limiting distribution of the estimator is unchanged if the specification errors are small.
However, observe that the condition $\tdepa + \tdepb = \littleO{\convseq^{-2}}$ used in Proposition~\ref{prop:limit-dist} is considerably stronger than Condition~\mainref{cond:limited-error-dependence} used for the consistency results.
Indeed, while the previous conditions allowed for considerable misspecification also in large samples, this stronger condition is saying that any misspecification is negligible asymptotically relative to the variability induced by the randomization of treatments.
Hence, the applicability of Proposition~\ref{prop:limit-dist} is limited, and experimenters should show caution before using the proposition to motivate any inferential statements.
But, in the few situations in which Proposition~\ref{prop:limit-dist} is applicable, experimenters may find the following special case of the proposition particularly useful.

\begin{corollary}
If the Horvitz--Thompson estimator would be root-$n$ consistent and asymptotically normal when the exposures are correctly specified, in the sense that $\convseq$ is of order $\sqrt{n}$ and $\limdist$ is a standard normal deviate, then the estimator is root-$n$ consistent and asymptotically normal also under misspecified exposures provided that $\tdepa + \tdepb = \littleO{n^{-1}}$.
\end{corollary}

\subsection{Asymptotic normality using Stein's method}\label{sec:limiting-distribution-v2}

Stein's method has been extensively used in the recent literature on interference.
An early example is \citet{Aronow2017Estimating}.
The result used here is due to \citet{Ross2011Fundamentals}.
This result has previously been used by \citet{Chin2019Central}.
\citet{Ogburn2022Causal} provide a stronger result, but the result in \citet{Ross2011Fundamentals} is used here due to its simplicity.
There is nothing in the current application of Stein's method that is new compared to the previous results in the interference literature, but it might be of interest to note that these well-known results extend also to situations with misspecification.

Define the dependency neighborhood for a unit $i \in \Sample$ to be the smallest subset $\set{N}_i \subset \Sample$ such that unit $i$'s exposure $\exoi$ and specification error $\seri$ are independent of the exposures and specification errors of all units not in $\set{N}_i$:
\begin{equation}
	\paren{ \exoi, \seri} \indep \braces[\big]{ \paren{ \exoj, \serj} : j \in \Sample \setminus \set{N}_i}.
\end{equation}
Let $d_{\max} = \max_{i\in\Sample} \abs{\set{N}_i}$ be the size of the largest dependency neighborhood in the sample.

\begin{condition}[Restricted dependency neighborhoods]\label{cond:rest-depend-neighbor}
	$d_{\max} = \littleO{n^{1/4}}$.
\end{condition}

Condition~\ref{cond:rest-depend-neighbor} is considerably stronger than the conditions used in the main paper.
It restricts the design to only exhibit local dependence with respect to the exposures, which rules out many common designs and some exposure mappings.
It also restrict the specification errors to only be locally dependent.
That is, the type of global error dependence that was found to be unproblematic in the main paper is not allowed here.
While the condition is strong, it might nevertheless be found to be reasonable in some experiments.
In that case, the sampling distribution of the estimator will be approximately normal in large samples, as the following proposition demonstrates.

\newcommand{\proplimitdiststein}{%
Provided that Conditions~\mainref{cond:bounded-pos},~\mainref{cond:positivity}~and~\ref{cond:rest-depend-neighbor} hold for exposures $\exea$ and $\exeb$ and that the estimator is not superefficient, ensuring that the variance does not diminish at a faster than parametric rate $\Var{\estd} = \bigOmega{n^{-1}}$, the Horvitz--Thompson estimator is asymptotically normal:
\begin{equation}
	\frac{\estd - \efd}{\sqrt{\Var{\estd}}} \darrow Z,
\end{equation}
where $Z$ denotes a standard normal deviate.
}

\begin{proposition}\label{prop:limit-dist-stein}
\proplimitdiststein
\end{proposition}

\section{Variance estimation}\label{sec:variance-estimation}

\subsection{Variance estimation under correctly specified exposures}

Variance estimation for exposure effect estimators is challenging because the variance consists of pair-wise products of potential outcomes, and some of those outcomes are not simultaneously observable.
The issue is not unique to exposure effects, but exposure mappings tend to induce complex distributions on the exposures, which exacerbates the problem.

The solution suggested by \citet{Aronow2017Estimating} is to use Young's inequality for products to bound the unobservable parts of the variance expression.
To better understand this idea, let $\preijot = \Pr{\exoi = \exeone, \exoj = \exetwo}$ be the joint probability of unit $i$ and $j$'s exposures.
If $\preijot = 0$, then the potential outcomes $\poexio$ and $\poexjt$ are never observed simultaneously, which will complicate variance estimation because the variance depends on the product $\poexio \poexjt$.
\citet{Aronow2017Estimating} tackle these products by using the bound
\begin{equation}
	\poexio \poexjt
	\leq \frac{\bracket{\poexio}^2 + \bracket{\poexjt}^2}{2}.
\end{equation}
After having applied the bound to all problematic terms in the variance of the point estimator, they arrive at the estimator
\begin{multline}
	\EstVarHTAS
	= \frac{1}{n^2} \sumin \sumjn \paren{\exiia - \exiib} \paren{\exija - \exijb} \varwei{\exoi, \exoj} \ooi \ooj
	\\
	+ \frac{1}{n^2} \sumin \sumjn \bracket[\bigg]{\frac{\exiia}{\preia} + \frac{\exiib}{\preib}} \bracket[\big]{\zpreij{\exoi, \exea} + \zpreij{\exoi, \exeb}} \ooi^2,
\end{multline}
where,
\begin{equation}
	\varwei{\exeone, \exetwo}
	= \frac{\preijot - \preio \prejt}{\preijot \preio \prejt + \zpreijot}
	\quadand
	\zpreijot
	= \indicator{\preijot = 0}.
\end{equation}

\citet{Aronow2017Estimating} show that this variance estimator is conservative in expectation when exposures are correctly specified.
However, what does not appear to be fully appreciated in the literature is that the bound on the problematic products may make the estimator excessively conservative.
In fact, unless the assumption of correctly specified exposures is complemented with
\begin{equation}
	\frac{1}{n^2} \sumin \sum_{j \neq i} \bracket[\big]{\zpreijaa + \zpreijbb + \zpreijab} = \bigO[\big]{n^{-1}},
\end{equation}
the normalized variance estimator $n \EstVarHTAS$ generally diverges to infinity.
I will not offer a solution to this problem.
The remark instead serves as an illustration of the difficulty of variance estimation for complex exposure effects.
It also provides insights about the mechanics of the estimator, which will aid our understanding of its behavior under misspecification.

\subsection{Variance estimation under misspecification}

The analysis of the variance estimator by \citet{Aronow2017Estimating} does not hold when the exposures are misspecified.
The expectation of the variance estimator could both increase and decrease under misspecification relative to when the exposures are correctly specified.
In some situations, the decrease is sizable, and the estimator may become anti-conservative, providing an unjustly optimistic estimate of the precision of the point estimator.
The following proposition exactly characterizes the bias of the variance estimator under misspecification.

\newcommand{\propbiasvarianceestimator}{%
Provided that Conditions~\mainref{cond:bounded-pos}~and~\mainref{cond:positivity} hold for exposures $\exea$ and $\exeb$, the bias of the variance estimator described by \citet{Aronow2017Estimating} is
\begin{multline}
	\E[\Big]{\EstVarHTAS} - \Var[\big]{\estd} = \varbiasf{1}{\exea, \exeb} + \varbiasf{2}{\exea, \exeb} + \varbiasf{2}{\exeb, \exea} + \varbiasf{3}{\exea, \exeb} + \varbiasf{3}{\exeb, \exea}
	\\
	+ 2\varbiasf{4}{\exea, \exeb} - \varbiasf{4}{\exea, \exea} - \varbiasf{4}{\exeb, \exeb},
\end{multline}
where
\begin{align}
	\varbiasf{1}{\exeone, \exetwo} &= \frac{1}{n^2} \sumin \bracket[\big]{\poexio - \poexit}^2,
	\\
	\varbiasf{2}{\exeone, \exetwo} &= \frac{1}{2n^2} \sumin \sum_{j \neq i} \paren[\Big]{\zpreij{\exeone, \exeone} \bracket[\big]{\poexio + \poexjo}^2 + \zpreij{\exeone, \exetwo} \bracket[\big]{\poexio - \poexit}^2},
	\\
	\varbiasf{3}{\exeone, \exetwo} &= \frac{1}{n^2} \sumin \sum_{j \neq i} \bracket[\big]{\zpreij{\exeone, \exeone} + \zpreij{\exeone, \exetwo}} \Var{\seri \given \exoi = \exeone},
	\\
	\varbiasf{4}{\exeone, \exetwo} &= \frac{1}{n^2} \sumin \sum_{j \neq i} \bracket[\big]{1 - \zpreij{\exeone, \exetwo}} \Big[\poexio\serefji{\exetwo, \exeone} + \poexit\serefij{\exeone, \exetwo}
	\\
		&\qquad\qquad\qquad\qquad + \serefij{\exeone, \exetwo}\serefji{\exetwo, \exeone} + \Cov{\seruij, \seruji \given \exoi = \exeone, \exoj = \exetwo}\Big].
\end{align}%
}

\begin{proposition}\label{prop:bias-variance-estimator}
\propbiasvarianceestimator
\end{proposition}

The terms on the right-hand side capture aspects that introduce bias of the variance estimator.
These bias terms help us understand when variance estimation is possible under misspecification.
The term $\varbiasf{1}{\exea, \exeb}$ stems from what \citet{Holland1986Statistics} describes as the fundamental problem of causal inference, namely that a unit cannot simultaneously be assigned to two different treatments.
The joint distribution of the potential outcomes affects the variance, but the distribution can only be estimated if both potential outcomes are observed simultaneously.
Such simultaneous observations would require simultaneous assignment of two different treatments to the same unit, but this is not possible.
The first term captures the bias arising from our inability to estimate this aspect of the potential outcomes.
The issue is not unique to the current setting, and similar bias terms arise for most causal inference problems in finite populations, including when the exposures are correctly specified.

As already noted, the distribution of the exposures is often complex, and the joint exposure probabilities may be zero for a considerable number of pairs of units.
This issue is similar to the first source of bias, but it is now introduced by the design rather than being fundamental.
The use of Young's inequality to tackle this problem introduces bias, and this bias is captured by the terms $\varbiasf{2}{\exea, \exeb}$ and $\varbiasf{2}{\exeb, \exea}$ in Proposition~\ref{prop:bias-variance-estimator}.
These biases arise also when the exposures are correctly specified.

At this point, we have replicated the result in \citet{Aronow2017Estimating}.
In particular, the remaining terms are zero by construction when the exposures are correctly specified.
In other words, the bias of the variance estimator when the exposures are correctly specified is
\begin{equation}
	\E[\Big]{\EstVarHTAS} - \Var[\big]{\estd} = \varbiasf{1}{\exea, \exeb} + \varbiasf{2}{\exea, \exeb} + \varbiasf{2}{\exeb, \exea},
\end{equation}
which the same result as in \citet{Aronow2017Estimating}, but presented in a different form.
The fact that these three terms are non-negative by construction confirms that the estimator is conservative when the exposures are correctly specified.
What Proposition~\ref{prop:bias-variance-estimator} shows is that we need to consider additional sources of bias when the exposures are misspecified.

The remaining biases stem from two sources.
The first also arises from the use of Young's inequality in the construction of the estimator.
If the probability that two units are simultaneously assigned to a certain combination of exposures is zero, then the corresponding specification errors cannot interact, and the errors do not affect the variance of the point estimator.
This means that we do not need to adjust for any dependence between such errors, because we know it is zero.
However, these pairs of units are exactly the ones for which we apply the bound on the unobserved product of potential outcomes to ensure conservativeness under correctly specified exposures.
This inadvertently leads to that the variance estimator is affected by the magnitude of the corresponding errors.
The terms $\varbiasf{3}{\exea}$ and $\varbiasf{3}{\exeb}$ capture this part of the bias.
Like the previous terms, these terms are non-negative by construction.

The terms of real concern are the last three: $\varbiasf{4}{\exea, \exeb}$, $\varbiasf{4}{\exea, \exea}$, and $\varbiasf{4}{\exeb, \exeb}$.
These capture the bias introduced by our inability to estimate the dependence in the specification errors.
Unlike the previous terms, the signs of the terms are unknown, so they may introduce negative bias.
The consequence is that we could systematically underestimate the variance when the specification errors are negatively dependent, and our inferences would then be anti-conservative.

The problem has no immediate solution, but some progress can be made.
Similar to the approach taken in Section~\ref{sec:limiting-distribution}, if the specification errors can be assumed to be negligible relative to the sample size, the terms $\varbiasf{4}{\exeone, \exetwo}$ are negligible relative to the other terms, and the variance estimator is ensured to be asymptotically conservative.

An alternative approach is to incorporate more information about the structure of the interference in the variance estimator.
In particular, the anti-conservative behavior of the estimator stems from negative interactions of errors in $\varbiasf{4}{\exeone, \exetwo}$.
One may remove such interactions by setting $\zpreij{\exeone, \exetwo} = 1$ for the corresponding pairs of units, even if $\preij{\exeone, \exetwo} > 0$ holds.
This will move the corresponding terms from $\varbiasf{4}{\exeone, \exetwo}$, where negative interactions are possible, to $\varbiasf{3}{\exe}$, where no interactions exist.
It may be hard to discern whether the interaction between two specific units' errors is negative or positive.
A conservative approach is to set $\zpreij{\exeone, \exetwo} = 1$ for all pairs of units where an interaction of any type is suspected.

As an example, consider when units only interfere with each other within known disjoint groups, which is an interference structure often referred to as ``partial interference.''
One may here set $\zpreij{\exeone, \exetwo} = 1$ if either $\preij{\exeone, \exetwo}$ is zero or if unit $i$ and $j$ belong to the same group.
A redefinition of $\zpreij{\exeone, \exetwo}$ along these lines would ensure that $\varbiasf{4}{\exeone, \exetwo}=0$, so the variance estimator remains conservative.
Of course, such knowledge about the interference would allow for the definition of exposure mappings that are correctly specified, which would remove any concerns about misspecification.
However, experimenters may want to keep the main exposure mapping simple to facilitate interpretation.
They can then proceed with misspecified exposures for point estimation, and use the more intricate information about the interference structure only when estimating variance.

A third approach is a combination of the previous two.
One may set $\zpreij{\exeone, \exetwo} = 1$ for pairs of units where negative interaction are suspected to be particularly large.
Even if one does not capture all pairs with negative interactions, the variance estimator will still be conservative as long as the missed interactions are small.
In other words, setting $\zpreij{\exeone, \exetwo} = 1$ for the terms deemed most problematic makes the assumption that the remaining errors are small more reasonable.
It should also be noted that the other bias terms tend to be large and positive, and they therefore provide considerable leeway in the case $\varbiasf{4}{\exeone, \exetwo}$ is negative.

The expectation of the variance estimator should be seen as a rough representation of its behavior more generally.
The precision of the variance estimator will be poor if the joint exposure probabilities are small even if they are never exactly zero.
Experimenters should be aware that the variance estimator could be imprecise even when the bias is positive, particularly when the number of exposures is large.
This concern is not specific to misspecification.

\section{Proofs}

\subsection{Miscellaneous lemmas}

\begin{lemma}\label{lem:var-avg-bound}
For any $N$ random variables $X_1, \dotsc, X_N$,
\begin{equation}
	\Var[\big]{X_1 + \dotsb + X_N}
	\leq \paren[\Big]{\sqrt{\Var{X_1}} + \dotsb + \sqrt{\Var{X_N}}}^2.
\end{equation}
\end{lemma}

\begin{proof}
Write the variance of the sum as a double sum of covariances:
\begin{equation}
	\Var[\Bigg]{\sum_{i=1}^N X_i}
	= \sum_{i=1}^N \sum_{j=1}^N \Cov{X_i, X_j}.
\end{equation}
Separate the covariances using the Cauchy--Schwarz inequality and reorder:
\begin{equation}
	\sum_{i=1}^N \sum_{j=1}^N \Cov{X_i, X_j}
	\leq \sum_{i=1}^N \sum_{j=1}^N\sqrt{\Var{X_i}\Var{X_j}}
	= \paren[\Bigg]{\sum_{i=1}^N\sqrt{\Var{X_i}}}^2.
	\tag*{\qedhere}
\end{equation}
\end{proof}

\begin{lemma}\label{lem:var-sum-bound}
For any $N$ random variables $X_1, \dotsc, X_N$,
\begin{equation}
	\Var{X_1 + \dotsb + X_N}
	\leq N \Var{X_1} + \dotsb + N \Var{X_N}.
\end{equation}
\end{lemma}

\begin{proof}
Apply Lemma~\ref{lem:var-avg-bound} to get
\begin{equation}
	\Var[\Bigg]{\sum_{i=1}^N X_i}
	\leq \paren[\Bigg]{\sum_{i=1}^N\sqrt{\Var{X_i}}}^2
	= N^2 \paren[\Bigg]{\frac{1}{N} \sum_{i=1}^N\sqrt{\Var{X_i}}}^2.
\end{equation}
The square is a convex function, so Jensen's inequality gives
\begin{equation}
	N^2 \paren[\Bigg]{\frac{1}{N} \sum_{i=1}^N\sqrt{\Var{X_i}}}^2
	\leq N^2 \paren[\Bigg]{\frac{1}{N} \sum_{i=1}^N \Var{X_i}}
	= N \sum_{i=1}^N \Var{X_i}.
	\tag*{\qedhere}
\end{equation}
\end{proof}

\begin{lemma}\label{lem:outcome-bounds}
If Condition~\mainref{cond:bounded-pos} holds, then for all $i\in\Sample$ and $\exe \in \exall$,
\begin{enumerate}[label=\roman*.,ref=\ref*{lem:outcome-bounds}\roman*]
\item $\abs{\poexie} \leq \constpo$,\label{lem:outcome-bounds-po}
\item $\E{\ooi^2 \given \exoi = \exe} \leq \constpo^2$,\label{lem:outcome-bounds-ysq}
\item $\abs[\big]{\E{\seri\serj \given \exoi = \exoj = \exe}} \leq 4\constpo^2$.\label{lem:outcome-bounds-errors}
\end{enumerate}
\end{lemma}

\begin{proof}
Consider each statement in turn:
\begin{enumerate}[label=\roman*.]
\item
Recall the definition $\poexie = \E{\poi{\Z} \given \exoi = \exe}$.
Condition~\mainref{cond:bounded-pos} gives
\begin{equation}
	\abs{\poexie}
	= \abs[\big]{\E{\poi{\Z} \given \exoi = \exe}}
	\leq \E[\big]{\abs{\poi{\Z}} \given \exoi = \exe}
	\leq \E{\constpo \given \exoi = \exe}
	= \constpo.
\end{equation}

\item
Condition~\mainref{cond:bounded-pos} ensures that $\abs{\ooi} = \abs{\poi{\Z}} \leq \constpo$ with probability one.
It follows that
\begin{equation}
	\E[\big]{\ooi^2 \given \exoi = \exe}
	= \E[\big]{\paren[\big]{\poi{\Z}}^2 \given \exoi = \exe}
	\leq \constpo^2.
\end{equation}

\item
Recall that $\seri = \ooi - \poexi{\exoi}$, so we can write
\begin{equation}
	\seri\serj
	= \ooi \ooj - \ooi \poexj{\exoj} - \ooj \poexi{\exoi} + \poexi{\exoi} \poexj{\exoj}.
\end{equation}
Condition~\mainref{cond:bounded-pos} ensures that $\abs{\ooi} \leq \constpo$ and $\abs{\poexi{\exoi}} \leq \constpo$ with probability one.
It follows that
\begin{equation}
	\abs{\seri\serj}
	\leq 4\constpo^2,
\end{equation}
with probability one.
Hence,
\begin{equation}
	\abs[\big]{\E{\seri\serj \given \exoi = \exoj = \exe}}
	\leq \E{\abs{\seri\serj} \given \exoi = \exoj = \exe}
	\leq 4\constpo^2.
	\tag*{\qedhere}
\end{equation}
\end{enumerate}
\end{proof}

\newcommand{\lemunbiasedness}{%
Provided that Condition~\mainref{cond:positivity} holds for exposures $\exea$ and $\exeb$, the Horvitz--Thompson estimator is unbiased for the expected exposure effect: $\E{\estd} = \efd$.
}

\begin{lemma}\label{lem:unbiasedness}
\lemunbiasedness
\end{lemma}

\begin{proof}
For any exposure $\exe \in \exall$ satisfying Condition~\mainref{cond:positivity}, we have
\begin{equation}
	\avgin \frac{\E{\exiie\ooi}}{\preie}
	= \avgin \frac{\preie\E{\ooi \given \exoi = \exe}}{\preie}
	= \avgin \frac{\preie\poexie}{\preie}
	= \avgin \poexie.
\end{equation}
It follows that
\begin{equation}
	\E{\estd}
	= \avgin \frac{\E{\exiia\ooi}}{\preia} - \avgin \frac{\E{\exiib\ooi}}{\preib}
	= \avgin \poexia - \avgin \poexib
	= \efd.
	\tag*{\qedhere}
\end{equation}
\end{proof}

\subsection{Proposition~\mainref{prop:variance-bound}: Variance bound}

\begin{refproposition}{\mainref{prop:variance-bound}}
\propvariancebound
\end{refproposition}

\begin{proof}
Lemma~\ref{lem:initial-var-bound} gives
\begin{equation}
	\Var[\big]{\estd}
	\leq \frac{8 \constpo^2 \constpi}{n}
	+ 4 \constpo^2 \constpi^2 \paren{\ddep{\exea} + \ddep{\exeb}}
	+ \frac{4}{n^2} \sum_{\exe \in \braces{\exea, \exeb}} \sumin \sum_{j \neq i} \frac{\Cov[\big]{\exiie \seri, \exije \serj}}{\preie \preje}.
\end{equation}
Using Lemma~\ref{lem:bound-error-cov}, we can bound the third term for $\exe \in \braces{\exea, \exeb}$ by
\begin{multline}
	\frac{4}{n^2} \sumin \sum_{j \neq i} \frac{\Cov[\big]{\exiie \seri, \exije \serj}}{\preie \preje}
	\leq \frac{16 \constpo^2 \constpi^2}{n^2} \sumin \sum_{j \neq i} \abs[\big]{\Cov{\exiie, \exije}}
	\\
	+ \frac{4}{n^2} \sumin \sum_{j \neq i} \E[\big]{\seri\serj \given \exoi = \exoj = \exe}.
\end{multline}
Using the definitions of $\ddepe$ and $\odepe$, we can write this as
\begin{equation}
	\frac{4}{n^2} \sumin \sum_{j \neq i} \frac{\Cov[\big]{\exiie \seri, \exije \serj}}{\preie \preje}
	\leq 16 \constpo^2 \constpi^2 \ddepe + 4 \odepe.
\end{equation}
Collecting terms completes the proof.
\end{proof}

\begin{lemma}\label{lem:initial-var-bound}
Provided that Conditions~\mainref{cond:bounded-pos}~and~\mainref{cond:positivity} hold for exposures $\exea$ and $\exeb$, the variance of the Horvitz--Thompson estimator is upper bounded by
\begin{equation}
	\Var[\big]{\estd}
	\leq \frac{8 \constpo^2 \constpi}{n}
	+ 4 \constpo^2 \constpi^2 \paren{\ddep{\exea} + \ddep{\exeb}}
	+ \frac{4}{n^2} \sum_{\exe \in \braces{\exea, \exeb}} \sumin \sum_{j \neq i} \frac{\Cov[\big]{\exiie \seri, \exije \serj}}{\preie \preje}.
\end{equation}
\end{lemma}

\begin{proof}
Recall the definition of the estimator:
\begin{equation}
	\estd
	= \avgin \frac{\exiia}{\preia} \ooi - \avgin \frac{\exiib}{\preib} \ooi.
\end{equation}
Apply Lemma~\ref{lem:var-sum-bound} to get
\begin{equation}
	\Var[\big]{\estd}
	\leq 2\Var[\bigg]{\avgin \frac{\exiia}{\preia} \ooi} + 2\Var[\bigg]{\avgin \frac{\exiib}{\preib} \ooi}.
\end{equation}
Note that Definition~\mainref{def:specification-error} implies that $\ooi = \poexi{\exoi} + \seri$, so for any exposure $\exe \in \exall$,
\begin{equation}
	\avgin \frac{\exiie}{\preie} \ooi
	= \avgin \frac{\exiie}{\preie} \poexie + \avgin \frac{\exiie}{\preie} \seri.
\end{equation}
Apply Lemma~\ref{lem:var-sum-bound} again to get
\begin{equation}
	2\Var[\bigg]{\avgin \frac{\exiie}{\preie} \ooi}
	\leq 4 \Var[\bigg]{\avgin \frac{\exiie}{\preie} \poexie} + 4 \Var[\bigg]{\avgin \frac{\exiie}{\preie} \seri}.
\end{equation}
We can write the two terms as
\begin{align}
	4 \Var[\bigg]{\avgin \frac{\exiie}{\preie} \poexie}
	&= \frac{4}{n^2} \sumin \bracket[\bigg]{\frac{\poexie}{\preie}}^2 \Var{\exiie}
	+ \frac{4}{n^2} \sumin \sum_{j \neq i} \frac{\poexie \poexje}{\preie \preie} \Cov{\exiie, \exiie},
	\\
	4 \Var[\bigg]{\avgin \frac{\exiie}{\preie} \seri}
	&= \frac{4}{n^2} \sumin \bracket[\bigg]{\frac{1}{\preie}}^2 \Var[\big]{\exiie \seri}
	 + \frac{4}{n^2} \sumin \sum_{j \neq i} \frac{\Cov[\big]{\exiie \seri, \exije \serj}}{\preie \preje},
\end{align}
because $\poexie$ and $\preie$ are not random.

Lemma~\ref{lem:var-prod-decomp} implies that
\begin{equation}
	\frac{4}{n^2} \sumin \bracket[\bigg]{\frac{\poexie}{\preie}}^2 \Var[\big]{\exiie} + \frac{4}{n^2} \sumin \bracket[\bigg]{\frac{1}{\preie}}^2 \Var[\big]{\exiie \seri}
	= \frac{4}{n^2} \sumin \bracket[\bigg]{\frac{1}{\preie}}^2 \Var[\big]{\exiie \ooi}.
\end{equation}
Lemma~\ref{lem:outcome-bounds-ysq} gives
\begin{equation}
	\Var[\big]{\exiie \ooi}
	\leq \E[\big]{\exiie \ooi^2}
	= \preie \E[\big]{\ooi^2 \given \exoi = \exe}
	\leq \preie \constpo^2.
\end{equation}
Together with Condition~\mainref{cond:positivity}, this implies
\begin{equation}
	\frac{4}{n^2} \sumin \bracket[\bigg]{\frac{1}{\preie}}^2 \Var[\big]{\exiie \ooi}
	\leq \frac{4}{n^2} \sumin \constpo^2 \constpi
	= \frac{4 \constpo^2 \constpi}{n}.
\end{equation}

Using Condition~\mainref{cond:positivity} and Lemma~\ref{lem:outcome-bounds-po}, we can bound the following term by
\begin{equation}
	\frac{4}{n^2} \sumin \sum_{j \neq i} \frac{\poexie \poexje}{\preie \preie} \Cov{\exiie, \exiie}
	\leq \frac{4 \constpo^2 \constpi^2}{n^2} \sumin \sum_{j \neq i} \abs[\big]{\Cov{\exiie, \exiie}}
	= 4 \constpo^2 \constpi^2 \ddepe.
\end{equation}
Collecting terms completes the proof.
\end{proof}

\begin{lemma}\label{lem:var-prod-decomp}
For all $i\in\Sample$ and $\exe \in \exall$,
\begin{equation}
	\Var[\big]{\exiie \ooi}
	= \bracket{\poexie}^2 \Var{\exiie} + \Var[\big]{\exiie \seri}.
\end{equation}
\end{lemma}

\begin{proof}
Definition~\mainref{def:specification-error} implies that
\begin{equation}
	\E{\seri \given \exoi = \exe}
	= \E{\ooi - \poexi{\exoi} \given \exoi = \exe}
	= \E{\poi{\Z} \given \exoi = \exe} - \poexie
	= 0.
\end{equation}
It follows that
\begin{equation}
	\E[\big]{\exiie^2 \poexie \seri}
	= \preie \E[\big]{\exiie^2 \poexie \seri \given \exoi = \exe}
	= \poexie \preie \E{\seri \given \exoi = \exe}
	= 0,
\end{equation}
and
\begin{equation}
	\E[\big]{\exiie \seri}
	= \preie \E[\big]{\exiie \seri \given \exoi = \exe}
	= \preie \E{\seri \given \exoi = \exe}
	= 0.
	\label{eq:error-prod-zero}
\end{equation}
We therefore have
\begin{equation}
	\Cov[\big]{\exiie \poexie, \exiie \seri}
	= \E[\big]{\exiie^2 \poexie \seri} - \E[\big]{\exiie \poexie} \E[\big]{\exiie \seri}
	= 0.
\end{equation}
It follows that
\begin{equation}
	\Var[\big]{\exiie \ooi}
	= \Var[\big]{\exiie \poexie + \exiie \seri}
	= \bracket{\poexie}^2 \Var[\big]{\exiie} + \Var[\big]{\exiie \seri}.
	\tag*{\qedhere}
\end{equation}
\end{proof}

\begin{lemma}\label{lem:bound-error-cov}
If Conditions~\mainref{cond:bounded-pos}~and~\mainref{cond:positivity} hold, then for all $i, j \in \Sample$ and $\exe \in \exall$,
\begin{equation}
	\frac{\Cov[\big]{\exiie \seri, \exije \serj}}{\preie \preje}
	\leq 4 \constpo^2 \constpi^2 \abs[\big]{\Cov{\exiie, \exije}} + \E[\big]{\seri\serj \given \exoi = \exoj = \exe}.
\end{equation}
\end{lemma}

\begin{proof}
Equation~\eqref{eq:error-prod-zero} in the proof of Lemma~\ref{lem:var-prod-decomp} showed that $\E{\exiie \seri} = 0$, so we have
\begin{equation}
	\Cov[\big]{\exiie \seri, \exije \serj}
	= \E[\big]{\exiie \exije \seri \serj}
\end{equation}
The rest of the proof proceeds separately for the cases $\preij{\exe, \exe} > 0$ and $\preij{\exe, \exe} = 0$.

Starting with the case $\preij{\exe, \exe} > 0$, note that we can write the expression to be bounded as
\begin{align}
	\frac{\Cov[\big]{\exiie \seri, \exije \serj}}{\preie \preje}
	&= \frac{\E[\big]{\exiie \exije \serj \seri}}{\preie \preje}
	= \frac{\Pr{\exoi = \exoj = \exe}}{\preie\preje} \E[\big]{\seri\serj \given \exoi = \exoj = \exe}
	\\
	&= \frac{\Pr{\exoi = \exoj = \exe} - \preie\preje + \preie\preje}{\preie\preje} \E[\big]{\seri\serj \given \exoi = \exoj = \exe}
	\\
	&= \frac{\Cov[\big]{\exiie, \exije}}{\preie\preje} \E[\big]{\seri\serj \given \exoi = \exoj = \exe} + \E[\big]{\seri\serj \given \exoi = \exoj = \exe},
\end{align}
because $\preij{\exe, \exe} > 0$ ensures that $\preie > 0$ and $\preje > 0$.
Consider the first term through the lens of Condition~\mainref{cond:positivity} and Lemma~\ref{lem:outcome-bounds-errors}:
\begin{equation}
	\frac{\Cov[\big]{\exiie, \exije}}{\preie\preje} \E[\big]{\seri\serj \given \exoi = \exoj = \exe}
	\leq 4\constpo^2\constpi^2\abs[\big]{\Cov{\exiie, \exije}}.
\end{equation}
Thus, when $\preij{\exe, \exe} > 0$,
\begin{equation}
	\frac{\Cov[\big]{\exiie \seri, \exije \serj}}{\preie \preje}
	\leq 4 \constpo^2 \constpi^2 \abs[\big]{\Cov{\exiie, \exije}} + \E[\big]{\seri\serj \given \exoi = \exoj = \exe}.
\end{equation}

Continuing with the case $\preij{\exe, \exe} = 0$, note that $\exiie \exije$ is zero with probability one in this case, so
\begin{equation}
	\frac{\Cov[\big]{\exiie \seri, \exije \serj}}{\preie \preje}
	= \frac{\E[\big]{\exiie \exije \serj \seri}}{\preie \preje}
	= 0.
\end{equation}
Recall that we defined $\E[\big]{\seri\serj \given \exoi = \exoj = \exe} = 0$ in the case $\preij{\exe, \exe} = 0$.
It follows that
\begin{equation}
	\frac{\Cov[\big]{\exiie \seri, \exije \serj}}{\preie \preje}
	\leq 4 \constpo^2 \constpi^2 \abs[\big]{\Cov{\exiie, \exije}} + \E[\big]{\seri\serj \given \exoi = \exoj = \exe},
\end{equation}
also when $\preij{\exe, \exe} = 0$.
\end{proof}

\subsection{Proposition~\ref*{prop:consistency-decompose}: Bound based on error decomposition}\label{sec:consistency-decompose-proof}

\begin{refproposition}{\ref*{prop:consistency-decompose}}
\propconsistencydecomp
\end{refproposition}

\begin{proof}
The proof presumes that $\edep{\exea}$, $\edep{\exeb}$, $\udep{\exea}$ and $\udep{\exeb}$ are non-negative.
If they are not, the contravening quantities can be set to zero, and the proof would apply.

Consider the root mean square error:
\begin{multline}
	\sqrt{\E[\bigg]{\paren[\Big]{\estd - \efd}^2}} = \sqrt{\Var[\big]{\estd}}
	\\
	\leq \sqrt{8 \constpo^2 \constpi / n + 20 \constpo^2 \constpi^2 \paren{\ddep{\exea} + \ddep{\exeb}} + 4 \paren{\edep{\exea} + \edep{\exeb} + \udep{\exea} + \udep{\exeb}}}
\end{multline}
where the equality follows from Lemma~\ref{lem:unbiasedness}, and the inequality from Lemma~\ref{lem:variance-bound-decompose}.
Concavity of the square root gives
\begin{equation}
	\sqrt{\Var[\big]{\estd}}
	\leq 3 \constpo \constpi^{0.5} / n^{0.5} + 5 \constpo \constpi \paren[\big]{\ddep{\exea}^{0.5} + \ddep{\exeb}^{0.5}} + 2 \paren[\big]{\edep{\exea}^{0.5} + \edep{\exeb}^{0.5} + \udep{\exea}^{0.5} + \udep{\exeb}^{0.5}},
\end{equation}
which gives the rate of convergence in the $L^2$-norm:
\begin{equation}
	\sqrt{\E[\bigg]{\paren[\Big]{\estd - \efd}^2}} =
	\bigO[\big]{n^{-0.5} + \ddep{\exea}^{0.5} + \ddep{\exeb}^{0.5} + \edep{\exea}^{0.5} + \edep{\exeb}^{0.5} + \udep{\exea}^{0.5} + \udep{\exeb}^{0.5}}.
\end{equation}
In turn, this gives the rate of convergence in probability by Markov's inequality.
\end{proof}

\begin{lemma}\label{lem:variance-bound-decompose}
Provided that Conditions~\mainref{cond:bounded-pos}~and~\mainref{cond:positivity} hold for exposures $\exea$ and $\exeb$, the variance of the Horvitz--Thompson estimator is upper bounded by
\begin{equation}
	\Var[\big]{\estd} \leq 8 \constpo^2 \constpi / n + 20 \constpo^2 \constpi^2 \paren{\ddep{\exea} + \ddep{\exeb}} + 4 \paren{\edep{\exea} + \edep{\exeb} + \udep{\exea} + \udep{\exeb}}.
\end{equation}
\end{lemma}

\begin{proof}
Following the proof of Proposition~\mainref{prop:variance-bound}, use Lemma~\ref{lem:initial-var-bound} to write
\begin{equation}
	\Var[\big]{\estd}
	\leq \frac{8 \constpo^2 \constpi}{n}
	+ 4 \constpo^2 \constpi^2 \braces{\ddep{\exea} + \ddep{\exeb}}
	+ \frac{4}{n^2} \sum_{\exe \in \braces{\exea, \exeb}} \sumin \sum_{j \neq i} \frac{\Cov[\big]{\exiie \seri, \exije \serj}}{\preie \preje}.
\end{equation}
Using Lemma~\ref{lem:bound-error-cov-decomp}, we can bound the third term for $\exe \in \braces{\exea, \exeb}$ by
\begin{multline}
	\frac{4}{n^2} \sumin \sum_{j \neq i} \frac{\Cov[\big]{\exiie \seri, \exije \serj}}{\preie \preje}
	\leq \frac{16 \constpo^2 \constpi^2}{n^2} \sumin \sum_{j \neq i} \abs[\big]{\Cov{\exiie, \exije}}
	\\
	+ \frac{4}{n^2} \sumin \sum_{j \neq i} \serefije\serefjie
	+ \frac{4}{n^2} \sumin \sum_{j \neq i} \Cov{\seruij, \seruji \given \exoi = \exoj = \exe}.
\end{multline}
Using the definitions of $\ddepe$, $\edepe$ and $\udepe$, we can write this as
\begin{equation}
	\frac{4}{n^2} \sumin \sum_{j \neq i} \frac{\Cov[\big]{\exiie \seri, \exije \serj}}{\preie \preje}
	\leq 16 \constpo^2 \constpi^2 \ddepe + 4 \edepe + 4 \udepe.
\end{equation}
Collecting terms completes the proof.
\end{proof}

\begin{lemma}\label{lem:bound-error-cov-decomp}
If Conditions~\mainref{cond:bounded-pos}~and~\mainref{cond:positivity} hold, then for all $i, j \in \Sample$ and $\exe \in \exall$,
\begin{multline}
	\frac{\Cov[\big]{\exiie \seri, \exije \serj}}{\preie \preje}
	\leq 4 \constpo^2 \constpi^2 \abs[\big]{\Cov{\exiie, \exije}} + \serefije\serefjie
	\\
	+ \Cov{\seruij, \seruji \given \exoi = \exoj = \exe}.
\end{multline}
\end{lemma}

\begin{proof}
Equation~\eqref{eq:error-prod-zero} in the proof of Lemma~\ref{lem:var-prod-decomp} showed that $\E{\exiie \seri} = 0$, so we have
\begin{equation}
	\Cov[\big]{\exiie \seri, \exije \serj}
	= \E[\big]{\exiie \exije \serj \seri}
\end{equation}
The rest of the proof proceeds separately for the cases $\preij{\exe, \exe} > 0$ and $\preij{\exe, \exe} = 0$.

Starting with the case $\preij{\exe, \exe} > 0$, note that we can write the expression to be bounded as
\begin{align}
	\frac{\Cov[\big]{\exiie \seri, \exije \serj}}{\preie \preje}
	&= \frac{\E[\big]{\exiie \exije \serj \seri}}{\preie \preje}
	= \frac{\Pr{\exoi = \exoj = \exe}}{\preie\preje} \E[\big]{\seri\serj \given \exoi = \exoj = \exe}
	\\
	&= \frac{\Pr{\exoi = \exoj = \exe} - \preie\preje + \preie\preje}{\preie\preje} \E[\big]{\seri\serj \given \exoi = \exoj = \exe}
	\\
	&= \frac{\Cov[\big]{\exiie, \exije}}{\preie\preje} \E[\big]{\seri\serj \given \exoi = \exoj = \exe} + \E[\big]{\seri\serj \given \exoi = \exoj = \exe},
\end{align}
because $\preij{\exe, \exe} > 0$ ensures that $\preie > 0$ and $\preje > 0$.
Consider the first term through the lens of Condition~\mainref{cond:positivity} and Lemma~\ref{lem:outcome-bounds-errors}:
\begin{equation}
	\frac{\Cov[\big]{\exiie, \exije}}{\preie\preje} \E[\big]{\seri\serj \given \exoi = \exoj = \exe}
	\leq 4\constpo^2\constpi^2\abs[\big]{\Cov{\exiie, \exije}}.
\end{equation}

For the second term, let $\sereij = \serefij{\exoi, \exoj}$, so that $\seri = \sereij + \seruij$.
Then,
\begin{multline}
	\E[\big]{\seri\serj \given \exoi = \exoj = \exe}
	= \E[\Big]{\paren[\big]{\sereij + \seruij}\paren[\big]{\sereji + \seruji} \given \exoi = \exoj = \exe}
	\\
	= \E[\big]{\sereij\sereji + \sereij\seruji + \seruij\sereji + \seruij\seruji \given \exoi = \exoj = \exe}.
\end{multline}
Note that $\sereij = \serefij{\exoi, \exoj}$ and $\sereji = \serefij{\exoj, \exoi}$ are constant conditional on $\exoi = \exoj = \exe$, and that
\begin{equation}
	\E{\seruij \given \exoi, \exoj}
	= \E[\big]{\ooi - \porexij{\exoi, \exoj} \given \exoi, \exoj}
	= \E{\ooi \given \exoi, \exoj} - \porexij{\exoi, \exoj}
	= 0.
\end{equation}
It follows that
\begin{equation}
	\E[\big]{\sereij\seruji \given \exoi = \exoj = \exe}
	= \serefije\E[\big]{\seruji \given \exoi = \exoj = \exe}
	= 0,
\end{equation}
and
\begin{multline}
	\E[\big]{\seri\serj \given \exoi = \exoj = \exe}
	= \serefije\serefjie + \E[\big]{\seruij\seruji \given \exoi = \exoj = \exe}
	\\
	= \serefije\serefjie + \Cov{\seruij, \seruji \given \exoi = \exoj = \exe}.
\end{multline}
Taken together we have showed that, when $\preij{\exe, \exe} > 0$,
\begin{multline}
	\frac{\Cov[\big]{\exiie \seri, \exije \serj}}{\preie \preje}
	\leq 4 \constpo^2 \constpi^2 \abs[\big]{\Cov{\exiie, \exije}} + \serefije\serefjie
	\\
	+ \Cov{\seruij, \seruji \given \exoi = \exoj = \exe}.
\end{multline}

Continuing with the case $\preij{\exe, \exe} = 0$, note that $\exiie \exije$ is zero with probability one in this case, so
\begin{equation}
	\frac{\Cov[\big]{\exiie \seri, \exije \serj}}{\preie \preje}
	= \frac{\E[\big]{\exiie \exije \serj \seri}}{\preie \preje}
	= 0.
\end{equation}
Recall that when $\preij{\exe, \exe} = 0$,
\begin{equation}
	\serefije = 0
	\qquadand
	\Cov{\seruij, \seruji \given \exoi = \exoj = \exe} = 0.
\end{equation}
It follows that, when $\preij{\exe, \exe} = 0$,
\begin{multline}
	\frac{\Cov[\big]{\exiie \seri, \exije \serj}}{\preie \preje}
	\leq 4 \constpo^2 \constpi^2 \abs[\big]{\Cov{\exiie, \exije}} + \serefije\serefjie
	\\
	+ \Cov{\seruij, \seruji \given \exoi = \exoj = \exe}.
	\tag*{\qedhere}
\end{multline}
\end{proof}

\subsection{Proposition~\ref*{prop:limit-dist}: Limit distribution}

\begin{refproposition}{\ref*{prop:limit-dist}}
\proplimitdist
\end{refproposition}

\DeclarePairedDelimiterXPP\estorc[1]{\widehat{\bar{\tau}}}{\lparen}{\rparen}{}{#1}
\newcommand{\estorcd}{\estorc{\exea, \exeb}}

\DeclarePairedDelimiterXPP\esterr[1]{\hat{\varepsilon}}{\lparen}{\rparen}{}{#1}
\newcommand{\esterrd}{\esterr{\exea, \exeb}}

\begin{proof}
Note that Definition~\mainref{def:specification-error} allows us to write the observed outcome as $\ooi = \poexi{\exoi} + \seri$, which means that the estimator can be written as
\begin{equation}
	\estd = \estorcd + \esterrd,
\end{equation}
where
\begin{equation}
\estorcd = \frac{1}{n} \sumin \frac{\exiia \poexia}{\preia} - \frac{1}{n} \sumin \frac{\exiib \poexib}{\preib}
\quadand
\esterrd = \frac{1}{n} \sumin \frac{\exiia \seri}{\preia} - \frac{1}{n} \sumin \frac{\exiib \seri}{\preib}.
\end{equation}
The random variable $\estorcd$ is the Horvitz--Thompson estimator we would use if we observed $\poexie$ when $\exoi = \exe$.
That is, if we observed the expected potential outcome without the specification error.
The random variable $\esterrd$ is the Horvitz--Thompson estimator of the treatment effect on the expected specification error.

The decomposition of the estimator allows us to write
\begin{equation}
	\convseq \bracket[\big]{\estd - \efd}
	= \convseq \bracket[\Big]{\estorcd - \efd} + \convseq \esterrd.
\end{equation}
Condition~\mainref{cond:bounded-pos} implies that $\poexia$ and $\poexib$ are bounded for all $i$ and $n$.
Condition~\ref{cond:design-regularity} then implies that
\begin{equation}
\convseq \bracket[\Big]{\estorcd - \efd} \darrow \limdist.
\end{equation}
Lemma~\ref{lem:limit-error-ht} gives
\begin{equation}
	\convseq \esterrd
	= \bigOp[\big]{\convseq \tdepa^{0.5} + \convseq \tdepb^{0.5}}
	= \littleOp{1}.
\end{equation}
As a result, we can apply Slutsky's theorem to get
\begin{equation}
	\convseq \bracket[\big]{\estd - \efd}
	\darrow \limdist.
	\tag*{\qedhere}
\end{equation}
\end{proof}

\begin{lemma}\label{lem:limit-error-ht}
If Condition~\mainref{cond:positivity} holds, then
\begin{equation}
	\esterrd
	= \frac{1}{n} \sumin \frac{\exiia \seri}{\preia} - \frac{1}{n} \sumin \frac{\exiib \seri}{\preib}
	= \bigOp[\big]{\tdepa^{0.5} + \tdepb^{0.5}}.
\end{equation}
\end{lemma}

\begin{proof}
Note that the expectation of $\esterrd$ is zero:
\begin{equation}
	\E[\big]{\esterrd}
	= \frac{1}{n} \sumin \E{\seri \given \exoi = \exea} - \frac{1}{n} \sumin \E{\seri \given \exoi = \exeb}
	= 0,
\end{equation}
because
\begin{equation}
	\E{\seri \given \exoi = \exea}
	= \E[\big]{\ooi - \poexi{\exoi} \given \exoi = \exea}
	= \E{\ooi \given \exoi = \exea} - \poexi{\exea}
	= 0.
\end{equation}
Using Lemma~\ref{lem:var-sum-bound}, the variance can be bounded as
\begin{equation}
	\Var[\big]{\esterrd}
	\leq \frac{2}{n^2} \Var[\bigg]{ \sumin \frac{\exiia \seri}{\preia} } + \frac{2}{n^2} \Var[\bigg]{ \sumin \frac{\exiib \seri}{\preib} }.
\end{equation}
Each of the two terms can be written as
\begin{multline}
	\frac{2}{n^2} \Var[\bigg]{ \sumin \frac{\exiie \seri}{\preie} }
	= \frac{2}{n^2} \sumin \sumjn \frac{\Cov[\big]{ \exiie \seri, \exije \serj }}{\preie \preje}
	\\
	= \frac{2}{n^2} \sumin \sumjn \frac{\Pr{\exoi = \exoj = \exe}}{\preie \preje} \E{ \seri \serj \given \exoi = \exoj = \exe},
\end{multline}
because $\E[\big]{ \exiie \seri } = 0$.
Using Condition~\mainref{cond:positivity}, this may in turn be bounded by
\begin{multline}
	\frac{2}{n^2} \sumin \sumjn \frac{\Pr{\exoi = \exoj = \exe}}{\preie \preje} \E{ \seri \serj \given \exoi = \exoj = \exe}
	\\
	\leq \frac{2}{n^2} \sumin \sumjn \frac{\Pr{\exoi = \exoj = \exe}}{\preie \preje} \bracket[\big]{\E{\seri \serj \given \exoi = \exoj = \exe}}^+
	\\
	\leq \frac{2 \constpi}{n^2} \sumin \sumjn \bracket[\big]{\E{\seri \serj \given \exoi = \exoj = \exe}}^+
	= 2 \constpi \tdepe.
\end{multline}
Hence, the variance is bounded by
\begin{equation}
	\Var[\big]{\esterrd}
	\leq 2 \constpi \paren[\big]{\tdepa + \tdepb}.
\end{equation}
This provides convergence to zero in mean square:
\begin{equation}
	\sqrt{\E[\Big]{\paren[\big]{\esterrd}^2}}
	= \bigO[\big]{\tdepa^{0.5} + \tdepb^{0.5}},
\end{equation}
which in turn provides convergence to zero in probability by Markov's inequality.
\end{proof}

\subsection{Proposition~\ref*{prop:limit-dist-stein}: Asymptotic normality using Stein's method}

\begin{refproposition}{\ref*{prop:limit-dist-stein}}
\proplimitdiststein
\end{refproposition}

\begin{proof}
Let $\delta = \estd - \efd$.
Note that $\Var{\estd} = \Var{\delta}$ because $\efd$ is not random.
Hence,
\begin{equation}
	\frac{\estd - \efd}{\sqrt{\Var{\estd}}} = \frac{\delta}{\sqrt{\Var{\delta}}}.
\end{equation}
Let
\begin{equation}
	\delta_i = \frac{\exiia \ooi}{\preia} - \frac{\exiib \ooi}{\preib} - \braces[\big]{\poexia - \poexib},
\end{equation}
so that $\delta = n^{-1} \sum_{i=1}^n \delta_i$.
Note that the dependency neighborhood for $\delta_i$ is the same as the dependency neighborhood for $\paren{ \exoi, \seri}$, because we can use Definition~\mainref{def:specification-error} to write
\begin{equation}
	\delta_i = \frac{\exiia \poexia}{\preia} + \frac{\exiia \seri}{\preia} - \frac{\exiib \poexib}{\preib} - \frac{\exiib \seri}{\preib} - \braces[\big]{\poexia - \poexib},
\end{equation}
and $\exoi$ and $\seri$ are the only random variables in the expression.
Using the argument in Lemma~\ref{lem:unbiasedness}, one can show that $\E{\delta_i} = 0$.
Lemma~\ref{lem:delta-bounded} implies that $\E{\delta_i^4} < \infty$, and that
\begin{equation}
	\frac{1}{n} \sum_{i=1}^n \E[\big]{\abs{\delta_i}^3} = \bigO{1},
	\qquad\qquad
	\frac{1}{n} \sum_{i=1}^n \E[\big]{\delta_i^4} = \bigO{1}.
\end{equation}
This means that we can apply Theorem 3.6 from \citet{Ross2011Fundamentals} to get
\begin{equation}
	d_W\paren[\bigg]{\frac{\delta}{\sqrt{\Var{\delta}}}, Z} = \bigO[\bigg]{ \frac{d_{\max}^2}{n^2 \bracket{\Var{\delta}}^{3/2}}
	+ \frac{d_{\max}^{3/2} }{n^{3/2} \Var{\delta}} },
\end{equation}
where $d_W$ denotes the Wasserstein distance between the two distributions in its arguments.

Condition~\ref{cond:rest-depend-neighbor} and the fact that $\Var{\delta} = \bigOmega{n^{-1}}$ implies that
\begin{equation}
	d_W\paren[\bigg]{\frac{\delta}{\sqrt{\Var{\delta}}}, Z} = \littleO{1}.
	\tag*{\qedhere}
\end{equation}
\end{proof}

\begin{lemma}\label{lem:delta-bounded}
	For all $i \in \Sample$, $\abs{\delta_i} \leq 4 \constpo \constpi$ with probability one.
\end{lemma}

\begin{proof}
Distribute the absolute value among the terms:
\begin{equation}
	\abs{\delta_i} \leq \frac{\abs{\ooi}}{\preia} + \frac{\abs{\ooi}}{\preib} + \abs{\poexia} + \abs{\poexib},
\end{equation}
where the facts that $\exiia, \exiib \in \bset$ and that $\preia,\preib \geq 0$ were used.
Conditions~\mainref{cond:bounded-pos} and~\mainref{cond:positivity} then give
\begin{equation}
	\frac{\abs{\ooi}}{\preia} + \frac{\abs{\ooi}}{\preib} + \abs{\poexia} + \abs{\poexib}
	\leq 2 \constpo \constpi + 2 \constpo
	\leq 4 \constpo \constpi.
\end{equation}
\end{proof}

\subsection{Proposition~\ref*{prop:zht-consistency}: Consistency without positivity}

\begin{refproposition}{\ref*{prop:zht-consistency}}
\propzhtconsistency
\end{refproposition}

\begin{proof}
Decompose the mean square error of the estimator into its squared bias and variance:
\begin{equation}
	\E[\Big]{\paren[\big]{\estd - \efd}^2} = \paren[\big]{\E{\estd} - \efd}^2 + \Var[\big]{\estd}.
\end{equation}
Applying the bounds in Lemmas~\ref{lem:non-positivity-bias}~and~\ref{lem:non-positivity-var}, we get
\begin{multline}
	\E[\Big]{\paren[\big]{\estd - \efd}^2}
	\leq \constpo^2 \paren[\big]{\zpravga + \zpravgb}^2 + \frac{4 \constpo^2}{n} \bracket[\big]{\zprmoma + \zprmomb}
	\\
	+ 20 \constpo^2 \braces[\Big]{\bracket[\big]{\zprmoma}^2 \ddepextstda + \bracket[\big]{\zprmomb}^2 \ddepextstdb}
	+ 4 \paren[\big]{\edepa + \edepb}  + 4 \paren[\big]{\udepa + \udepb}.
\end{multline}
A premise of the proposition was $\zprmoma \leq \constarb$ and $\zprmomb \leq \constarb$, which implies
\begin{multline}
	\E[\Big]{\paren[\big]{\estd - \efd}^2}
	\leq \constpo^2 \paren[\big]{\zpravga + \zpravgb}^2 + \frac{8 \constarb \constpo^2}{n}
	+ 20 \constarb^2 \constpo^2 \paren[\big]{\ddepextstda + \ddepextstdb}
	\\
	+ 4 \paren[\big]{\edepa + \edepb}  + 4 \paren[\big]{\udepa + \udepb}.
\end{multline}
Because the square root is a concave function on the positive number line, we can apply Jensen's inequality to get
\begin{multline}
	\sqrt{\E[\Big]{\paren[\big]{\estd - \efd}^2}}
	\leq \constpo \paren[\big]{\zpravga + \zpravgb} + \frac{3 \constarb \constpo}{n}
	+ 5 \constarb \constpo \paren[\big]{\ddepextstda^{0.5} + \ddepextstdb^{0.5}}
	\\
	+ 2 \paren[\big]{\edepa^{0.5} + \edepb^{0.5}}  + 2 \paren[\big]{\udepa^{0.5} + \udepb^{0.5}},
\end{multline}
which implies that the root mean square error is asymptotically bounded as
\begin{equation}
	\sqrt{\E[\Big]{\paren[\big]{\estd - \efd}^2}}
	= \bigO[\big]{n^{-0.5} + \zpravga + \zpravgb + \ddepextstda^{0.5} + \ddepextstdb^{0.5} + \edepa^{0.5} + \edepb^{0.5} + \udepa^{0.5} + \udepb^{0.5}}.
\end{equation}
Thus, the estimator is consistent in mean square under the conditions of the proposition.
Markov's inequality then gives consistency and the rate of convergence in probability.
\end{proof}

\begin{lemma}\label{lem:non-positivity-bias}
Provided that Condition~\mainref{cond:bounded-pos} holds for exposures $\exea$ and $\exeb$, the absolute bias of the Horvitz--Thompson estimator is upper bounded by
\begin{equation}
	\abs[\big]{\E{\estd} - \efd}
	\leq \constpo \paren[\big]{\zpravga + \zpravgb}.
\end{equation}
\end{lemma}

\begin{proof}
Write the estimator as in Lemma~\ref{lem:non-positivity-est}, and take its expectation
\begin{multline}
	\E{\estd}
	= \avgin \frac{\bracket{1 - \zpreia} \E{\exiia \ooi}}{\preia + \zpreia} - \avgin \frac{\bracket{1 - \zpreib} \E{\exiib \ooi}}{\preib + \zpreib}
	\\
	= \avgin \bracket{1 - \zpreia} \poexia - \avgin \bracket{1 - \zpreib} \poexib.
\end{multline}
The bias is therefore
\begin{equation}
	\E{\estd} - \efd
	= - \avgin \zpreia \poexia + \avgin \zpreib \poexib.
\end{equation}
Using Condition~\mainref{cond:bounded-pos}, we can upper bound the absolute value of the bias as
\begin{equation}
	\abs[\big]{\E{\estd} - \efd}
	\leq \avgin \zpreia \abs{\poexia} + \avgin \zpreib \abs{\poexib}
	= \constpo \paren[\big]{\zpravga + \zpravgb}.
\end{equation}
where $\zpravge$ is defined in Section~\ref{sec:lack-positivity}.
\end{proof}

\begin{lemma}\label{lem:non-positivity-est}
When $0/0$ is defined to be zero, the Horvitz--Thompson estimator can be written as
\begin{equation}
	\estd = \avgin \frac{\bracket{1 - \zpreia} \exiia \ooi}{\preia + \zpreia} - \avgin \frac{\bracket{1 - \zpreib} \exiib \ooi}{\preib + \zpreib},
\end{equation}
where $\zpreie = \indicator{\preie = 0}$.
\end{lemma}

\begin{proof}
The estimator is
\begin{equation}
	\estd = \avgin \frac{\exiia \ooi}{\preia} - \avgin \frac{\exiib \ooi}{\preib}.
\end{equation}
Note that $\exiie$ is zero with probability one when $\zpreie = 1$.
Furthermore, $\preie + \zpreie$ is one when $\preie = 0$, and $\preie + \zpreie = \preie$ when $\preie > 0$.
Taken together, it follows that
\begin{equation}
	\frac{\exiia \ooi}{\preia}
	= \frac{\exiia \ooi}{\preie + \zpreie}
\end{equation}
with probability one, because both expressions are zero when $\zpreie = 1$ given that we define $0/0$ to be zero.
The fact that $\exiie$ is zero when $\zpreie = 1$ also implies that
\begin{equation}
	\exiie \ooi
	= \bracket{1 - \zpreie} \exiie \ooi.
	\tag*{\qedhere}
\end{equation}
\end{proof}

\begin{lemma}\label{lem:non-positivity-var}
\begin{multline}
	\Var[\big]{\estd}
	\leq \frac{4 \constpo^2}{n} \bracket[\big]{\zprmoma + \zprmomb} + 20 \constpo^2 \braces[\Big]{\bracket[\big]{\zprmoma}^2 \ddepextstda + \bracket[\big]{\zprmomb}^2 \ddepextstdb}
	\\
	+ 4 \paren[\big]{\edepa + \edepb}  + 4 \paren[\big]{\udepa + \udepb},
\end{multline}
where $\ddepextstde = \ddepext[\big]{\exe}{\exppimom / (\exppimom - 2)}$.
\end{lemma}

\begin{proof}
Use Lemmas~\ref{lem:var-sum-bound}~and~\ref{lem:non-positivity-est} to bound the variance of the estimator as
\begin{equation}
	\Var[\big]{\estd}
	\leq 2 \Var[\bigg]{\avgin \frac{\bracket{1 - \zpreia} \exiia \ooi}{\preia + \zpreia}} + 2 \Var[\bigg]{\avgin \frac{\bracket{1 - \zpreib} \exiib \ooi}{\preib + \zpreib}}.
\end{equation}
Using a similar decomposition as in the proof of Proposition~\mainref{prop:variance-bound}, we can bound each term of this expression as
\begin{align}
	& 2 \Var[\bigg]{\avgin \frac{\bracket{1 - \zpreie} \exiie \ooi}{\preie + \zpreie}}
	\\
	&\qqqquad\leq \frac{4}{n^2} \sumin \bracket[\bigg]{\frac{1 - \zpreie}{\preie + \zpreie}}^2 \Var{\exiie \ooi}
	\\
	&\qqqqquad + \frac{4}{n^2} \sumin \sum_{j \neq i} \bracket[\bigg]{\frac{\bracket{1 - \zpreie} \poexie}{\preie + \zpreie}} \bracket[\bigg]{\frac{\bracket{1 - \zpreje} \poexje}{\preje + \zpreje}} \Cov[\big]{\exiie, \exije}
	\\
	&\qqqqquad + \frac{4}{n^2} \sumin \sum_{j \neq i} \bracket[\bigg]{\frac{1 - \zpreie}{\preie + \zpreie}} \bracket[\bigg]{\frac{1 - \zpreje}{\preje + \zpreje}} \Cov[\big]{\exiie \seri, \exije \serj}.
\end{align}
The terms are bounded, respectively, by Lemmas~\ref{lem:non-positivity-var-first-term},~\ref{lem:non-positivity-var-middle-term}~and~\ref{lem:non-positivity-var-final-term}, resulting in
\begin{equation}
	2 \Var[\bigg]{\avgin \frac{\bracket{1 - \zpreie} \exiie \ooi}{\preie + \zpreie}}
	\leq \frac{4 \constpo^2 \zprmome}{n} + 20 \constpo^2 \bracket[\big]{\zprmome}^2 \ddepextstde + 4 \edepe + 4 \udepe
\end{equation}
The lemma follows when this bound is applied to the two terms in the first expression.
\end{proof}

\begin{lemma}\label{lem:non-positivity-var-first-term}
\begin{equation}
	\frac{4}{n^2} \sumin \bracket[\bigg]{\frac{1 - \zpreie}{\preie + \zpreie}}^2 \Var{\exiie \ooi}
	\leq \frac{4 \constpo^2 \zprmome}{n}.
\end{equation}
\end{lemma}

\begin{proof}
The expression can be bounded as
\begin{multline}
	\frac{4}{n^2} \sumin \bracket[\bigg]{\frac{1 - \zpreie}{\preie + \zpreie}}^2 \Var{\exiie \ooi}
	\leq \frac{4}{n^2} \sumin \bracket[\bigg]{\frac{1 - \zpreie}{\preie + \zpreie}}^2 \E[\big]{\exiie \ooi^2}
	\\
	= \frac{4}{n^2} \sumin \bracket[\bigg]{\frac{1 - \zpreie}{\preie + \zpreie}}^2 \preie \E[\big]{\ooi^2 \given \exoi = \exe}
	\leq \frac{4}{n^2} \sumin \bracket[\bigg]{\frac{1 - \zpreie}{\preie + \zpreie}} \constpo^2,
\end{multline}
where the last inequality follows from Lemma~\ref{lem:outcome-bounds-ysq} and the fact that
\begin{equation}
	\bracket[\bigg]{\frac{1 - \zpreie}{\preie + \zpreie}}^2
	= \begin{cases}
		\bracket{\preie}^{-2} & \text{if } \zpreie = 0,
		\\
		0 & \text{if } \zpreie = 1.
	\end{cases}
\end{equation}
When $\exppimom \geq 1$, we can use Jensen's inequality to get
\begin{equation}
	\frac{1}{n} \sumin \bracket[\bigg]{\frac{1 - \zpreie}{\preie + \zpreie}}
	\leq \bracket[\Bigg]{\frac{1}{n} \sumin \bracket[\bigg]{\frac{1 - \zpreie}{\preie + \zpreie}}^{\exppimom}}^{1 / \exppimom}
	= \zprmome,
\end{equation}
where $\zprmome$ is defined in Section~\ref{sec:lack-positivity}.
\end{proof}

\begin{lemma}\label{lem:non-positivity-var-middle-term}
\begin{equation}
	\frac{4}{n^2} \sumin \sum_{j \neq i} \bracket[\bigg]{\frac{\bracket{1 - \zpreie} \poexie}{\preie + \zpreie}} \bracket[\bigg]{\frac{\bracket{1 - \zpreje} \poexje}{\preje + \zpreje}} \Cov[\big]{\exiie, \exije}
	\leq 4 \constpo^2 \bracket[\big]{\zprmome}^2 \ddepextstde,
\end{equation}
where $\ddepextstde = \ddepext[\big]{\exe}{\exppimom / (\exppimom - 2)}$.
\end{lemma}

\begin{proof}
Using Lemma~\ref{lem:outcome-bounds-po}, the expression can be bounded as
\begin{multline}
	\frac{4}{n^2} \sumin \sum_{j \neq i} \bracket[\bigg]{\frac{\bracket{1 - \zpreie} \poexie}{\preie + \zpreie}} \bracket[\bigg]{\frac{\bracket{1 - \zpreje} \poexje}{\preje + \zpreje}} \Cov[\big]{\exiie, \exije}
	\\
	\leq \frac{4 \constpo^2}{n^2} \sumin \sum_{j \neq i} \bracket[\bigg]{\frac{1 - \zpreie}{\preie + \zpreie}} \bracket[\bigg]{\frac{1 - \zpreje}{\preje + \zpreje}} \abs[\big]{\Cov{\exiie, \exije}}.
\end{multline}
Apply Hölder's inequality with conjugates $\exppimom / 2$ and $\exppimom / (\exppimom - 2)$ to get
\begin{align}
	&\frac{4 \constpo^2}{n^2} \sumin \sum_{j \neq i} \bracket[\bigg]{\frac{1 - \zpreie}{\preie + \zpreie}} \bracket[\bigg]{\frac{1 - \zpreje}{\preje + \zpreje}} \abs[\big]{\Cov{\exiie, \exije}}
	\\
	&\qquad \leq \frac{4 \constpo^2}{n^2} \bracket[\Bigg]{\sumin \sum_{j \neq i} \bracket[\bigg]{\frac{1 - \zpreie}{\preie + \zpreie}}^{\exppimom / 2} \bracket[\bigg]{\frac{1 - \zpreje}{\preje + \zpreje}}^{\exppimom / 2}}^{2 / \exppimom}
	\\
	&\qqqquad \times \bracket[\Bigg]{\sumin \sum_{j \neq i} \abs[\big]{\Cov{\exiie, \exije}}^{\exppimom / (\exppimom - 2)}}^{(\exppimom - 2) / \exppimom}
	\\
	&\qquad \leq 4 \constpo^2 \bracket[\Bigg]{\frac{1}{n^2} \sumin \sum_{j \neq i} \bracket[\bigg]{\frac{1 - \zpreie}{\preie + \zpreie}}^{\exppimom / 2} \bracket[\bigg]{\frac{1 - \zpreje}{\preje + \zpreje}}^{\exppimom / 2}}^{2 / \exppimom}
	\\
	&\qqqquad \times \bracket[\Bigg]{\frac{1}{n^2} \sumin \sum_{j \neq i} \abs[\big]{\Cov{\exiie, \exije}}^{\exppimom / (\exppimom - 2)}}^{(\exppimom - 2) / \exppimom}
\end{align}
The middle factor can be bounded as
\begin{align}
	&\bracket[\Bigg]{\frac{1}{n^2} \sumin \sum_{j \neq i} \bracket[\bigg]{\frac{1 - \zpreie}{\preie + \zpreie}}^{\exppimom / 2} \bracket[\bigg]{\frac{1 - \zpreje}{\preje + \zpreje}}^{\exppimom / 2}}^{2 / \exppimom}
	\\
	&\qquad \leq \bracket[\Bigg]{\frac{1}{n^2} \sumin \sumjn \bracket[\bigg]{\frac{1 - \zpreie}{\preie + \zpreie}}^{\exppimom / 2} \bracket[\bigg]{\frac{1 - \zpreje}{\preje + \zpreje}}^{\exppimom / 2}}^{2 / \exppimom}
	\\
	&\qquad = \bracket[\Bigg]{\paren[\bigg]{\frac{1}{n} \sumin \bracket[\bigg]{\frac{1 - \zpreie}{\preie + \zpreie}}^{\exppimom / 2}}^2}^{2 / \exppimom}.
\end{align}
Apply Jensen's inequality on the square to get
\begin{equation}
	\bracket[\Bigg]{\paren[\bigg]{\frac{1}{n} \sumin \bracket[\bigg]{\frac{1 - \zpreie}{\preie + \zpreie}}^{\exppimom / 2}}^2}^{2 / \exppimom}
	\leq \bracket[\Bigg]{\frac{1}{n} \sumin \bracket[\bigg]{\frac{1 - \zpreie}{\preie + \zpreie}}^{\exppimom}}^{2 / \exppimom}
	= \bracket[\big]{\zprmome}^2.
\end{equation}
The final factor is
\begin{equation}
	\ddepextstde
	= \ddepext[\big]{\exe}{\exppimom / (\exppimom - 2)}
	= \bracket[\Bigg]{\frac{1}{n^2} \sumin \sum_{j \neq i} \abs[\big]{\Cov{\exiie, \exije}}^{\exppimom / (\exppimom - 2)}}^{(\exppimom - 2) / \exppimom},
\end{equation}
so when taken together with the bound on the middle factor, we get
\begin{equation}
	\frac{4}{n^2} \sumin \sum_{j \neq i} \bracket[\bigg]{\frac{\bracket{1 - \zpreie} \poexie}{\preie + \zpreie}} \bracket[\bigg]{\frac{\bracket{1 - \zpreje} \poexje}{\preje + \zpreje}} \Cov[\big]{\exiie, \exije}
	\leq 4 \constpo^2 \bracket[\big]{\zprmome}^2 \ddepextstde.
	\tag*{\qedhere}
\end{equation}
\end{proof}

\begin{lemma}\label{lem:non-positivity-var-final-term}
\begin{multline}
	\frac{4}{n^2} \sumin \sum_{j \neq i} \bracket[\bigg]{\frac{1 - \zpreie}{\preie + \zpreie}} \bracket[\bigg]{\frac{1 - \zpreje}{\preje + \zpreje}} \Cov[\big]{\exiie \seri, \exije \serj}
	\\
	\leq 16 \constpo^2 \bracket[\big]{\zprmome}^2 \ddepextstde + 4 \edepe + 4 \udepe,
\end{multline}
where $\ddepextstde = \ddepext[\big]{\exe}{\exppimom / (\exppimom - 2)}$.
\end{lemma}

\begin{proof}
If either $\zpreie = 1$ or $\zpreje = 1$, then
\begin{equation}
	\bracket[\bigg]{\frac{1 - \zpreie}{\preie + \zpreie}} \bracket[\bigg]{\frac{1 - \zpreje}{\preje + \zpreje}} \Cov[\big]{\exiie \seri, \exije \serj}
	= 0.
\end{equation}
Furthermore, if $\zpreie = 0$ and $\zpreje = 0$, but $\zpreij{\exe, \exe} = \indicator{\preij{\exe, \exe} = 0} = 1$, then
\begin{equation}
	\Cov[\big]{\exiie \seri, \exije \serj}
	= \E[\big]{\exiie \exije \seri \serj} - \E[\big]{\exiie \seri} \E[\big]{\exije \serj}
	= 0,
\end{equation}
because $\exiie \exije$ is then zero with probability one, and $\E{\exiie \seri} = 0$ as shown in the proof of Lemma~\ref{lem:var-prod-decomp}.
It follows that the terms in the double sum are nonzero only when $\zpreij{\exe, \exe} = 0$.
Note that $\zpreij{\exe, \exe} = 0$ implies that $\zpreie = 0$ and $\zpreje = 0$, which means that we can rewrite the expression to be bounded as
\begin{multline}
	\frac{4}{n^2} \sumin \sum_{j \neq i} \bracket[\bigg]{\frac{1 - \zpreie}{\preie + \zpreie}} \bracket[\bigg]{\frac{1 - \zpreje}{\preje + \zpreje}} \Cov[\big]{\exiie \seri, \exije \serj}
	\\
	= \frac{4}{n^2} \sumin \sum_{j \neq i} \bracket[\bigg]{\frac{1 - \zpreij{\exe, \exe}}{\preie \preje + \zpreij{\exe, \exe}}} \Cov[\big]{\exiie \seri, \exije \serj}.
\end{multline}

When $\zpreij{\exe, \exe} = 0$, we can use the same decomposition as in the proof of Lemma~\ref{lem:bound-error-cov} to get
\begin{multline}
	\frac{\Cov[\big]{\exiie \seri, \exije \serj}}{\preie \preje}
	= \frac{\Cov[\big]{\exiie, \exije}}{\preie\preje} \E[\big]{\seri\serj \given \exoi = \exoj = \exe} + \E[\big]{\seri\serj \given \exoi = \exoj = \exe}
	\\
	\leq \frac{4 \constpo^2 \abs[\big]{\Cov{\exiie, \exije}}}{\preie\preje} + \E[\big]{\seri\serj \given \exoi = \exoj = \exe}
\end{multline}
where the final inequality follows from Lemma~\ref{lem:outcome-bounds-errors}.
This means that the expression can be bounded as
\begin{align}
	&\frac{4}{n^2} \sumin \sum_{j \neq i} \bracket[\bigg]{\frac{1 - \zpreij{\exe, \exe}}{\preie \preje + \zpreij{\exe, \exe}}} \Cov[\big]{\exiie \seri, \exije \serj}
	\\
	&\qqquad \leq \frac{16 \constpo^2}{n^2} \sumin \sum_{j \neq i} \bracket[\bigg]{\frac{1 - \zpreij{\exe, \exe}}{\preie \preje + \zpreij{\exe, \exe}}} \abs[\big]{\Cov{\exiie, \exije}}
	\\
	&\qqqqquad + \frac{4}{n^2} \sumin \sum_{j \neq i} \bracket[\big]{1 - \zpreij{\exe, \exe}} \E[\big]{\seri\serj \given \exoi = \exoj = \exe}.
\end{align}
Because $\zpreij{\exe, \exe} = 0$ implies that $\zpreie = 0$ and $\zpreje = 0$, we have
\begin{equation}
	\bracket[\bigg]{\frac{1 - \zpreij{\exe, \exe}}{\preie \preje + \zpreij{\exe, \exe}}}
	\leq \bracket[\bigg]{\frac{1 - \zpreie}{\preie + \zpreie}} \bracket[\bigg]{\frac{1 - \zpreje}{\preje + \zpreje}},
\end{equation}
so we also have
\begin{multline}
	\frac{16 \constpo^2}{n^2} \sumin \sum_{j \neq i} \bracket[\bigg]{\frac{1 - \zpreij{\exe, \exe}}{\preie \preje + \zpreij{\exe, \exe}}} \abs[\big]{\Cov{\exiie, \exije}}
	\\
	\leq \frac{16 \constpo^2}{n^2} \sumin \sum_{j \neq i} \bracket[\bigg]{\frac{1 - \zpreie}{\preie + \zpreie}} \bracket[\bigg]{\frac{1 - \zpreje}{\preje + \zpreje}} \abs[\big]{\Cov{\exiie, \exije}}.
\end{multline}
By applying Hölder's inequality in the same way as in the proof of Lemma~\ref{lem:non-positivity-var-middle-term}, we get
\begin{equation}
	\frac{16 \constpo^2}{n^2} \sumin \sum_{j \neq i} \bracket[\bigg]{\frac{1 - \zpreie}{\preie + \zpreie}} \bracket[\bigg]{\frac{1 - \zpreje}{\preje + \zpreje}} \abs[\big]{\Cov{\exiie, \exije}}
	\leq 16 \constpo^2 \bracket[\big]{\zprmome}^2 \ddepextstde,
\end{equation}

Continuing, we can write
\begin{equation}
	\bracket[\big]{1 - \zpreij{\exe, \exe}} \E[\big]{\seri\serj \given \exoi = \exoj = \exe}
	= \serefije\serefjie + \Cov{\seruij, \seruji \given \exoi = \exoj = \exe},
\end{equation}
because we can use the same derivation as in the proof of Lemma~\ref{lem:bound-error-cov} when $\zpreij{\exe, \exe} = 0$, and we have
\begin{equation}
	\serefije
	= \Cov{\seruij, \seruji \given \exoi = \exoj = \exe}
	= 0,
\end{equation}
when $\zpreij{\exe, \exe} = 1$.
Therefore,
\begin{multline}
	\frac{4}{n^2} \sumin \sum_{j \neq i} \bracket[\big]{1 - \zpreij{\exe, \exe}} \E[\big]{\seri\serj \given \exoi = \exoj = \exe}
	\\
	= \frac{4}{n^2} \sumin \sum_{j \neq i} \serefije\serefjie
	+ \frac{4}{n^2} \sumin \sum_{j \neq i} \Cov{\seruij, \seruji \given \exoi = \exoj = \exe},
\end{multline}
which is equal to $4 \edepe + 4 \udepe$ as defined by Definitions~\mainref{def:error-dependence}.
\end{proof}

\subsection{Propositions~\ref*{prop:ha-consistency},~\ref*{prop:de-unbiasedness},~\ref*{prop:de-consistency}~and~\ref*{prop:gr-consistency}: Other estimators}

\subsubsection{Proposition~\ref*{prop:ha-consistency}: The H{\'a}jek estimator}

The linearization used to prove consistency for the H{\'a}jek estimator requires an alternative representation of the estimator.

\begin{definition}[Components for the H{\'a}jek estimator]
\begin{equation}
	\mutre
	= \sumin \poexie,
	\qqqquad
	\mueste
	= \sumin \frac{\exiie \ooi}{\preie},
	\qqqquad
	\neste
	= \sumin \frac{\exiie}{\preie}.
\end{equation}
\end{definition}

\begin{lemma}\label{lem:bounded-ha-target}
If Condition~\mainref{cond:bounded-pos} holds, then $\mutre = \bigO{n}$ for all $\exe \in \exall$.
\end{lemma}

\begin{proof}
By Lemma~\ref{lem:outcome-bounds-po},
\begin{equation}
	\mutre
	= \sumin \poexie
	\leq \sumin \abs{\poexie}
	\leq \sumin \constpo
	= \constpo n.
	\tag*{\qedhere}
\end{equation}
\end{proof}

\begin{lemma}\label{lem:consistency-ha-est}
If Condition~\mainref{cond:bounded-pos}~and~\mainref{cond:positivity} hold, then
\begin{equation}
\paren{\mueste - \mutre} / n
= \bigOp[\big]{n^{-0.5} + \ddepe^{0.5} + \edepe^{0.5} + \udepe^{0.5}}.
\end{equation}
\end{lemma}

\begin{proof}
Note that $\estd = \paren{\muesta - \muestb} / n$ and $\efd = \paren{\mutra - \mutrb} / n$, so the proofs of Lemma~\ref{lem:unbiasedness} and Proposition~\mainref{prop:variance-bound} can be copied almost in verbatim to show
\begin{equation}
	\E[\big]{\mueste}
	= \mutre
	\qquadand
	\frac{\Var{\mueste}}{n^2}
	\leq \frac{2 \constpo^2 \constpi}{n} + 10 \constpo^2 \constpi^2 \ddepe + 2 \edepe + 2 \udepe.
\end{equation}
The logic of the proof of Proposition~\ref{prop:consistency-decompose} then gives
\begin{equation}
	\sqrt{\E[\big]{\paren{\mueste - \mutre}^2 / n^2}}
	= \bigO[\big]{n^{-0.5} + \ddepe^{0.5} + \edepe^{0.5} + \udepe^{0.5}}.
\end{equation}
Markov's inequality completes the proof.
\end{proof}

\begin{lemma}\label{lem:consistency-ha-n}
If Condition~\mainref{cond:positivity} holds, then $\paren{\neste - n} / n = \bigOp[\big]{n^{-0.5} + \ddepe^{0.5}}$.
\end{lemma}

\begin{proof}
The first step is to show that $\E{\neste} = n$ when $\exe$ satisfies Condition~\mainref{cond:positivity}:
\begin{equation}
	\E{\neste}
	= \sumin \frac{\E{\exiie}}{\preie}
	= \sumin \frac{\preie}{\preie}
	= n.
\end{equation}
Next consider the variance:
\begin{equation}
	\Var{\neste}
	= \sumin \frac{\Var{\exiie}}{\bracket{\preie}^2} + \sumin \sum_{j \neq i} \frac{\Cov{\exiie, \exije}}{\preie \preje}.
\end{equation}
By Condition~\mainref{cond:positivity},
\begin{equation}
	\frac{\Var{\exiie}}{\bracket{\preie}^2}
	= \frac{\preie \bracket{1 - \preie}}{\bracket{\preie}^2} \leq \constpi
	\qquadand
	\frac{\Cov{\exiie, \exije}}{\preie\preje}
	\leq \constpi^2\abs[\big]{\Cov{\exiie, \exije}},
\end{equation}
so $\Var{\neste} / n^2 \leq \constpi n^{-1} + \constpi^2 \ddepe$.
The logic of the proof of Proposition~\ref{prop:consistency-decompose} then gives
\begin{equation}
	\sqrt{\E[\big]{\paren{\neste - n}^2 / n^2}} = \bigO[\big]{n^{-0.5} + \ddepe^{0.5}},
\end{equation}
and Markov's inequality completes the proof.
\end{proof}

\begin{refproposition}{\ref*{prop:ha-consistency}}
\prophaconsistency
\end{refproposition}

\begin{proof}
Note that $\esthad = \muesta / \nesta - \muestb / \nestb$ and $\efd = \mutra / n - \mutrb / n$, so we can write
\begin{equation}
	\esthad - \efd
	= \paren[\bigg]{\frac{\muesta}{\nesta} - \frac{\mutra}{n}} - \paren[\bigg]{\frac{\muestb}{\nestb} - \frac{\mutrb}{n}}.
\end{equation}
For any $\exe \in \setb{\exea, \exeb}$, consider
\begin{equation}
	\frac{\mueste}{\neste} - \frac{\mutre}{n}
	= \frac{\mueste / n}{\neste / n} - \frac{\paren{\mutre / n}\paren{\neste / n}}{\neste / n}
	= \frac{\paren{\mueste - \mutre} / n}{\neste / n} - \frac{\paren{\mutre / n}\paren{\neste - n} / n}{\neste / n},
\end{equation}
where Lemma~\ref{lem:consistency-ha-n} ensures that we can ignore the event $\neste = 0$.

Let $f\paren{x, y} = x / y$ and consider a Taylor expansion of the two terms around $\paren{0, 1}$:
\begin{align}
	\frac{\paren{\mueste - \mutre} / n}{\neste / n}
	&= f\paren[\big]{\paren{\mueste - \mutre} / n, \neste / n}
	= \paren{\mueste - \mutre} / n + r_1
	\\
	\frac{\paren{\mutre / n}\paren{\neste - n} / n}{\neste / n}
	&= f\paren[\big]{\paren{\mutre / n}\paren{\neste - n} / n, \neste / n}
	= \paren{\mutre / n}\paren{\neste - n} / n + r_2
\end{align}
where $r_1 = \littleOp[\big]{\paren{\mueste - \mutre} / n + \paren{\neste - n} / n}$ and $r_2 = \littleOp[\big]{\paren{\mutre / n}\paren{\neste - n} / n +  \paren{\neste - n} / n}$ because Lemmas~\ref{lem:bounded-ha-target},~\ref{lem:consistency-ha-est} and~\ref{lem:consistency-ha-n} give convergence of $\paren{\mueste - \mutre} / n$ and $\paren{\mutre / n}\paren{\neste - n} / n$ to zero and of $\neste / n$ to one.
Lemma~\ref{lem:bounded-ha-target} gives $\paren{\mutre / n}\paren{\neste - n} / n = \bigOp[\big]{\paren{\neste - n} / n}$, so by Lemmas~\ref{lem:consistency-ha-est} and~\ref{lem:consistency-ha-n}
\begin{equation}
	\esthad - \efd
	= \bigOp[\big]{n^{-0.5} + \ddep{\exea}^{0.5} + \ddep{\exeb}^{0.5} + \edep{\exea}^{0.5} + \edep{\exeb}^{0.5} + \udep{\exea}^{0.5} + \udep{\exeb}^{0.5}}.
	\tag*{\qedhere}
\end{equation}
\end{proof}

\subsubsection{Propositions~\ref*{prop:de-unbiasedness}~and~\ref*{prop:de-consistency}: The difference estimator}

\begin{refproposition}{\ref*{prop:de-unbiasedness}}
\propdeunbiasedness
\end{refproposition}

\begin{proof}
Write the estimator as
\begin{align}
	\estded
	&= \avgin \bracket[\big]{\poexestia - \poexestib} + \avgin \frac{\paren{\exiia - \exiib} \bracket[\big]{\ooi - \poexesti{\exoi}}}{\prei{\exoi}}
	\\
	&= \avgin \bracket[\big]{\poexestia - \poexestib} + \avgin \bracket[\bigg]{\frac{\exiia \ooi}{\preia} - \frac{\exiib \ooi}{\preib} - \frac{\exiia \poexestia}{\preia} + \frac{\exiib \poexestib}{\preib}}
	\\
	&= \estd + \avgin \bracket[\big]{\poexestia - \poexestib} - \avgin \bracket[\bigg]{\frac{\exiia \poexestia}{\preia} - \frac{\exiib \poexestib}{\preib}}.
\end{align}
Take expectations to get
\begin{equation}
	\E[\big]{\estded}
	= \efd + \avgin \bracket[\Big]{\E[\big]{\poexestia} - \E[\big]{\poexestib}} - \avgin \bracket[\bigg]{\frac{\E[\big]{\exiia \poexestia}}{\preia} - \frac{\E[\big]{\exiib \poexestib}}{\preib}}.
\end{equation}
We have
\begin{equation}
	\frac{\E[\big]{\exiie \poexestie}}{\preie}
	= \frac{\E[\big]{\exiie} \E[\big]{\poexestie}}{\preie}
	= \frac{\preie \E[\big]{\poexestie}}{\preie}
	= \E[\big]{\poexestie},
\end{equation}
because the predictions are external and Condition~\mainref{cond:positivity}.
As a result, the two last terms in the expectation cancel.
\end{proof}

\begin{refproposition}{\ref*{prop:de-consistency}}
\propdeconsistency
\end{refproposition}

\begin{proof}
Write the estimator as
\begin{equation}
	\estded
	= \estd + \avgin \bracket[\big]{\poexestia - \poexestib} - \avgin \bracket[\bigg]{\frac{\exiia \poexestia}{\preia} - \frac{\exiib \poexestib}{\preib}}.
\end{equation}
Apply Lemma~\ref{lem:var-sum-bound} to get
\begin{multline}
	\Var{\estded}
	\leq 5\Var[\big]{\estd} + \frac{5}{n^2} \Var[\bigg]{\sumin \poexestia} + \frac{5}{n^2} \Var[\bigg]{\sumin \poexestib}
	\\
	+ \frac{5}{n^2} \Var[\bigg]{\sumin \frac{\exiia \poexestia}{\preia}} + \frac{5}{n^2} \Var[\bigg]{\sumin \frac{\exiib\poexestib}{\preib}}.
\end{multline}
The first term is bounded by Proposition~\mainref{prop:variance-bound} as
\begin{equation}
	5\Var[\big]{\estd}
	\leq 40 \constpo^2 \constpi / n + 100 \constpo^2 \constpi^2 \paren[\big]{\ddep{\exea} + \ddep{\exeb}} + 20 \paren[\big]{\edep{\exea} + \edep{\exeb} + \udep{\exea} + \udep{\exeb}}.
\end{equation}
Consider the two subsequent terms:
\begin{equation}
	\frac{5}{n^2} \Var[\bigg]{\sumin \poexestie}
	= \frac{5}{n^2} \sumin \Var[\big]{\poexestie} + \frac{5}{n^2} \sumin \sum_{j \neq i} \Cov[\big]{\poexestie, \poexestje}
	\leq \frac{5 \constpred}{n} + 5 \pdepe,
\end{equation}
where Condition~\ref{cond:prediction-moments} and Definition~\ref{def:average-prediction-dependence} were applied in the last inequality.

Next, consider the last two terms in the variance expression:
\begin{equation}
	\frac{5}{n^2} \Var[\bigg]{\sumin \frac{\exiie\poexestie}{\preie}}
	= \frac{5}{n^2} \sumin \frac{\Var[\big]{\exiie \poexestie}}{\bracket{\preie}^2}
	+ \frac{5}{n^2} \sumin \sum_{j \neq i} \frac{\Cov[\big]{\exiie \poexestie, \exije \poexestje}}{\preie \preje}.
\end{equation}
Consider the first term, and recall that the predictions are external (i.e., independent of the assignemnts), so
\begin{equation}
	\frac{\Var[\big]{\exiie \poexestie}}{\bracket{\preie}^2}
	= \frac{\bracket{\preie}^2 \Var[\big]{\poexestie} + \paren[\big]{\E[\big]{\poexestie}}^2 \Var[\big]{\exiie} + \Var[\big]{\exiie} \Var[\big]{\poexestie}}{\bracket{\preie}^2}.
\end{equation}
Note that $\Var{\exiie} / \preie \leq 1$ because $\exiie$ is binary.
Furthermore, Condition~\ref{cond:prediction-moments} gives
\begin{equation}
	\Var[\big]{\poexestie}
	\leq \E[\Big]{\bracket[\big]{\poexestie}^2}
	\leq \constpred
	\qquadand
	\paren[\big]{\E[\big]{\poexestie}}^2
	\leq \E[\Big]{\bracket[\big]{\poexestie}^2}
	\leq \constpred,
\end{equation}
where the second result uses Jensen's inequality.
Together with Condition~\mainref{cond:positivity}, it follows that
\begin{equation}
	\frac{5}{n^2} \sumin \frac{\Var[\big]{\exiie \poexestie}}{\bracket{\preie}^2}
	\leq \frac{5}{n^2} \sumin \paren[\big]{\constpred + 2 \constpi^2 \constpred}
	\leq \frac{5 \constpred}{n} + \frac{10 \constpi \constpred}{n}.
\end{equation}

Next, consider the second term.
Recall again that the predictions are external and apply the covariance decomposition in \citet{Bohrnstedt1969Exact} to get
\begin{multline}
	\Cov[\big]{\exiie \poexestie, \exije \poexestje}
	= \preie \preje \Cov[\big]{\poexestie, \poexestje} + \E[\big]{\poexestie} \E[\big]{\poexestje} \Cov[\big]{\exiie, \exije}
	\\
	+ \Cov[\big]{\exiie, \exije} \Cov[\big]{\poexestie, \poexestje}
\end{multline}
Note that
\begin{equation}
	0
	\leq \preie \preje + \Cov[\big]{\exiie, \exije}
	= \E[\big]{\exiie \exije}
	\leq 1,
\end{equation}
so
\begin{equation}
	\preie \preje \Cov[\big]{\poexestie, \poexestje} + \Cov[\big]{\exiie, \exije} \Cov[\big]{\poexestie, \poexestje}
	\leq \abs[\big]{\Cov{\poexestie, \poexestje}}.
\end{equation}
By Jensen's inequality and Condition~\ref{cond:prediction-moments},
\begin{equation}
	\abs[\Big]{\E[\big]{\poexestie} \E[\big]{\poexestje}}
	\leq \sqrt{\E[\Big]{\bracket[\big]{\poexestie}^2} \E[\Big]{\bracket[\big]{\poexestje}^2}}
	\leq \sqrt{\constpred^2}
	= \constpred,
\end{equation}
so
\begin{equation}
	\E[\big]{\poexestie} \E[\big]{\poexestje} \Cov[\big]{\exiie, \exije}
	\leq \constpred \abs[\big]{\Cov{\exiie, \exije}}.
\end{equation}
Together with Condition~\mainref{cond:positivity}, this yields
\begin{equation}
	\frac{5}{n^2} \sumin \sum_{j \neq i} \frac{\Cov[\big]{\exiie \poexestie, \exije \poexestje}}{\preie \preje}
	\leq 5 \constpi^2 \pdepe + 5 \constpi^2 \constpred \ddepe.
\end{equation}
Together with the bound on the first term, we have
\begin{equation}
	\frac{5}{n^2} \Var[\bigg]{\sumin \frac{\exiie\poexestie}{\preie}}
	\leq \frac{5 \constpred}{n} + \frac{10 \constpi \constpred}{n} + 5 \constpi^2 \pdepe + 5 \constpi^2 \constpred \ddepe
\end{equation}

It follows that
\begin{multline}
	\Var{\estded}
	\leq \frac{40 \constpo^2 \constpi + 20 \constpi \constpred + 20 \constpred}{n}
	+ \paren[\big]{100 \constpo^2 \constpi^2 + 5 \constpi^2 \constpred} \bracket[\big]{\ddepa + \ddepb}
	\\
	+ 20 \bracket[\big]{\edepa + \edepb + \udepa + \udepb}
	+ \paren[\big]{5 + 5 \constpi^2} \bracket[\big]{\pdepa + \pdepb}
\end{multline}
Together with Proposition~\ref{prop:de-unbiasedness}, this gives convergence in mean square:
\begin{equation}
	\sqrt{\E[\Big]{\paren[\big]{\estded - \estd}^2}}
	= \bigO{n^{-0.5} + \ddepa^{0.5} + \ddepb^{0.5} + \edepa^{0.5} + \edepb^{0.5} + \udepa^{0.5} + \udepb^{0.5} + \pdepa^{0.5} + \pdepb^{0.5}},
\end{equation}
and Markov's inequality completes the proof.
\end{proof}

\subsubsection{Proposition~\ref*{prop:gr-consistency}: The generalized regression estimator}

\begin{refproposition}{\ref*{prop:gr-consistency}}
\propgrconsistency
\end{refproposition}

\begin{proof}
Write the estimator as
\begin{equation}
	\estgrd
	= \estd + \bracket[\Bigg]{\avgin \paren[\bigg]{\covi^\tran - \frac{\exiia \covi^\tran}{\preia}}} \coefesta - \bracket[\Bigg]{\avgin \paren[\bigg]{\covi^\tran - \frac{\exiib \covi^\tran}{\preib}}} \coefestb.
\end{equation}
Using the same approach as in the proofs of Proposition~\ref{prop:consistency-decompose} and Lemma~\ref{lem:unbiasedness}, it can be shown that
\begin{equation}
	\bracket[\Bigg]{\avgin \paren[\bigg]{\covi^\tran - \frac{\exiie\covi^\tran}{\preie}}} = \bigOp[\big]{n^{-0.5} + \ddepe^{0.5}}.
\end{equation}
By Markov's inequality, $\E[\big]{\norm{\coefeste}} = \bigO{1}$ implies $\coefeste = \bigOp{1}$, so
\begin{equation}
	\bracket[\Bigg]{\avgin \paren[\bigg]{\covi^\tran - \frac{\exiie\covi^\tran}{\preie}}} \coefeste = \bigOp[\big]{n^{-0.5} + \ddepe^{0.5}}.
\end{equation}
The proof is then completed by Proposition~\ref{prop:consistency-decompose}.
\end{proof}

\subsection{Proposition~\ref*{prop:bias-variance-estimator}: Bias of variance estimator}

\begin{refproposition}{\ref*{prop:bias-variance-estimator}}
	\propbiasvarianceestimator
\end{refproposition}

\begin{proof}
Recall that the variance estimator is
\begin{multline}
	\EstVarHTAS
	= \frac{1}{n^2} \sumin \sumjn \paren{\exiia - \exiib} \paren{\exija - \exijb} \varwei{\exoi, \exoj} \ooi \ooj
	\\
	+ \frac{1}{n^2} \sumin \sumjn \bracket[\bigg]{\frac{\exiia}{\preia} + \frac{\exiib}{\preib}} \bracket[\big]{\zpreij{\exoi, \exea} + \zpreij{\exoi, \exeb}} \ooi^2,
\end{multline}
where
\begin{equation}
	\varwei{\exeone, \exetwo}
	= \frac{\preijot - \preio \prejt}{\preijot \preio \prejt + \zpreijot},
	\qqquad
	\zpreijot
	= \indicator{\preijot = 0},
\end{equation}
and $\preijot = \Pr{\exoi = \exeone, \exoj = \exetwo}$.
Taking expectations of each terms and applying Lemmas~\ref{lem:bias-var-est-first-term}~and~\ref{lem:bias-var-est-second-term} gives
\begin{align}
	&\E[\Big]{\EstVarHTAS}
	\\
	&\qquad = \E[\Big]{\paren[\big]{\estd}^2} + \varbiasf{3}{\exea, \exeb} + \varbiasf{3}{\exeb, \exea} + 2\varbiasf{4}{\exea, \exeb} - \varbiasf{4}{\exea, \exea} - \varbiasf{4}{\exeb, \exeb}
	\\
	&\qqquad
	+ \frac{1}{n^2} \sumin \sum_{j \neq i} \bracket[\big]{\zpreij{\exea, \exea} + \zpreij{\exea, \exeb}} \bracket{\poexi{\exea}}^2
	- \frac{1}{n^2} \sumin \sum_{j \neq i} \bracket[\big]{1 - \zpreijaa} \poexia \poexja
	\\
	&\qqquad
	+ \frac{1}{n^2} \sumin \sum_{j \neq i} \bracket[\big]{\zpreij{\exeb, \exea} + \zpreij{\exeb, \exeb}} \bracket{\poexi{\exeb}}^2
	- \frac{1}{n^2} \sumin \sum_{j \neq i} \bracket[\big]{1 - \zpreijbb} \poexib \poexjb
	\\
	&\qqquad
	+ \frac{2}{n^2} \sumin \sum_{j \neq i} \bracket[\big]{1 - \zpreijab} \poexia \poexjb
\end{align}
Applying Lemma~\ref{lem:bias-var-est-rest-terms} on the five last terms gives
\begin{multline}
	\E[\Big]{\EstVarHTAS} = \E[\Big]{\paren[\big]{\estd}^2} - \paren[\big]{\efd}^2 + \varbiasf{1}{\exea, \exeb} + \varbiasf{2}{\exea, \exeb} + \varbiasf{2}{\exeb, \exea}
	\\
	+ \varbiasf{3}{\exea, \exeb} + \varbiasf{3}{\exeb, \exea} + 2\varbiasf{4}{\exea, \exeb} - \varbiasf{4}{\exea, \exea} - \varbiasf{4}{\exeb, \exeb}.
\end{multline}
The proof is completed by
\begin{equation}
	\Var[\big]{\estd}
	= \E[\Big]{\paren[\big]{\estd}^2} - \paren[\Big]{\E[\big]{\estd}}^2
	= \E[\Big]{\paren[\big]{\estd}^2} - \paren[\big]{\efd}^2,
\end{equation}
where the final equality follows from Lemma~\ref{lem:unbiasedness}.
\end{proof}

\begin{lemma}\label{lem:bias-var-est-first-term}
\begin{align}
	&\E[\bigg]{\frac{1}{n^2} \sumin \sumjn \paren{\exiia - \exiib} \paren{\exija - \exijb} \varwei{\exoi, \exoj} \ooi \ooj}
	\\
	&\qqquad = \E[\Big]{\paren[\big]{\estd}^2} -
	\frac{1}{n^2} \sumin \paren[\Big]{\E[\big]{\ooi^2 \given \exoi = \exea} + \E[\big]{\ooi^2 \given \exoi = \exeb}}
	\\
	&\qqqquad -
	\frac{1}{n^2} \sumin \sum_{j \neq i} \bracket[\big]{1 - \zpreijaa} \poexia \poexja - \varbiasf{4}{\exea, \exea}
	\\
	&\qqqquad -
	\frac{1}{n^2} \sumin \sum_{j \neq i} \bracket[\big]{1 - \zpreijbb} \poexib \poexjb - \varbiasf{4}{\exeb, \exeb}
	\\
	&\qqqquad +
	\frac{2}{n^2} \sumin \sum_{j \neq i} \bracket[\big]{1 - \zpreijab} \poexia \poexjb + 2\varbiasf{4}{\exea, \exeb}.
\end{align}
\end{lemma}

\begin{proof}
Note that when $\zpreijot = 0$, we have
\begin{equation}
	\varwei{\exeone, \exetwo}
	= \frac{1}{\preio \prejt} - \frac{1}{\preijot},
\end{equation}
and when $\zpreijot = 1$, we have
\begin{equation}
	\varwei{\exeone, \exetwo}
	= - \preio \prejt.
\end{equation}
Hence, for $\zpreijot \in \setb{0, 1}$, we have
\begin{equation}
	\varwei{\exeone, \exetwo}
	= \frac{1 - \zpreijot}{\preio \prejt} - \frac{1 - \zpreijot}{\preijot + \zpreijot} - \zpreijot \preio \prejt.
\end{equation}
This allows us to decompose the first term of the variance estimator as
\begin{align}
	&\frac{1}{n^2} \sumin \sumjn \paren{\exiia - \exiib} \paren{\exija - \exijb} \varwei{\exoi, \exoj} \ooi \ooj
	\\
	&\qqquad =
	\frac{1}{n^2} \sumin \sumjn \frac{\paren{\exiia - \exiib} \paren{\exija - \exijb} \bracket[\big]{1 - \zpreij{\exoi, \exoj}} \ooi \ooj}{\prei{\exoi} \prei{\exoj}}
	\\
	&\qqqquad - \frac{1}{n^2} \sumin \sumjn \frac{\paren{\exiia - \exiib} \paren{\exija - \exijb} \bracket[\big]{1 - \zpreij{\exoi, \exoj}} \ooi \ooj}{\preij{\exoi, \exoj} + \zpreij{\exoi, \exoj}}
	\\
	&\qqqquad - \frac{1}{n^2} \sumin \sumjn \paren{\exiia - \exiib} \paren{\exija - \exijb} \zpreij{\exoi, \exoj} \prei{\exoi} \prei{\exoj} \ooi \ooj
\end{align}
Note that $\zpreij{\exoi, \exoj} = 0$ with probability one by construction.
It follows that
\begin{multline}
	\frac{1}{n^2} \sumin \sumjn \frac{\paren{\exiia - \exiib} \paren{\exija - \exijb} \bracket[\big]{1 - \zpreij{\exoi, \exoj}} \ooi \ooj}{\prei{\exoi} \prei{\exoj}}
	= \paren[\Bigg]{\frac{1}{n} \sumin \frac{\paren{\exiia - \exiib} \ooi}{\prei{\exoi}}}^2
	\\
	= \paren[\Bigg]{\frac{1}{n} \sumin \bracket[\bigg]{\frac{\exiia \ooi}{\preia} - \frac{\exiib \ooi}{\preib}}}^2
	= \paren[\big]{\estd}^2.
\end{multline}
Also because $\zpreij{\exoi, \exoj} = 0$ with probability one, the third term is zero:
\begin{equation}
	\frac{1}{n^2} \sumin \sumjn \paren{\exiia - \exiib} \paren{\exija - \exijb} \zpreij{\exoi, \exoj} \prei{\exoi} \prei{\exoj} \ooi \ooj
	= 0.
\end{equation}

Next, consider the expectation of the second term:
\begin{equation}
	\frac{1}{n^2} \sumin \sumjn \E[\bigg]{\frac{\paren{\exiia - \exiib} \paren{\exija - \exijb} \bracket[\big]{1 - \zpreij{\exoi, \exoj}} \ooi \ooj}{\preij{\exoi, \exoj} + \zpreij{\exoi, \exoj}}}.
\end{equation}
Using the law of total expectation, we can write
\begin{align}
	&\E[\bigg]{\frac{\paren{\exiia - \exiib} \paren{\exija - \exijb} \bracket[\big]{1 - \zpreij{\exoi, \exoj}} \ooi \ooj}{\preij{\exoi, \exoj} + \zpreij{\exoi, \exoj}}}
	\\
	&\qqquad = \bracket[\big]{1 - \zpreijaa} \E[\big]{\ooi \ooj \given \exoi = \exea, \exoj = \exea}
	\\
	&\qqqquad - \bracket[\big]{1 - \zpreijab} \E[\big]{\ooi \ooj \given \exoi = \exea, \exoj = \exeb}
	\\
	&\qqqquad - \bracket[\big]{1 - \zpreijba} \E[\big]{\ooi \ooj \given \exoi = \exeb, \exoj = \exea}
	\\
	&\qqqquad + \bracket[\big]{1 - \zpreijbb} \E[\big]{\ooi \ooj \given \exoi = \exeb, \exoj = \exeb}.
\end{align}
Note that $\zpre{ii}{\exea, \exea} = \zpre{ii}{\exeb, \exeb} = 0$ due to marginal positivity as stipulated by Condition~\mainref{cond:positivity}, and that $\zpre{ii}{\exea, \exeb} = 1$ by the fundamental problem of causal inference.
Hence, when $i = j$, the expression reduces to
\begin{equation}
	\E[\bigg]{\frac{\paren{\exiia - \exiib} \paren{\exija - \exijb} \bracket[\big]{1 - \zpreij{\exoi, \exoj}} \ooi \ooj}{\preij{\exoi, \exoj} + \zpreij{\exoi, \exoj}}}
	= \E[\big]{\ooi^2 \given \exoi = \exea} + \E[\big]{\ooi^2 \given \exoi = \exeb}.
\end{equation}
This means that we can write
\begin{align}
	&\frac{1}{n^2} \sumin \sumjn \E[\bigg]{\frac{\paren{\exiia - \exiib} \paren{\exija - \exijb} \bracket[\big]{1 - \zpreij{\exoi, \exoj}} \ooi \ooj}{\preij{\exoi, \exoj} + \zpreij{\exoi, \exoj}}}
	\\
	&\qqquad =
	\frac{1}{n^2} \sumin \paren[\Big]{\E[\big]{\ooi^2 \given \exoi = \exea} + \E[\big]{\ooi^2 \given \exoi = \exeb}}
	\\
	&\qqqquad +
	\frac{1}{n^2} \sumin \sum_{j \neq i} \bracket[\big]{1 - \zpreijaa} \E[\big]{\ooi \ooj \given \exoi = \exea, \exoj = \exea}
	\\
	&\qqqquad +
	\frac{1}{n^2} \sumin \sum_{j \neq i} \bracket[\big]{1 - \zpreijbb} \E[\big]{\ooi \ooj \given \exoi = \exeb, \exoj = \exeb}
	\\
	&\qqqquad -
	\frac{2}{n^2} \sumin \sum_{j \neq i} \bracket[\big]{1 - \zpreijab} \E[\big]{\ooi \ooj \given \exoi = \exea, \exoj = \exeb}.
\end{align}

Note that $\ooi = \poexi{\exoi} + \sereij + \seruij$.
Therefore, provided that $i \neq j$ and $\zpreijot = 0$, we can write
\begin{align}
	&\E[\big]{\ooi \ooj \given \exoi = \exeone, \exoj = \exetwo}
	\\
	&\qquad =
	\E[\Big]{\paren[\big]{\poexi{\exoi} + \sereij + \seruij} \paren[\big]{\poexj{\exoj} + \sereji + \seruji} \given \exoi = \exeone, \exoj = \exetwo}
	\\
	&\qquad = \poexio \poexjt + \poexjt \serefijot + \poexjt \E{\seruij \given \exoi = \exeone, \exoj = \exetwo}
	\\
	&\qqquad + \poexio \serefjito + \serefijot \serefjito + \serefjito \E{\seruij \given \exoi = \exeone, \exoj = \exetwo}
	\\
	&\qqquad + \poexio \E{\seruji \given \exoi = \exeone, \exoj = \exetwo} + \serefijot \E{\seruji \given \exoi = \exeone, \exoj = \exetwo}
	\\
	&\qqquad + \E{\seruij \seruji \given \exoi = \exeone, \exoj = \exetwo}
\end{align}
Recall that $\seruij = \ooi - \porexij{\exoi, \exoj}$ and $\porexij{\exeone, \exetwo} = \E{\ooi \given \exoi = \exeone, \exoj = \exetwo}$, so
\begin{multline}
	\E{\seruji \given \exoi = \exeone, \exoj = \exetwo}
	= \E[\big]{\ooi - \porexij{\exoi, \exoj} \given \exoi = \exeone, \exoj = \exetwo}
	\\
	= \E{\ooi \given \exoi = \exeone, \exoj = \exetwo} - \porexij{\exeone, \exetwo}
	= 0.
\end{multline}
It follows that
\begin{equation}
	\E{\seruij \seruji \given \exoi = \exeone, \exoj = \exetwo}
	= \Cov{\seruij, \seruji \given \exoi = \exeone, \exoj = \exetwo},
\end{equation}
and
\begin{multline}
	\E[\big]{\ooi \ooj \given \exoi = \exeone, \exoj = \exetwo}
	= \poexio \poexjt + \poexio \serefjito + \poexjt \serefijot
	\\
	+ \serefijot \serefjito + \Cov{\seruij, \seruji \given \exoi = \exeone, \exoj = \exetwo}.
\end{multline}
Recall the definition of $\varbiasf{4}{\exeone, \exetwo}$:
\begin{multline}
	\varbiasf{4}{\exeone, \exetwo} = \frac{1}{n^2} \sumin \sum_{j \neq i} \bracket[\big]{1 - \zpreij{\exeone, \exetwo}} \Big[\poexio\serefji{\exetwo, \exeone} + \poexit\serefij{\exeone, \exetwo}
	\\
	+ \serefij{\exeone, \exetwo}\serefji{\exetwo, \exeone} + \Cov{\seruij, \seruji \given \exoi = \exeone, \exoj = \exetwo}\Big],
\end{multline}
which means we can write
\begin{multline}
	\frac{1}{n^2} \sumin \sum_{j \neq i} \bracket[\big]{1 - \zpreijot} \E[\big]{\ooi \ooj \given \exoi = \exeone, \exoj = \exetwo}
	\\
	= \frac{1}{n^2} \sumin \sum_{j \neq i} \bracket[\big]{1 - \zpreijot} \poexio \poexjt + \varbiasf{4}{\exeone, \exetwo}.
\end{multline}
Taken together, this allows us to write
\begin{align}
	&\frac{1}{n^2} \sumin \sumjn \E[\bigg]{\frac{\paren{\exiia - \exiib} \paren{\exija - \exijb} \bracket[\big]{1 - \zpreij{\exoi, \exoj}} \ooi \ooj}{\preij{\exoi, \exoj} + \zpreij{\exoi, \exoj}}}
	\\
	&\qqquad =
	\frac{1}{n^2} \sumin \paren[\Big]{\E[\big]{\ooi^2 \given \exoi = \exea} + \E[\big]{\ooi^2 \given \exoi = \exeb}}
	\\
	&\qqqquad +
	\frac{1}{n^2} \sumin \sum_{j \neq i} \bracket[\big]{1 - \zpreijaa} \poexia \poexja + \varbiasf{4}{\exea, \exea}
	\\
	&\qqqquad +
	\frac{1}{n^2} \sumin \sum_{j \neq i} \bracket[\big]{1 - \zpreijbb} \poexib \poexjb + \varbiasf{4}{\exeb, \exeb}
	\\
	&\qqqquad -
	\frac{2}{n^2} \sumin \sum_{j \neq i} \bracket[\big]{1 - \zpreijab} \poexia \poexjb - 2\varbiasf{4}{\exea, \exeb}.
	\tag*{\qedhere}
\end{align}
\end{proof}

\begin{lemma}\label{lem:bias-var-est-second-term}
\begin{align}
	&\E[\Bigg]{\frac{1}{n^2} \sumin \sumjn \bracket[\bigg]{\frac{\exiia}{\preia} + \frac{\exiib}{\preib}} \bracket[\big]{\zpreij{\exoi, \exea} + \zpreij{\exoi, \exeb}} \ooi^2}
	\\
	&\qqquad = \frac{1}{n^2} \sumin \paren[\Big]{\E[\big]{\ooi^2 \given \exoi = \exea} + \E[\big]{\ooi^2 \given \exoi = \exeb}}
	\\
	&\qqqquad + \frac{1}{n^2} \sumin \sum_{j \neq i} \bracket[\big]{\zpreij{\exea, \exea} + \zpreij{\exea, \exeb}} \bracket{\poexi{\exea}}^2
	+ \varbiasf{3}{\exea, \exeb}
	\\
	&\qqqquad + \frac{1}{n^2} \sumin \sum_{j \neq i} \bracket[\big]{\zpreij{\exeb, \exea} + \zpreij{\exeb, \exeb}} \bracket{\poexi{\exeb}}^2
	+ \varbiasf{3}{\exeb, \exea}.
\end{align}
\end{lemma}

\begin{proof}
The law of total expectation gives
\begin{align}
	&\E[\Bigg]{\frac{1}{n^2} \sumin \sumjn \bracket[\bigg]{\frac{\exiia}{\preia} + \frac{\exiib}{\preib}} \bracket[\big]{\zpreij{\exoi, \exea} + \zpreij{\exoi, \exeb}} \ooi^2}
	\\
	&\qqquad = \preia \E[\Bigg]{\frac{1}{n^2} \sumin \sumjn \bracket[\bigg]{\frac{\exiia}{\preia} + \frac{\exiib}{\preib}} \bracket[\big]{\zpreij{\exoi, \exea} + \zpreij{\exoi, \exeb}} \ooi^2 \given \exoi = \exea}
	\\
	&\qqqquad + \preib \E[\Bigg]{\frac{1}{n^2} \sumin \sumjn \bracket[\bigg]{\frac{\exiia}{\preia} + \frac{\exiib}{\preib}} \bracket[\big]{\zpreij{\exoi, \exea} + \zpreij{\exoi, \exeb}} \ooi^2 \given \exoi = \exeb}
	\\
	&\qqquad = \frac{1}{n^2} \sumin \sumjn \bracket[\big]{\zpreij{\exea, \exea} + \zpreij{\exea, \exeb}} \E[\big]{\ooi^2 \given \exoi = \exea}
	\\
	&\qqqquad + \frac{1}{n^2} \sumin \sumjn \bracket[\big]{\zpreij{\exeb, \exea} + \zpreij{\exeb, \exeb}} \E[\big]{\ooi^2 \given \exoi = \exeb}
\end{align}
Recall from the proof of Lemma~\ref{lem:bias-var-est-first-term} that $\zpre{ii}{\exea, \exea} = \zpre{ii}{\exeb, \exeb} = 0$ due to marginal positivity as stipulated by Condition~\mainref{cond:positivity}, and that $\zpre{ii}{\exea, \exeb} = 1$ by the fundamental problem of causal inference.
This allows us to write
\begin{align}
	&\E[\Bigg]{\frac{1}{n^2} \sumin \sumjn \bracket[\bigg]{\frac{\exiia}{\preia} + \frac{\exiib}{\preib}} \bracket[\big]{\zpreij{\exoi, \exea} + \zpreij{\exoi, \exeb}} \ooi^2}
	\\
	&\qqquad = \frac{1}{n^2} \sumin \paren[\Big]{\E[\big]{\ooi^2 \given \exoi = \exea} + \E[\big]{\ooi^2 \given \exoi = \exeb}}
	\\
	&\qqqquad + \frac{1}{n^2} \sumin \sum_{j \neq i} \bracket[\big]{\zpreij{\exea, \exea} + \zpreij{\exea, \exeb}} \E[\big]{\ooi^2 \given \exoi = \exea}
	\\
	&\qqqquad + \frac{1}{n^2} \sumin \sum_{j \neq i} \bracket[\big]{\zpreij{\exeb, \exea} + \zpreij{\exeb, \exeb}} \E[\big]{\ooi^2 \given \exoi = \exeb}
\end{align}
Note that $\ooi^2 = \paren[\big]{\poexi{\exoi} + \seri}^2$ and $\E{\seri \given \exoi = \exe} = 0$, so
\begin{equation}
	\E[\big]{\ooi^2 \given \exoi = \exe}
	= \bracket{\poexi{\exe}}^2
	+ 2 \poexi{\exe} \E{\seri \given \exoi = \exe}
	+ \E[\big]{\seri^2 \given \exoi = \exe}
	= \bracket{\poexi{\exe}}^2 + \Var{\seri \given \exoi = \exe},
\end{equation}
and
\begin{multline}
	\frac{1}{n^2} \sumin \sum_{j \neq i} \bracket[\big]{\zpreij{\exea, \exea} + \zpreij{\exea, \exeb}} \E[\big]{\ooi^2 \given \exoi = \exea}
	\\
	= \frac{1}{n^2} \sumin \sum_{j \neq i} \bracket[\big]{\zpreij{\exea, \exea} + \zpreij{\exea, \exeb}} \bracket{\poexi{\exea}}^2
	+ \varbiasf{3}{\exea, \exeb}
\end{multline}
where
\begin{equation}
	\varbiasf{3}{\exeone, \exetwo}
	= \frac{1}{n^2} \sumin \sum_{j \neq i} \bracket[\big]{\zpreij{\exeone, \exeone} + \zpreij{\exeone, \exetwo}} \Var{\seri \given \exoi = \exeone}.
\end{equation}
Taken together, we have
\begin{align}
	&\E[\Bigg]{\frac{1}{n^2} \sumin \sumjn \bracket[\bigg]{\frac{\exiia}{\preia} + \frac{\exiib}{\preib}} \bracket[\big]{\zpreij{\exoi, \exea} + \zpreij{\exoi, \exeb}} \ooi^2}
	\\
	&\qqquad = \frac{1}{n^2} \sumin \paren[\Big]{\E[\big]{\ooi^2 \given \exoi = \exea} + \E[\big]{\ooi^2 \given \exoi = \exeb}}
	\\
	&\qqqquad + \frac{1}{n^2} \sumin \sum_{j \neq i} \bracket[\big]{\zpreij{\exea, \exea} + \zpreij{\exea, \exeb}} \bracket{\poexi{\exea}}^2
	+ \varbiasf{3}{\exea, \exeb}
	\\
	&\qqqquad + \frac{1}{n^2} \sumin \sum_{j \neq i} \bracket[\big]{\zpreij{\exeb, \exea} + \zpreij{\exeb, \exeb}} \bracket{\poexi{\exeb}}^2
	+ \varbiasf{3}{\exeb, \exea}.
	\tag*{\qedhere}
\end{align}
\end{proof}

\begin{lemma}\label{lem:bias-var-est-rest-terms}
\begin{align}
	&\varbiasf{1}{\exea, \exeb} + \varbiasf{2}{\exea, \exeb} + \varbiasf{2}{\exeb, \exea} - \paren[\big]{\efd}^2
	\\
	&\qquad
	= \frac{1}{n^2} \sumin \sum_{j \neq i} \bracket[\big]{\zpreij{\exea, \exea} + \zpreij{\exea, \exeb}} \bracket{\poexi{\exea}}^2
	- \frac{1}{n^2} \sumin \sum_{j \neq i} \bracket[\big]{1 - \zpreijaa} \poexia \poexja
	\\
	&\qqquad
	+ \frac{1}{n^2} \sumin \sum_{j \neq i} \bracket[\big]{\zpreij{\exeb, \exea} + \zpreij{\exeb, \exeb}} \bracket{\poexi{\exeb}}^2
	- \frac{1}{n^2} \sumin \sum_{j \neq i} \bracket[\big]{1 - \zpreijbb} \poexib \poexjb
	\\
	&\qqquad
	+ \frac{2}{n^2} \sumin \sum_{j \neq i} \bracket[\big]{1 - \zpreijab} \poexia \poexjb
\end{align}
\end{lemma}

\begin{proof}
Recall that $\zpreijot = \indicator{\preijot = 0}$, so $\zpreij{\exe, \exe} = \zpreji{\exe, \exe}$.
It follows that
\begin{multline}
\frac{1}{n^2} \sumin \sum_{j \neq i} \zpreij{\exe, \exe} \bracket{\poexi{\exe}}^2
= \frac{1}{2n^2} \sumin \sum_{j \neq i} \bracket[\big]{\zpreij{\exe, \exe} + \zpreji{\exe, \exe}} \bracket{\poexi{\exe}}^2
\\
= \frac{1}{2n^2} \sumin \sum_{j \neq i} \zpreij{\exe, \exe} \bracket[\Big]{\bracket{\poexi{\exe}}^2 + \bracket{\poexj{\exe}}^2},
\end{multline}
where the final equality follows from a reversal of the indices $i$ and $j$ for one of the terms.
Note that
\begin{equation}
	\bracket[\big]{\poexi{\exe} + \poexj{\exe}}^2 = \bracket{\poexi{\exe}}^2 + \bracket{\poexj{\exe}}^2 + 2 \poexi{\exe} \poexj{\exe},
\end{equation}
so we can write
\begin{multline}
\frac{1}{2n^2} \sumin \sum_{j \neq i} \zpreij{\exe, \exe} \bracket[\Big]{\bracket{\poexi{\exe}}^2 + \bracket{\poexj{\exe}}^2}
= \frac{1}{2n^2} \sumin \sum_{j \neq i} \zpreij{\exe, \exe} \bracket[\big]{\poexi{\exe} + \poexj{\exe}}^2
\\
- \frac{1}{n^2} \sumin \sum_{j \neq i} \zpreij{\exe, \exe} \poexi{\exe} \poexj{\exe}.
\end{multline}

Next, note that $\zpreijot = \zprejito$ for the same reason as above.
This allows us to write
\begin{equation}
	\frac{1}{n^2} \sumin \sum_{j \neq i} \zpreijba \bracket{\poexib}^2
	= \frac{1}{n^2} \sumin \sum_{j \neq i} \zpreijab \bracket{\poexjb}^2,
\end{equation}
which means that
\begin{multline}
\frac{1}{n^2} \sumin \sum_{j \neq i} \zpreijab \bracket{\poexia}^2
+ \frac{1}{n^2} \sumin \sum_{j \neq i} \zpreijba \bracket{\poexib}^2
\\
= \frac{1}{n^2} \sumin \sum_{j \neq i} \zpreijab \bracket[\Big]{\bracket{\poexia}^2 + \bracket{\poexjb}^2}.
\end{multline}
Note that
\begin{equation}
	\bracket[\big]{\poexia - \poexjb}^2
	= \bracket{\poexia}^2 + \bracket{\poexjb}^2 - 2 \poexia \poexjb,
\end{equation}
so
\begin{multline}
\frac{1}{n^2} \sumin \sum_{j \neq i} \zpreijab \bracket[\Big]{\bracket{\poexia}^2 + \bracket{\poexjb}^2}.
\\
= \frac{1}{n^2} \sumin \sum_{j \neq i} \zpreijab \bracket[\big]{\poexia - \poexjb}^2
+ \frac{2}{n^2} \sumin \sum_{j \neq i} \zpreijab \poexia \poexjb.
\end{multline}
Taken together, we can write
\begin{align}
	&\frac{1}{n^2} \sumin \sum_{j \neq i} \bracket[\big]{\zpreij{\exea, \exea} + \zpreij{\exea, \exeb}} \bracket{\poexi{\exea}}^2
	+ \frac{1}{n^2} \sumin \sum_{j \neq i} \bracket[\big]{\zpreij{\exeb, \exea} + \zpreij{\exeb, \exeb}} \bracket{\poexi{\exeb}}^2
	\\
	&\qqquad
	=
	\frac{1}{2n^2} \sumin \sum_{j \neq i} \zpreijaa \bracket[\big]{\poexia + \poexja}^2
	- \frac{1}{n^2} \sumin \sum_{j \neq i} \zpreijaa \poexia \poexja
	\\
	&\qqqquad
	+ \frac{1}{2n^2} \sumin \sum_{j \neq i} \zpreijbb \bracket[\big]{\poexib + \poexjb}^2
	- \frac{1}{n^2} \sumin \sum_{j \neq i} \zpreijbb \poexib \poexjb
	\\
	&\qqqquad
	+ \frac{1}{n^2} \sumin \sum_{j \neq i} \zpreijab \bracket[\big]{\poexia - \poexjb}^2
	+ \frac{2}{n^2} \sumin \sum_{j \neq i} \zpreijab \poexia \poexjb.
\end{align}

Recall that
\begin{equation}
	\varbiasf{2}{\exeone, \exetwo}
	= \frac{1}{2n^2} \sumin \sum_{j \neq i} \paren[\Big]{\zpreij{\exeone, \exeone} \bracket[\big]{\poexio + \poexjo}^2 + \zpreij{\exeone, \exetwo} \bracket[\big]{\poexio - \poexit}^2},
\end{equation}
which allows us to write
\begin{align}
	&\frac{1}{n^2} \sumin \sum_{j \neq i} \bracket[\big]{\zpreij{\exea, \exea} + \zpreij{\exea, \exeb}} \bracket{\poexi{\exea}}^2
	+ \frac{1}{n^2} \sumin \sum_{j \neq i} \bracket[\big]{\zpreij{\exeb, \exea} + \zpreij{\exeb, \exeb}} \bracket{\poexi{\exeb}}^2
	\\
	&\qqquad
	= \varbiasf{2}{\exea, \exeb} + \varbiasf{2}{\exeb, \exea}
	+ \frac{2}{n^2} \sumin \sum_{j \neq i} \zpreijab \poexia \poexjb
	\\
	&\qqqquad
	- \frac{1}{n^2} \sumin \sum_{j \neq i} \zpreijaa \poexia \poexja
	- \frac{1}{n^2} \sumin \sum_{j \neq i} \zpreijbb \poexib \poexjb.
\end{align}
In turn, this gives
\begin{align}
	&\frac{1}{n^2} \sumin \sum_{j \neq i} \bracket[\big]{\zpreij{\exea, \exea} + \zpreij{\exea, \exeb}} \bracket{\poexi{\exea}}^2
	- \frac{1}{n^2} \sumin \sum_{j \neq i} \bracket[\big]{1 - \zpreijaa} \poexia \poexja
	\\
	&\quad
	+ \frac{1}{n^2} \sumin \sum_{j \neq i} \bracket[\big]{\zpreij{\exeb, \exea} + \zpreij{\exeb, \exeb}} \bracket{\poexi{\exeb}}^2
	- \frac{1}{n^2} \sumin \sum_{j \neq i} \bracket[\big]{1 - \zpreijbb} \poexib \poexjb
	\\
	&\quad
	+ \frac{2}{n^2} \sumin \sum_{j \neq i} \bracket[\big]{1 - \zpreijab} \poexia \poexjb
	\\
	&\qqquad
	= \varbiasf{2}{\exea, \exeb} + \varbiasf{2}{\exeb, \exea}
	+ \frac{2}{n^2} \sumin \sum_{j \neq i} \zpreijab \poexia \poexjb
	\\
	&\qqqquad
	- \frac{1}{n^2} \sumin \sum_{j \neq i} \zpreijaa \poexia \poexja
	- \frac{1}{n^2} \sumin \sum_{j \neq i} \zpreijbb \poexib \poexjb
	\\
	&\qqqquad
	- \frac{1}{n^2} \sumin \sum_{j \neq i} \bracket[\big]{1 - \zpreijaa} \poexia \poexja
	- \frac{1}{n^2} \sumin \sum_{j \neq i} \bracket[\big]{1 - \zpreijbb} \poexib \poexjb
	\\
	&\qqqquad
	+ \frac{2}{n^2} \sumin \sum_{j \neq i} \bracket[\big]{1 - \zpreijab} \poexia \poexjb
	\\
	&\qqquad
	= \varbiasf{2}{\exea, \exeb} + \varbiasf{2}{\exeb, \exea}
	- \frac{1}{n^2} \sumin \sum_{j \neq i} \bracket[\big]{\poexia \poexja + \poexib \poexjb - 2 \poexia \poexjb}
\end{align}

Note that
\begin{equation}
	\varbiasf{1}{\exeone, \exetwo}
	= \frac{1}{n^2} \sumin \bracket[\big]{\poexio - \poexit}^2
	= \frac{1}{n^2} \sumin \bracket[\big]{\poexio\poexio - 2 \poexio\poexit + \poexit\poexit},
\end{equation}
so
\begin{multline}
	- \frac{1}{n^2} \sumin \sum_{j \neq i} \bracket[\big]{\poexia \poexja + \poexib \poexjb - 2 \poexia \poexjb}
	\\
	= \varbiasf{1}{\exea, \exeb} - \frac{1}{n^2} \sumin \sumjn \bracket[\big]{\poexia \poexja + \poexib \poexjb - 2 \poexia \poexjb}.
\end{multline}

Finally, note that
\begin{multline}
	\frac{1}{n^2} \sumin \sumjn \bracket[\big]{\poexia \poexja + \poexib \poexjb - 2 \poexia \poexjb}
	\\
	= \frac{1}{n^2} \sumin \sumjn \bracket[\big]{\poexia \poexja + \poexib \poexjb - \poexia \poexjb - \poexib \poexja}
\end{multline}
and
\begin{equation}
	\poexia \poexja + \poexib \poexjb - \poexia \poexjb - \poexib \poexja
	= \bracket[\big]{\poexia - \poexib} \bracket[\big]{\poexja - \poexjb},
\end{equation}
which implies that
\begin{multline}
	\frac{1}{n^2} \sumin \sumjn \bracket[\big]{\poexia \poexja + \poexib \poexjb - 2 \poexia \poexjb}
	\\
	= \frac{1}{n^2} \sumin \sumjn \bracket[\big]{\poexia - \poexib} \bracket[\big]{\poexja - \poexjb}
	= \paren[\bigg]{\frac{1}{n} \sumin \bracket[\big]{\poexia - \poexib}}^2
	= \paren[\big]{\efd}^2,
\end{multline}
which completes the proof.
\end{proof}

\end{document}